\documentclass[a4paper,10pt]{article}
\usepackage{geometry}
\usepackage{stmaryrd}
\usepackage{bbm}
\usepackage{amsfonts}
\usepackage{amsmath}
\usepackage{amssymb}
\usepackage{mathrsfs}
\usepackage{accents}
\usepackage{amscd}

\newtheorem{atheorem}{\bf \temp}[section]
\newtheorem{thm}[atheorem]{Theorem}
\newtheorem{cor}[atheorem]{Corollary}
\newtheorem{lem}[atheorem]{Lemma}
\newtheorem{prop}[atheorem]{Proposition}
\newtheorem{de}[atheorem]{Definition}
\newtheorem{rem}[atheorem]{Remark}

\numberwithin{equation}{section}

\title{\textbf{Landau Damping in a weakly collisional regime}}
\author{ Xixia Ma \footnote{Yau Mathematical Sciences Center, Tsinghua University. E-mail addresses:kfmaxixia@tsinghua.edu.cn} \ \ \ \ \ \ \
 \\}
\date{}

\begin{document}

\maketitle
\textbf{Abstract.}   In this paper, we consider the nonlinear Vlasov-Poisson equations in a weakly collisional regime
and  study the linear Boltzmann collision operator. We prove that Landau damping still occurs in this case.

%\textbf{Keywords.}
\begin{center}
%\item\section{Linear}
{\bf\large 0. \quad Introduction }
\end{center}

In this paper, it is assumed that the plasma system is weakly collisional, nonrelativistic, hot. The kinetic theory is an effective method of studying the hot plasma particles. Perhaps the most widely used formulation of kinetic theory is the Boltzmann equation, for which the nonrelativistic form  for particles of the $s $ species is
\begin{align}
\partial_{t}f_{s}+v\cdot\nabla_{x}f_{s}+\frac{q_{s}}{m_{s}}(E+v\times B)\cdot\nabla_{v}f_{s}=\frac{df_{s}}{dt}|_{\textmd{collisions}}.
\end{align}

In Eq.(0.1), $f_{s},E,$ and $B$ may be thought of as the $s-$particle species distribution function $f_{s}(x,v,t),$ and the electric and magnetic fields in the plasma averaged over a spatial volume that contains many particles. It was A.A.Vlasov (1945) who first pointed out that Eq.(0.1) is dominated by the term on its left-hand side for a hot plasma. And for much of the study of waves in a hot plasma, it suffices to use the set of Vlasov equations in many situations,
\begin{align}
\partial_{t}f_{s}+v\cdot\nabla_{x}f_{s}+\frac{q_{s}}{m_{s}}(E+v\times B)\cdot\nabla_{v}f_{s}=0.
\end{align}
 For Eq.(0.2), it is well known that Mouhot and Villani[28] made a  ground-breaking work when $B\equiv0$.  And recent we [26] prove Landau damping on Eq.(0.2) in a uniformly  magnetic field case. In this paper, we consider the unmagnetized plasma in the weakly collisional case, that is,
 \begin{align}
\partial_{t}f_{s}+v\cdot\nabla_{x}f_{s}+\frac{q_{s}}{m_{s}}E\cdot\nabla_{v}f_{s}=\nu\frac{df_{s}}{dt}|_{\textmd{collisions}}.
\end{align}
where $\nu\in (0,\nu_{0}],\nu_{0}>0$  some small constant.

First we start with  the linearized Vlasov equation in an unmagnetized plasma to  analyze the effect on  between Landau damping and collision.
We write the linearized Vlasov equation  in collisionless case as
\begin{align}
\left\{
\begin{array}{l}
\partial_{t}f+v\cdot\nabla_{x}f+\frac{q}{m}E\cdot\nabla_{v}f^{0}=0,\\
f(x,v,t_{0})=f_{0}(x,v).
\end{array}\right.
\end{align}
We can solve Eq.(0.4) with a simple integration,
\begin{align}
f(x,v,t)=f_{0}(x-vt,v)-\int^{t}_{t_{0}}\frac{q}{m}E(x-v(t-t'),t')\cdot\nabla_{v}f^{0}dt'.
\end{align}
In order to simplify  the analysis, we assume that $E(x,t)$ is known and we represent it as $\textmd{Re}E_{1}\exp(ikx-i\omega t).$ Then  the integration in Eq.(0.5) becomes
\begin{align}
&f(x,v,t)=f_{0}(x-vt,v)-\textmd{Re}\frac{iqE_{1}}{m}\frac{df^{0}(v)}{dv}\exp(ikx-i\omega t)\frac{1-e^{i(\omega-kv)(t-t_{0})}}{\omega-kv},\notag\\
&f(k,v,\omega)=-\textmd{Re}\frac{iqE_{1\omega}}{m}\frac{df^{0}(v)}{dv}\frac{1}{\omega-kv}
\end{align}
 Then from the first equality of Eq.(0.6), we can observe an important feature  that $f(x,v,t)$ remains finite even at exact wave-particle resonance, $\omega-kv=0.$ On the other hand, the amplitude of $f(x,v,t)$ at resonance grows linearly with $t-t_{0}$ and for $t-t_{0}$ large, $f(x,v,t)$ becomes strongly oscillatory near resonance and desplays a large peak exactly at resonance. However, when we consider the  collisions among particles, we have to limit the magnitude of  $t-t_{0}.$  Now we use a really simple-minded model to simulate the collisional effect on $f(x,v,t) $ as follows:
 \begin{align}
\left\{
\begin{array}{l}
\partial_{t}f+v\cdot\nabla_{x}f+\frac{q}{m}E\cdot\nabla_{v}f^{0}=-\nu f,\\
f(x,v,t_{0})=f_{0}(x,v).
\end{array}\right.
\end{align}
 Then we have \begin{align}
f(x,v,t)=f_{0}(x-vt,v)e^{-\nu(t-t_{0})}-\nu\textmd{Re}\frac{iqE_{1}}{m}\frac{df^{0}(v)}{dv}e^{ikx-i\omega t}e^{-\nu(t-t_{0})}\frac{1-e^{i(\omega-kv)(t-t_{0})}}{\omega-kv}.
\end{align}
For Eq.(0.8), we can regard $e^{-\nu(t-t_{0})}$ as the probability that any single particle, now at $x,v,t,$ suffered a collision at $t_{0}$ in the past, here
$\nu$ is the collision frequency.
Then setting $s=t-t_{0},$ and averaging over the collision times for all particles that have reached $x,v,t,$ we obtain
\begin{align}
&\langle f(x,v,t)\rangle=-\nu\textmd{Re}\frac{iqE_{1}}{m}\frac{df^{0}(v)}{dv}e^{ikx-i\omega t}\int^{\infty}_{0}e^{-\nu s}\frac{1-e^{i(\omega-kv)s}}{\omega-kv}ds\notag\\
&=-\textmd{Re}\frac{iqE_{1}}{m}\frac{df^{0}(v)}{dv}e^{ikx-i\omega t}\frac{1}{\omega-kv+i\nu}.
\end{align}

From the second line of  Eq.(0.9), first  of all,  we shall observe that the effect of collisions in this model has been to transform the appearance of
$\omega,$ namely, $\omega\rightarrow\omega+i\nu.$  The second line of Eq.(0.9) also shows that $f(x,v,t)$ is the product of two peaking functions, one depending
on $\frac{df^{0}(v)}{dv}$ and the other on the resonance denominator, $\omega-kv+i\nu.$ And now we show that a mathematical representation that leads directly to Landau damping is to write the Fourier amplitude for $ f(x,v,t)$ as $\nu\rightarrow0,$
\begin{align}
&f(\omega,k,v)=\lim_{\nu\rightarrow 0^{+}}\frac{iqE(\omega,k)}{m}\frac{df^{0}(v)}{dv}\frac{1}{\omega-kv+i\nu}\notag\\
&=-\frac{iqE(\omega,k)}{m}\frac{df^{0}(v)}{dv}\bigg[P\bigg(\frac{1}{\omega-kv}\bigg)-\frac{i\pi}{|k|}\delta(v-\frac{\omega}{k})\bigg].
\end{align}
Although the $ f$ peak is infinitely sharp in Eq.(0.10), the moments of Eqs.(0.6) and (0.10) will be approximately the same provided that $\frac{df^{0}(v)}{dv}$ in Eq.(0.6) does not change appreciably over the range of $ v$ through which $(\omega-kv+i\nu)^{-1}$ is large. That is, the
collisional model, Eq.(0.6), will lead to the same moments of $f(x,v,t)$ as the collisionfree model Eq.(0.10), provided that the resonance denominator in Eq.(0.6) supplies the dominant peaking effect.

In this paper, based on Mouhot and Villani' work in [28],  we consider the following model,
\begin{equation}
\left\{\begin{array}{l}
\partial_{t}f+v\cdot\nabla_{x}f+\frac{q}{m}E\cdot\nabla_{v}f=\nu\mathcal{C}(f),\\
E=\nabla W(x)\ast\rho(t,x),\\
\rho(t,x)=\int_{\mathbb{R}^{3}}f(t,x,v)dv,\\
f(0,x,v)=f_{0}(x,v),f^{0}=f^{0}(v),
\end{array}\right.
 \end{equation}
 where $\nu\in (0,\nu_{0}],\nu_{0}$   some small constant, and  here $\mathcal{C}(f)$ represents the linear Boltzmann collisional case, namely
 $$\mathcal{C}(f)=\rho f^{0}-f.$$

 We recall the related results on Landau damping on  weakly collisional plasma as follows. First, if $\nu=0,$ some earlier results of the  linearized Vlasov-Poisson
 equation were obtained in [8] by Caglioti and Maffei and in [21] by Hwang and Vel$\acute{\textmd{a}}$zquez. We also refer to the work of Mouhot and Villani [28], they prove Landau damping (linear and nonlinear) in analytic or Gevrey regularity. Later Bedrossian, Masmoudi and Mouhot [6] give a simplified proof in
Gevrey norm. Let us mention that lots of  literature is  devoted to the study of the Vlasov-Poisson-Boltzmann equation with a general Boltzmann collision operator, for example, the paper of Dolbeault and Desvillettes[13]that deals with the large time behavior of solutions and two papers of Guo on the Vlasov-Poisson-Bolzmann equation [18,19] which is about the large time behaviour of solutions: the first one is in a near-vacuum regime, the second is in a near Maxwellian setting. Meanwhile, there are other many references which are concerned with the lager time behavior of solutions,such as Duan and Strain[15], Duan, and Liu [14] and so on.  However, as far as the case of the linear Boltzmann equation, the literature on the stability  is  very scarce, even in a weakly collisional regime, namely, if $\nu\rightarrow0,$ the first paper on Landau damping in this case  is by I.Tristani [31] for the linearized Vlasov-Poisson equation. It should be relevant to compare this kind of question with the one studied by Bedrossian, Masmoudi and Vivol[3] about the two-dimensional Euler equation where the equivalent of $\nu$ should be viscosity. I.Tristani [31], Bedrossian[1],Bendrossian and Wang[7] also study this kind of problematic of uniform analysis of large time behaviour in a weakly collisional regime through another different form : the Vlasov-Fokker-Plank model.

In this paper, we will consider the  nonlinear Vlasov-Poisson equation on weakly collisional plasma. On  one hand, different from the linear case, we have to face the resonance that the nonlinear term brings. In order to deal with  this difficulty, the method  we use is based on the one in [28]. On the other hand, comparing with the collisionless case in which  the index of the decay rate becomes smaller, we find that for the linear case in weakly collisional case,  the index of the decay rate is the same with the initial data, that is due to the effect of the weak collision. However, for the nonlinear case, the decay rate still becomes slow comparing with the initial time because of  the  nonlinear term. In other words, weak collision does not change  the dynamical behavior of the plasma.

 This paper is organized as follows: Section 1 mainly introduces  hybrid analytic norms and the related properties. Section 2 we will prove Landau damping at the linear level.
 We will state  sketch the proof of main theorem at the nonlinear level in section 3. And section 4 will show the deflection estimates of the particles
 trajectory, section 5 is the key section, it will
 state the phenomena of plasma echo. We will control the error terms in section  6, and give the iteration in section 7.

 Now we state our main result as follows.
 \begin{thm} Let $f^{0}:\mathbb{R}^{3}\rightarrow\mathbb{R}_{+}$ be an analytic velocity profile,
 and let $W(x):\mathbb{T}^{3}\rightarrow\mathbb{R}$ satisfy $$ \hat{W}(0)=0,\quad
  |\hat{W}(k)|\leq\frac{1}{1+|k|^{\gamma}},\gamma>1.$$ Further assume that, for some constant $\lambda_{0}>0,$
 \begin{align}
 \sup_{\eta\in\mathbb{R}^{3}}e^{2\pi\lambda_{0}|\eta|}|\tilde{f}^{0}(\eta)|\leq C_{0},\quad \sum_{n\in\mathbb{N}_{0}^{3}}\frac{\lambda_{0}^{n}}{n!}\|\nabla^{n}_{v}f^{0}\|_{L_{dv}^{1}}\leq C_{0}<\infty.
 \end{align}

Consider equations(0.11),  there is $\varepsilon=\varepsilon(\lambda_{0},\mu_{0},\beta,\gamma,\lambda'_{0},\mu'_{0})$ with the following property:
if $f_{0}=f_{0}(x,v)$ is an initial data satisfying
\begin{align}
\sup_{k\in\mathbb{Z}^{3},\eta\in\mathbb{R}^{3}}e^{2\pi\lambda_{0}|\eta|}e^{2\pi\mu_{0}|k|}|f^{0}-f_{0}|+\int_{\mathbb{T}^{3}}\int_{\mathbb{R}^{3}}
|f^{0}-f_{0}|e^{\beta|v|}dvdx\leq\varepsilon,
\end{align}
where any $\beta>0.$

At the same time,we also assume that the following stability condition holds:
%$\quad\mathbf{ how} \quad\mathbf{ to}\quad \mathbf{make} \quad \mathbf{assumptions }$
%$\quad \mathbf{on}\quad f^{0}$ is suitable
\begin{itemize}
       %\item[(i)]$f^{0}(v)$ is  Maxwellian, that is, $f^{0}(v)=f^{0}(v_{\perp},v_{z})=C_{M}exp(-\alpha(v^{2}_{\perp}+v^{2}_{z}))$ for some $\alpha>0$  ;
         %\item[(ii)] $|\hat{f}^{0}(\eta)|\leq C^{0}e^{-\lambda_{1}|\eta_{3}|}e^{-\lambda_{1}|\eta_{1}|}e^{-\lambda_{1}|\eta_{2}|},|\partial_{\eta_{3}}\hat{f}^{0}(\eta)|\leq C^{1}(\lambda_{1})e^{-\lambda_{1}|\eta_{1}|}e^{-\lambda_{1}|\eta_{2}|}e^{-\lambda_{1}|\eta_{3}|},$
            % $|f^{0}(\eta_{1},\eta_{2},v_{3})|\leq Ce^{-\lambda_{1}|\eta_{1}|}e^{-\lambda_{1}|\eta_{2}|}e^{-\alpha_{1}|v_{3}|},$
             % for some constant $\lambda_{1},\alpha_{1}>0;$

        % \item[(iii)]$|\hat{f}_{0}(k,\eta)|\leq C_{0} e^{-\lambda_{0}|\eta_{1}|}e^{-\lambda_{0}|\eta_{2}|}e^{-\lambda_{0}|\eta_{3}|};$

             \item[] $\mathbf{Stability}$ $ \mathbf{condition}:$ for any  velocity $v\in\mathbb{R}^{3},$
             there exists some positive constant $v_{Te}$ such that if $v=\frac{\omega}{k}+i\frac{\nu}{2\pi k},\omega,k$ are frequencies of time and space $t,x,$
             respectively,  then $|v|\gg v_{Te}.$
                  %\item[(v)] $f^{0}(\eta_{1},\eta_{2},v_{3})\leq Ce^{-\lambda_{1}|\eta_{1}|}e^{-\lambda_{1}|\eta_{2}|}e^{-\alpha_{1}|v_{3}|},$ otherwise $f^{0}(\eta_{1},\eta_{2},v_{3})<\varepsilon(v_{3}),$ for some $0<\varepsilon(v_{3})<1$ sufficient small.%$\inf_{k\in\mathbb{Z}^{3}}|\tilde{\mathcal{L}}(\omega,k)-1|>\kappa,$ for some $0<\kappa<1.$
         \end{itemize}
         Then %if $f_{0}(x,v),f^{0}(v)$ are axisymetric, then the solution $f$ of (1.1) is also axisymmetric.
for any fixed $\eta,k,\forall r\in\mathbb{N},$ as $\nu\rightarrow0,|t|\rightarrow\infty,$ we have
         \begin{align}
         &|\hat{f}(t,k,\eta)-\hat{f}_{0}(k,\eta)|\leq  e^{-\lambda'_{0}|\eta+kt|},\quad\|\rho(t,\cdot)-\rho_{0}\|_{C^{r}(\mathbb{T}^{3})}
         =O(e^{-2\pi\lambda'_{0}|t|}),\notag\\
          &\|E(t,\cdot)\|_{C^{r}(\mathbb{T}^{3})}=O(e^{-2\pi\lambda'_{0}|t|}),\notag\\
         \end{align}
         where $\rho_{0}=\int_{\mathbb{T}^{3}}\int_{\mathbb{R}^{3}}f_{0}(x,v)dvdx,$ for any $ 0<\lambda'_{0}<\lambda_{0} .$
         %(这里可能不需要假设初值为轴对称在速度方向 %maybe not assume $f_{0}$ axisymmetric in $v$)
         \end{thm}

         %Here we sketch the difficulties  and methods that the proof of Theorem 0.3 may meet. From Faraday Law of
         % Electromagnetic induction,  the force $F$ generated by the magnetic  field is that $F=BLV,$ where $L$ is a constant and $V$ is the velocity relative to the
         % magnetic field. Since $B(t,x)$ is independent of $V$ and $V$ is a unbounded variable, we almost have no hope to estimate $F,$ but this is the key point to
          %prove cyclotron damping. The idea to solve this problem  is from a new observation on the basis of Lenz's Law, then we reduce the inhomogeneous dynamical
         % equations of particles trajectory to the homogeneous equations  that is similar to that of Landau damping[]. Another difficulty is that the equation that
         % the density
          %$\rho[f]$ function satisfies  does not form the closed equation because of the additional term $v\times B\cdot\nabla_{v} f.$ For this, we have to estimate the
          %equations of the distribution equation and the density equation that bring  different kinds of resonances in different norms that is much more complicated
         % and difficult
          %than Landau damping in[].

          \begin{rem} $\gamma>1$ of Theorem 0.1  can be extended to $\gamma\geq1,$ the difference between $\gamma>1$ and $\gamma=1$ is the proof of the growth
          integral in section 7. The proof of $\gamma=1$ is similar to section 7 in [28], here we omit this case.
          \end{rem}
         \begin{rem}  First,  from the physics viewpoint,  for the collisionless case, the condition that the damping occur is that the number of particles that the wave velocity greatly
         exceeds  their velocity is much larger than the number of particles whose velocity is slower than the wave velocity. However, when considering the collision among particles, from the above theorem, we know that if very little energy due to collision loss, then when the stability condition of the collisionless case is satisfied, the damping still occurs.
         From the dynamical behavior viewpoint,  it can also be understood  that when the collision is very weak, the electric field play main role on the change of  the plasma' trajectories.
         \end{rem}
         \begin{rem} During the proof, it is easily observed that the regularity loss become smaller because of the weak collision, because the collision term provides a regularity $e^{-\nu t}.$
         \end{rem}
         \begin{rem} Combining the results in this paper  and the idea in our previous paper on Cyclotron damping in a uniform bounded magnetic field, Landau damping on Vlasov-Maxwell equations  in a weakly collisional regime may be proved.
         \end{rem}

           \begin{center}
\item\section{Notation and Hybrid analytic norm}
%{\bf\large 1. \quad Introduction }
\end{center}
         %$\mathbf{Step }\mathbf{2}.$  $\mathbf{Hybrid}$ $\mathbf{analytic}$ $\mathbf{norm}$

        Now we introduce some notations. We denote $\mathbb{T}^{3}=\mathbb{R}^{3}/\mathbb{Z}^{3}.$ For function $f(x,v),$ we define the Fourier transform as follows.

        For a function $f=f(x),x\in\mathbb{T}^{d},$ we define its Fourier transform as follows:
 $$\hat{f}(k)=\frac{1}{(2\pi)^{d}}\int_{\mathbb{T}^{d}}f(x)e^{-ix\cdot k}dx,\quad k\in\mathbb{Z}^{d}.$$
 Similarly, for a function $f=f(v),v\in\mathbb{R}^{d},$ we define its Fourier transform by:
 $$\hat{f}(\xi)=\frac{1}{(2\pi)^{d}}\int_{\mathbb{R}^{d}}f(v)e^{-iv\cdot \xi}dv,\quad \xi\in\mathbb{R}^{d}.$$

 Finally, if $f=f(x,v), (x,v)\in\mathbb{T}^{d}\times\mathbb{R}^{d},$ we define its Fourier transform through the following formula:
 $$\hat{f}(k,\xi)=\frac{1}{(2\pi)^{d}}\int_{\mathbb{T}^{d}\times\mathbb{R}^{d}}f(x,v)e^{-ik\cdot x-iv\cdot\xi}dxdv,\quad (k,\xi)\in\mathbb{Z}^{d}\times\mathbb{R}^{d}.$$

 We shall also use the Fourier transform in time, if $f=f(t),t\in\mathbb{R},$ we denote
 $$\tilde{f}(\omega)=\int_{\mathbb{R}}f(t)e^{-it\omega}dt,\quad \omega\in\mathbb{C}.$$

Now we start to introduce the very important tools in our paper. These are time-shift pure and hybrid analytic norms. They are the same with those
           in the paper [28] written by Mouhot and Villani.
            \begin{de}(Hybrid analytic norms)
            $$\|f\|_{\mathcal{C}^{\lambda,\mu}}=\sum_{m,n\in\mathbb{N}_{0}^{3}}\frac{\lambda^{n}}{n!}\frac{\mu^{m}}{m!}
            \|\nabla^{m}_{x}\nabla_{v}^{n}f\|_{L^{\infty}(\mathbb{T}_{x}^{3}\times\mathbb{R}^{3}_{v})},\quad
            \|f\|_{\mathcal{F}^{\lambda,\mu}}=\sum_{k\in\mathbb{Z}^{3}}\int_{\mathbb{R}^{3}}|\tilde{f}(k,\eta)|
           e^{2\pi\lambda|\eta|}e^{2\pi\mu|k|}d\eta,$$
            $$\|f\|_{\mathcal{Z}^{\lambda,\mu}}=\sum_{l\in\mathbb{Z}^{3}}\sum_{n\in\mathbb{N}_{0}^{3}}
            \frac{\lambda^{n}}{n!}e^{2\pi\mu|l|}\|\widehat{\nabla_{v}^{n}f(l,v)}\|_{L^{\infty}(\mathbb{R}_{v}^{3})}.$$
            \end{de}
           \begin{de} (Time-shift pure and hybrid analytic norms)
For any $\lambda,\mu\geq0,p\in[1,\infty],$ we define
            $$\|f\|_{\mathcal{C}^{\lambda,\mu}_{\tau}}=\|f\circ S^{0}_{\tau}(x,v)\|_{\mathcal{C}^{\lambda,\mu}}=\sum_{m,n\in\mathbb{N}_{0}^{3}}\frac{\lambda^{n}}{n!}\frac{\mu^{m}}{m!}
            \|\nabla^{m}_{x}(\nabla_{v}+\tau\nabla_{x})^{n}f\|_{L^{\infty}(\mathbb{T}_{x}^{3}\times\mathbb{R}^{3}_{v})},$$
           $$\|f\|_{\mathcal{F}^{\lambda,\mu}_{\tau}}=\|f\circ S^{0}_{\tau}(x,v)\|_{\mathcal{F}^{\lambda,\mu}}=\sum_{k\in\mathbb{Z}^{3}}\int_{\mathbb{R}^{3}}|\tilde{f}(k,\eta)|
           e^{2\pi\lambda|k\tau+\eta|}e^{2\pi\mu|k|}d\eta,$$
         $$\|f\|_{\mathcal{Z}^{\lambda,\mu}_{\tau}}=\|f\circ S^{0}_{\tau}(x,v)\|_{\mathcal{Z}^{\lambda,\mu}}=\sum_{l\in\mathbb{Z}^{3}}\sum_{n\in\mathbb{N}_{0}^{3}}
            \frac{\lambda^{n}}{n!}e^{2\pi\mu|l|}\|(\nabla_{v}+2i\pi\tau\cdot l)^{n}\hat{f}(l,v)\|
            _{L^{\infty}(\mathbb{R}_{v}^{3})},$$
               $$\|f\|_{\mathcal{Z}^{\lambda,\mu;p}_{\tau}}=\sum_{l\in\mathbb{Z}^{3}}\sum_{n\in\mathbb{N}_{0}^{3}}
            \frac{\lambda^{n}}{n!}e^{2\pi\mu|l|}\|(\nabla_{v}+2i\pi\tau\cdot l)^{n}\hat{f}(l,v)\|_{L^{p}(\mathbb{R}_{v}^{3})},$$
            $$\|f\|_{\mathcal{Y}^{\lambda,\mu}_{\tau}}=\|f\|_{\mathcal{F}^{\lambda,\mu;\infty}_{\tau}}
            =\sup_{k\in\mathbb{Z}^{3},\eta\in\mathbb{R}^{3}}e^{2\pi\mu|k|}e^{2\pi\lambda|\eta+k\tau|}|\hat{f}(k,\eta)|.$$
           \end{de}

           From the above  definitions, we can state some simple and important propositions, and the related proofs can be found in [26.28], so we remove the proofs.
\begin{prop} For any $\tau\in\mathbb{R},\lambda,\mu\geq 0,$
\begin{itemize}
 \item[{(i)}]   if $f$ is a function only of $x,$ then
           $\|f\|_{\mathcal{C}^{\lambda,\mu}_{\tau}}=\|f\|_{\mathcal{C}^{\lambda|\tau|+\mu}},\|f\|_{\mathcal{F}^{\lambda,\mu}_{\tau}}
           =\|f\|_{\mathcal{Z}^{\lambda,\mu}_{\tau}}=\|f\|_{\mathcal{F}^{\lambda|\tau|+\mu}};$
          \item[{(ii)}]    if $f$ is a function only of $v,$ then
           $\|f\|_{\mathcal{C}^{\lambda,\mu;p}_{\tau}}=\|f\|_{\mathcal{Z}^{\lambda,\mu;p}_{\tau}}=\|f\|_{\mathcal{C}^{\lambda,;p}};$
            \item[{(iii)}] for any $\lambda>0,$ then
            $\|f\circ(Id+G)\|_{\mathcal{F}^{\lambda}}\leq\|f\|_{\mathcal{F}^{\lambda+\nu}},\nu=\|G\|_{\dot{\mathcal{F}}^{\lambda}};$
            \item[{(iv)}] for any $\bar{\lambda}>\lambda,p\in[1,\infty],$ $\|\nabla f\|_{\mathcal{C}^{\lambda;p}}\leq
            \frac{1}{\lambda e\log(\frac{\bar{\lambda}}{\lambda} )}
            \| f\|_{\mathcal{C}^{\bar{\lambda};p}},$
            $\|\nabla f\|_{\mathcal{F}^{\lambda;p}}\leq\frac{1}{2\pi e(\bar{\lambda}-\lambda )}
            \| f\|_{\mathcal{F}^{\bar{\lambda};p}},$
            \item[{(v)}]  for any $\bar{\lambda}>\lambda>0,\mu>0,$ then
            $\|vf\|_{\mathcal{Z}^{\lambda,\mu;1}_{\tau}}\leq\|f\|_{\mathcal{Z}^{\bar{\lambda},\mu;1}_{\tau}};$
           \item[{(vi)}]  for any $\bar{\lambda}>\lambda,\bar{\mu}>\mu,$
          $\|\nabla_{v}f\|_{\mathcal{Z}_{\tau}^{\lambda,\mu;p}}
          \leq C\bigg(\frac{1}{\lambda\log(\frac{\bar{\lambda}}{\lambda})}\|f\|_{\mathcal{Z}_{\tau}^{\bar{\lambda},\bar{\mu};p}}+\frac{\tau}{\bar{\mu}-\mu}
 \|f\|_{\dot{\mathcal{Z}}_{\tau}^{\bar{\lambda},\bar{\mu};p}}\bigg);$
 \item[{(vii)}] for any $\bar{\lambda}>\lambda,$
 $\|(\nabla_{v}+\tau\nabla_{x})f\|_{\mathcal{Z}_{\tau}^{\lambda,\mu;p}}\leq \frac{1}{C\lambda\log(\frac{\bar{\lambda}}{\lambda})}
 \|f\|_{\mathcal{Z}_{\tau}^{\bar{\lambda},\mu;p}};$
  \item[{(viii)}] for any $\bar{\lambda}\geq\lambda\geq0,\bar{\mu}\geq\mu\geq0,$ then $\|f\|_{\mathcal{Z}_{\tau}^{\lambda,\mu}}
  \leq_{\mathcal{Z}_{\tau}^{\bar{\lambda},\bar{\mu}}}.$ Moreover, for any $\tau,\bar{\tau}\in\mathbb{R},$ $p\in[1,\infty],$ we have
  $\|f\|_{\mathcal{Z}_{\tau}^{\lambda,\mu;p}}\leq\|f\|_{\mathcal{Z}_{\bar{\tau}}^{\lambda,\mu+\lambda|\tau-\bar{\tau}|;p}};$
  \item[{(viiii)}] $\|f\|_{\mathcal{Y}^{\lambda,\mu}_{\tau}}\leq\|f\|_{\mathcal{Z}^{\lambda,\mu;1}_{\tau}};$
  \item[{(iX)}] for any function $f=f(x,v),$ $\|\int_{\mathbb{R}^{3}}fdv\|_{\mathcal{F}^{\lambda|\tau|+\mu}}\leq\|f\|_{\mathcal{Z}^{\lambda,\mu;1}_{\tau}}.$
 \end{itemize}
\end{prop}

           \begin{prop} For any $X\in\{\mathcal{C},\mathcal{F},\mathcal{Z}\}$ and any $t,\tau\in\mathbb{R},$
           $$\|f\circ S^{0}_{\tau}\|_{X^{\lambda,\mu}_{\tau}}=\|f\|_{X^{\lambda,\mu}_{t+\tau}}.$$

           \end{prop}
           \begin{lem} Let $\lambda,\mu\geq0,t\in\mathbb{R},$ and consider two functions $F,G:\mathbb{T}^{3}\times\mathbb{R}^{3}\rightarrow\mathbb{T}^{3}\times\mathbb{R}^{3}.$
           Then there is $ \varepsilon\in(0,\frac{1}{2})$ such that  if $F,G$ satisfy
\begin{align}
\|\nabla(F-Id)\|_{\mathcal{Z}^{\lambda,\mu}_{\tau}}\leq\varepsilon,
\end{align}
where $\lambda=\lambda +2\|F-G\|_{\mathcal{Z}^{\lambda ,\mu}_{\tau}},\quad\mu=\mu+2(1+|\tau|)\|F-G\|_{\mathcal{Z}^{\lambda ,\mu}_{\tau}} ,$ then $F$ is invertible  and
\begin{align}
&\|F^{-1}\circ G-Id\|_{\mathcal{Z}^{\lambda ,\mu}_{\tau}}\leq2\|F-G\|_{\mathcal{Z}^{\lambda ,\mu}_{\tau}}.
\end{align}
\end{lem}

\begin{prop} For any $\lambda,\mu\geq0$ and any $p\in[1,\infty],\tau\in\mathbb{R},\sigma\in\mathbb{R},a\in\mathbb{R}\setminus\{0\}$ and $b\in\mathbb{R},$
we have
$$\|f(x+bv+X(x,v),av+V(x,v))\|_{\mathcal{Z}^{\lambda,\mu;p}_{\tau}}\leq|a|^{-\frac{3}{p}}\|f\|_{\mathcal{Z}^{\alpha,\beta;p}_{\sigma}},$$
where $\alpha=\lambda|a|+\|V\|_{\mathcal{Z}^{\lambda,\mu}_{\tau}},\quad \beta=\mu+\lambda|b+\tau-a\sigma|+\|X-\sigma V\|_{\mathcal{Z}^{\lambda,\mu}_{\tau}}.$
\end{prop}
\begin{lem} Let $G=G(x,v)$ and $R=R(x,v)$ be valued in $\mathbb{R},$ and
$\beta(x)=\int_{\mathbb{R}^{3}}(G\cdot R)(x,v)dv.$
Then for any $\lambda,\mu,t\geq0$ and any $b>-1,$ we have
$$\|\beta\|_{\mathcal{F}^{\lambda t+\mu}}\leq3\|G\|_{\mathcal{Z}_{\tau-\frac{bt}{1+b}}^{\lambda(1+b),\mu;1}}
\|R\|_{\mathcal{Z}_{\tau-\frac{bt}{1+b}}^{\lambda(1+b),\mu}}.$$
\end{lem}
%To explain our result in this paper, we need to quantify regularity corrections to the analytic regularity in the $ x $ variable.
%\begin{de} For $\lambda,\mu,\gamma,\eta\geq0,\tau\in\mathbb{R},p\in[1,\infty],$ we define
%$$\|f\|_{\mathcal{Z}^{(\lambda,\eta),(\mu,\gamma);p}_{\tau}}=\sum_{l\in\mathbb{Z}^{d}}\sum_{n\in\mathbb{N}_{0}^{d}}
%\frac{\lambda^{n}}{n!}e^{2\pi\mu|l|}(1+|l|)^{\gamma}\|(\nabla_{v}+2i\pi\tau l)^{n}\hat{f}(l,v)(1+|v|)^{\mu}\|_{L^{p}(\mathbb{R}_{v}^{d})},$$
%$$\|f\|_{\mathcal{F}^{(\mu,\gamma);p}_{\tau}}=\sum_{l\in\mathbb{Z}^{d}}
%e^{2\pi\mu|l|}(1+|l|)^{\gamma}\|\hat{f}(l,v)\|_{L^{p}}.$$
%\end{de}
%To explain our result in this paper, we need to quantify regularity corrections to the analytic regularity in the $ x $ variable.
%\begin{de} For $\lambda,\mu,\gamma,\eta\geq0,\tau\in\mathbb{R},p\in[1,\infty],$ we define
%$$\|f\|_{\mathcal{Z}^{(\lambda,\eta),(\mu,\gamma);p}_{\tau}}=\sum_{l\in\mathbb{Z}^{d}}\sum_{n\in\mathbb{N}_{0}^{d}}
%\frac{\lambda^{n}}{n!}e^{2\pi\mu|l|}(1+|l|)^{\gamma}\|(\nabla_{v}+2i\pi\tau l)^{n}\hat{f}(l,v)(1+|v|)^{\mu}\|_{L^{p}(\mathbb{R}_{v}^{d})},$$
%$$\|f\|_{\mathcal{F}^{(\mu,\gamma);p}_{\tau}}=\sum_{l\in\mathbb{Z}^{d}}
%e^{2\pi\mu|l|}(1+|l|)^{\gamma}\|\hat{f}(l,v)\|_{L^{p}}.$$
%\end{de}

\begin{center}
\item\section{Linearized Landau damping in  weakly collisional plasma}
%{\bf\large 1. \quad Introduction }
\end{center}

In this section, we consider the linearized Vlasov-Poisson equations in the weakly collisional case as follows:
         \begin{equation}
\left\{\begin{array}{l}
\partial_{t}f+v\cdot\nabla_{x}f+\frac{q}{m}E\cdot\nabla_{v}f^{0}=\nu\mathcal{C}(f),\\
E=\nabla W(x)\ast\rho(t,x),\quad \mathcal{C}(f)=\rho f^{0}-f, \\
\rho(t,x)=\int_{\mathbb{R}^{3}}f(t,x,v)dv,\\
f(0,x,v)=f_{0}(x,v),f^{0}=f^{0}(v),
\end{array}\right.
 \end{equation}
for any $\nu\in [0,\nu_{0}],\nu_{0} $ some small constant.

\begin{thm}
Consider equations (2.1). For any $\eta,v\in\mathbb{R}^{3},k\in\mathbb{N}_{0}^{3},$  assume that the following conditions hold:
%$\quad\mathbf{ how} \quad\mathbf{ to}\quad \mathbf{make} \quad \mathbf{assumptions }$
%$\quad \mathbf{on}\quad f^{0}$ is suitable
\begin{itemize}
        \item[(i)] $ \widehat{W}(0)=0$ where $|\widehat{W}(k)|\leq\frac{1}{1+|k|^{\gamma}},\gamma>1;$
%\item[(i)]$f^{0}(v)$ is  Maxwellian, that is, $f^{0}(v)=f^{0}(v_{\perp},v_{z})=C_{M}exp(-\alpha(v^{2}_{\perp}+v^{2}_{z}))$ for some $\alpha>0$  ;
         \item[(ii)] $||f^{0}||_{\mathcal{C}^{\lambda;1}}\leq C_{0},$
              for some constants $\lambda,C_{0}>0;$

         \item[(iii)]$||f_{0}||_{\mathcal{Z}^{\lambda,\mu;1}}\leq \delta_{0}$ for some constants $\mu>0,\delta_{0}>0$;

             \item[(iv)]  for any velocity $v\in\mathbb{R}^{3},$
             there exists some positive constant $v_{Te}\in\mathbb{R}$ such that if $v=\frac{\omega}{k}+i\frac{\nu}{2\pi k},$ then $|v|\gg v_{Te}.$
                  %\item[(v)] $f^{0}(\eta_{1},\eta_{2},v_{3})\leq Ce^{-\lambda_{1}|\eta_{1}|}e^{-\lambda_{1}|\eta_{2}|}e^{-\alpha_{1}|v_{3}|},$ otherwise $f^{0}(\eta_{1},\eta_{2},v_{3})<\varepsilon(v_{3}),$ for some $0<\varepsilon(v_{3})<1$ sufficient small.%$\inf_{k\in\mathbb{Z}^{3}}|\tilde{\mathcal{L}}(\omega,k)-1|>\kappa,$ for some $0<\kappa<1.$
         \end{itemize}
         Then %if $f_{0}(x,v),f^{0}(v)$ are axisymetric, then the solution $f$ of (1.1) is also axisymmetric.
for any fixed $\eta,k,$  we have
         \begin{align}
         &|\hat{f}(t,k,\eta)-\hat{f}_{0}(k,\eta)|\leq C(C_{0},\delta_{0}) e^{-2\pi\lambda|\eta+kt|}e^{-2\pi\mu |k|}\max\bigg\{e^{-\nu t},\frac{1-e^{-\nu t}}{\nu},1-e^{-\nu t}\bigg\},\notag\\
         &|\hat{\rho}(t,k)-\hat{\rho}_{0}|\leq C(C_{0},\delta_{0}) e^{-2\pi(\lambda t+\mu)|k|},\notag\\
        & |\hat{E}(t,k)|\leq C(C_{0},\delta_{0})e^{-2\pi\lambda|k|t}e^{-2\pi\mu|k| }.
         \end{align}
        where $\rho_{0}=\int_{\mathbb{T}^{3}}\int_{\mathbb{R}^{3}}f_{0}(x,v)dvdx.$
         %(这里可能不需要假设初值为轴对称在速度方向 %maybe not assume $f_{0}$ axisymmetric in $v$)
         \end{thm}

\begin{rem} In the linear case, before proving Theorem 2.1, we recall the result of I.Tristani in [31].
We define
   \begin{equation}
\left\{\begin{array}{l}
K^{0}_{\nu}(t,k):=\nu e^{-\nu t}\hat{f}^{0}(kt),\\
K^{1}_{\nu}(t,k):=-\hat{W}(k) e^{-\nu t}\hat{f}^{0}(kt)|k|^{2}t, \\
K_{\nu}(t,k):=K^{0}_{\nu}(t,k)+K^{1}_{\nu}(t,k),\\
\mathcal{L}_{\nu}(\eta,k)=\int^{\infty}_{0} e^{2\pi\eta^{\ast}|k|t}K_{\nu}(t,k)dt,
\end{array}\right.
 \end{equation}

where $\eta^{\ast}$ is the complex conjugate to $\eta.$ The stability condition of I.Tristani'work is as follows::

 there exists $\varepsilon_{0}>0$ such that  $\mathcal{L}_{\nu}$ satisfies the following condition: for some constant $ \kappa>0,$
 $\forall\nu\in[0,\nu_{0}],$

 $$\mathbf{(H)}\quad \quad \inf_{k\in\mathbb{Z}^{d}\setminus\{0\}}\inf_{\textmd{Im}\eta\leq0}|1-\mathcal{L}_{\nu}(\eta,k)|\geq\kappa.$$
 Comparing $\mathbf{(H)}$ condition, our stability condition in Theorem 2.1 are suitable for  the physical intuitive  from the energy viewpoint. From the condition of the classical KAM theory, our condition is in correspondence to the Diophantus condition in KAM theory  in some sense.
 \end{rem}
$$$$

   The proof of Theorem 2.1. Without loss of generality, we assume $t\geq0.$ we consider (2.1) as a perturbation of free transport and apply the Duhamel's formula to get
 $$f(t,x,v)=f_{0}(x-vt,v)e^{-\nu t}+\nu\int^{t}_{0}e^{-\nu(t-s)}\rho(s,x-v(t-s))f^{0}(v)ds$$
 $$-\int^{t}_{0}e^{-\nu(t-s)}(E\cdot\nabla_{v}f^{0})(s,x-v(t-s),v)ds.$$

Then we take the Fourier transform in both variables $(x,v),$
$$\hat{f}(t,k,\xi)=e^{-\nu t}\hat{f}_{0}(k,\xi+kt)$$
$$+\frac{\nu}{(2\pi)^{d}}\int^{t}_{0}\int_{\mathbb{T}^{d}\times\mathbb{R}^{d}}e^{-\nu(t-s)}\rho(s,x)f^{0}(v)e^{-ik\cdot x-iv\cdot(\xi+k(t-s))}dxdvds$$
$$-\frac{1}{(2\pi)^{d}}\int^{t}_{0}\int_{\mathbb{T}^{d}\times\mathbb{R}^{d}}e^{-\nu(t-s)}(E\cdot\nabla_{v}f^{0})(s,x,v)e^{-ik\cdot x-iv\cdot(\xi+k(t-s))}dxdvds,$$
and from which we can deduce
$$\hat{f}(t,k,\xi)=e^{-\nu t}\hat{f}_{0}(k,\xi+kt)+\nu\int^{t}_{0}e^{-\nu(t-s)}\hat{\rho}(s,k)\hat{f}^{0}(\xi+k(t-s))ds$$
\begin{equation}
-\int^{t}_{0}e^{-\nu(t-s)}k\cdot(\xi+k(t-s))\widehat{W}(k)\hat{\rho}(s,k)
\hat{f}^{0}(\xi+k(t-s))ds.
\end{equation}

Then taking $\xi=0,$ we obtain the closed equation on $\hat{\rho}(t,k):$
$$\hat{\rho}(t,k)=e^{-\nu t}\hat{f}_{0}(k,kt)+\nu\int^{t}_{0}e^{-\nu(t-s)}\rho(s,k)\hat{f}^{0}(k(t-s))ds$$
\begin{equation}
-\int^{t}_{0}e^{-\nu(t-s)}k\cdot(k(t-s))\widehat{W}(k)\hat{\rho}(s,k)
\hat{f}^{0}(k(t-s))ds.
\end{equation}

Recall the definition of $K_{\nu},$ we have
\begin{align}
\hat{\rho}(t,k)=e^{-\nu t}\hat{f}_{0}(k,kt)+\int^{t}_{0}K_{\nu}(t-s,k)\hat{\rho}(s,k)ds.
\end{align}

First we assume $k\neq0,$ and consider $\lambda>0,$ write $$\Phi(t,k)=\hat{\rho}(t,k)e^{2\pi\lambda|k|t}\quad\textmd{and}\quad A(t,k)=\hat{f}_{0}(t,k)e^{-\nu t}e^{2\pi\lambda|k|t};$$
then (0.6) becomes
\begin{align}
\Phi(t,k)=A(t,k)+\int^{t}_{0}K_{\nu}(t-s,k)e^{2\pi\lambda'|k|(t-s)}\Phi(s,k)ds.
\end{align}

We take the Fourier transform in  time variable, after extending $ K,$
$A$ and $\Phi$ by 0 at negative times. We have,
$$\tilde{\Phi}(\omega,k)=\tilde{A}(\omega,k)+\mathcal{L}_{\nu}(\omega,k)\tilde{\Phi}(\omega,k).$$

By the Stability condition, let  $\eta=\frac{\omega}{k}+i\frac{\nu}{2\pi k },$

\begin{align}
&\tilde{\mathcal{L}}_{\nu}(\omega,k)=-\frac{q}{m}\int_{\mathbb{R}^{+}}\int_{\mathbb{R}^{3}}e^{-\nu t}e^{2\pi\lambda|k|t}e^{2\pi it\omega}e^{-2\pi ik\cdot vt}
[k\widehat{W}(k)\partial_{v}f^{0}+\nu f^{0}(v)]dvdt\notag\\
\end{align}
$$\leq\sup_{\omega}\frac{q}{m}\bigg|\int_{\mathbb{R}^{+}}\int_{\mathbb{R}^{3}}e^{-\nu t}e^{2\pi it\omega}
e^{-2\pi ikvt}\cdot\sum_{n}\frac{|2\pi i\lambda|k|t|^{n}}{n!}
[k\widehat{W}(k)\nabla_{v}f^{0}(v)+\nu f^{0}(v)]dvdt\bigg|$$
%&=\frac{q}{m}\sum_{n}\frac{\lambda^{n}}{n!}\int^{\infty}_{0}|\frac{k_{3}}{\omega}\hat{W}_{2}\int_{\mathbb{R}^{+}}e^{it\omega}(i\eta_{k1})(ik_{3}t)^{n}
%e^{-ik_{3}v_{3}t}
%\widehat{v_{3}f^{0}}(\eta_{k1},\eta_{k2},v_{3})dv_{3}dt|d\omega\notag\\
%$$=\frac{q}{m}\sum_{n}\frac{\lambda^{n}}{n!}\int^{\infty}_{0}\bigg|\frac{k_{3}}{\omega}\hat{W}_{2}\int_{\mathbb{R}^{+}}\int_{\mathbb{R}}(2\pi i\eta_{k1})(2\pi ik_{3}t)^{n}
%e^{2\pi\lambda_{0}|\eta_{k1}|}e^{2\pi\lambda_{0}|\eta_{k2}|}e^{2\pi ik_{3}t(\frac{\omega}{k_{3}}-v_{3})}\cdot\widehat{v_{3}f^{0}}(\eta_{k1},\eta_{k2},v_{3})dv_{3}dt\bigg|d\omega$$
%$$=\frac{q}{m}\sum_{n}\frac{\lambda^{n}}{n!}\int^{\infty}_{0}\bigg|\frac{k_{3}}{\omega}\hat{W}_{2}\int_{\mathbb{R}^{+}}\int_{\mathbb{R}}(2\pi i\eta_{k1})
%(-1)^{n}e^{2\pi\lambda_{0}|\eta_{k1}|}e^{2\pi\lambda_{0}|\eta_{k2}|}\nabla_{v_{3}}^{n}e^{2\pi ik_{3}t(\frac{\omega}{k_{3}}-v_{3})}\cdot\widehat{v_{3}f^{0}}(\eta_{k1},\eta_{k2},v_{3})dv_{3}dt\bigg|d\omega$$
$$=\sup_{\omega}\frac{q}{m}\sum_{n}\frac{\lambda^{n}}{n!}\bigg|\int_{\mathbb{R}^{+}}\int_{\mathbb{R}^{3}}e^{-\nu t}
e^{2\pi ikt(\frac{\omega}{k}-v)}
\cdot[k\hat{W}(k)\nabla^{n+1}_{v}f^{0}(v)+\nu \nabla_{v}^{n}f^{0}(v)]dv dt\bigg|$$
$$=\sup_{\omega}\frac{q}{m}\sum_{n}\frac{\lambda^{n}}{n!}\bigg|
(-i(k\widehat{W}(k)\nabla_{\eta}^{n+1}f^{0}(\eta)+\nu \nabla_{\eta}^{n}f^{0}(\eta))\bigg |\leq\frac{q}{mv_{Te}}e^{-c_{0}v_{Te}},$$
where in the last inequality we use the stability condition (iv) that if $v=\frac{\omega}{k}+i\frac{\nu}{2\pi k},$ then $|v|\gg v_{Te},$ and   the assumption (i) and (ii).
Then  %for any $k,\omega,$ by the assumption (iv) of $Proposition$ 1.1, we obtain
%$$|1-\tilde{\mathcal{L}}(\omega,k)|\geq\kappa,$$
 there exists  some constant $0<\kappa<1$ such that
$\|\widetilde{\mathcal{L}}_{\nu}(\omega,k)\|_{L^{\infty}}\leq\kappa.$

Now we apply the Plancherel's identity to find (for each $k$)
$$\|\Phi\|_{L^{2}(dt)}\leq\frac{\|A\|_{L^{2}(dt)}}{\kappa}.$$

Then we plug this into the equation (2.7) to get
$$\|\Phi\|_{L^{\infty}(dt)}\leq\|A\|_{L^{\infty}(dt)}+\|K_{\nu}e^{2\pi\lambda|k|t}\|_{L^{2}(dt)}\|\Phi\|_{L^{2}(dt)}$$
\begin{align}
\leq\|A\|_{L^{\infty}(dt)}+\frac{\|K_{\nu}e^{2\pi\lambda|k|t}\|_{L^{2}(dt)}\|A\|_{L^{2}(dt)}}{\kappa}.
\end{align}
Through  simple computation,  we can obtain
$$\sup_{t\geq0}\bigg(\sum_{k\in\mathbb{Z}^{d}\setminus0}|\hat{\rho}(t,k)|e^{2\pi(\lambda t+\mu)|k|}\bigg)$$
$$\leq C(\lambda,\kappa)\sup_{t\geq0}\sup_{k\in\mathbb{Z}^{d}\setminus0}
|\widehat{f_{0}}(k,kt)|e^{2\pi(\lambda t+\mu)|k|}e^{-\nu t}$$
$$\leq C(\lambda,\kappa)\sup_{t\geq0}\sum_{k\in\mathbb{Z}^{d}\setminus0}|\widehat{f_{0}}(k,kt)|
e^{2\pi((\lambda t+\mu)|k|}e^{-\nu t}.$$

Equivalently, $$\sup_{t\geq0}\|\rho(t,\cdot)\|_{\dot{\mathcal{F}}^{\lambda t+\mu }}\leq C\sup_{t\geq0}\|\int_{\mathbb{R}^{d}}f_{0}\circ S_{-t}^{0}dv\|
_{\dot{\mathcal{F}}^{\lambda) t+\mu}}e^{-\nu t}\leq\sup_{t\geq0}\|f_{0}\circ S_{-t}^{0}\|_{\dot{\mathcal{Z}}_{t}^{\lambda,\mu;1}}e^{-\nu t}=\|f_{0}\|
_{\dot{\mathcal{Z}}_{0}^{\lambda,\mu;1}}\leq\delta_{0}.$$

We write
$$f(t,\cdot)=(f_{0}\circ S_{-t}^{0})e^{-\nu t}-\int^{t}_{0}e^{-(t-s)\nu}((\nabla W\ast\rho_{s})\circ S_{-(t-s)}^{0})\cdot\nabla_{v}f^{0}ds
+\nu\int^{t}_{0}e^{-(t-s)\nu}(\rho_{s}\circ S_{-(t-s)}^{0})\cdot\nabla_{v}f^{0}ds,$$
where $\rho_{s}=\rho(s,\cdot).$ Then  we have, for all $t\geq0,$
$$\|f\|_{\dot{\mathcal{Z}}^{\lambda,\mu;1}_{t}}\leq\|f_{0}\circ S_{-t}^{0}\|_{\dot{\mathcal{Z}}_{t}^{\lambda,\mu;1}}e^{-\nu t}
+\int^{t}_{0}e^{-(t-s)\nu}\|(\nabla W\ast\rho_{s})\circ S_{-(t-s)}^{0}\|_{\dot{\mathcal{Z}}^{\lambda,\mu;
\infty}_{t}}\|\nabla_{v}f^{0}\|_{\dot{\mathcal{Z}}^{\lambda,\mu;1}_{t}}ds$$
$$+\nu\int^{t}_{0}e^{-(t-s)\nu}\|\rho_{s}\circ S_{-(t-s)}^{0}\|_{\dot{\mathcal{Z}}^{\lambda,\mu;\infty}_{t}}
\|\nabla_{v}f^{0}\|_{\dot{\mathcal{Z}}^{\lambda,\mu;1}_{t}}ds$$
$$\quad\quad=\|f_{0}\|_{\dot{\mathcal{Z}}^{\lambda,\mu;1}}e^{-\nu t}+\|\nabla_{v}f^{0}\|_{\dot{\mathcal{C}}^{\lambda;1}}
\int^{t}_{0}e^{-(t-s)\nu}\|\nabla W\ast\rho_{s}\|_{\dot{\mathcal{F}}^{\lambda s+\mu}}ds$$
\begin{align}
+\nu\|\nabla_{v}f^{0}\|_{\dot{\mathcal{C}}^{\lambda;1}}\int^{t}_{0}e^{-(t-s)\nu}
\|\rho_{s}\|_{\dot{\mathcal{F}}^{\lambda s+\mu}}ds
\end{align}
Now we consider the case $k=0.$
\begin{align}
&\hat{\rho}(t,0)=e^{-\nu t}\hat{f}_{0}(0,0)+\nu\int^{t}_{0}e^{-\nu(t-s)}\hat{\rho}(s,0)ds,\\
&\hat{f}(t,0,\xi)=e^{-\nu t}\hat{f}_{0}(0,\xi)+\nu\int^{t}_{0}e^{-\nu(t-s)}\hat{\rho}(s,0)\hat{f}^{0}(\xi)ds.
\end{align}
If we assume $f_{0}$ a mean-zero distribution,
then we have $$\hat{\rho}(t,0)\equiv0,\quad \textmd{for all}\quad t\geq0,$$
and \begin{align}
\|\hat{f}(t,0,\xi)\|_{\mathcal{Z}^{\lambda,\mu;1}_{t}}\leq e^{-\nu t}\|\hat{f}_{0}(0,\xi)\|_{\mathcal{Z}^{\lambda,\mu;1}_{t}}.
\end{align}

         \begin{center}
\item\section{ Nonlinearized picture in  weakly collisional plasma}
%{\bf\large 1. \quad Introduction }
\end{center}

We next give the proof of the main theorem 0.1, stating the primary steps as propositions which are proved in subsections.

 \begin{center}
\item\subsection{The  Newton iteration}
%{\bf\large 1. \quad Introduction }
\end{center}
First of all, we write  a classical Newton iteration :
 Let
$$f^{0}=f^{0}(v) \quad\textrm{be\quad given},$$
and $$f^{n}=f^{0}+h^{1}+\ldots+h^{n},$$
where
 \begin{align}
  \left\{\begin{array}{l}
      \partial_{t}h^{1}+v\cdot\nabla_{x}h^{1}+E[h^{1}]\cdot\nabla_{v}f^{0}=\nu(\rho[h^{1}]f^{0}-h^{1}),\\
     h^{1}(0,x,v) =f_{0}-f^{0}, \\
   \end{array}\right.
           \end{align}
 and now we consider the Vlasov equation in step $n+1,n\geq1,$
 \begin{align}
  \left\{\begin{array}{l}
      \partial_{t}h^{n+1}+v\cdot\nabla_{x}h^{n+1}+E[f^{n}]\cdot\nabla_{v}h^{n+1} \\
      =-E[h^{n+1}]\cdot\nabla_{v}f^{n}-E[h^{n}]\cdot\nabla_{v}h^{n}+\nu(\rho[h^{n+1}]f^{0}-h^{n+1}),\\
      h^{n+1}(0,x,v) =0, \\
   \end{array}\right.
           \end{align}
the corresponding dynamical system is described as follows:
for any $(x,v)\in\mathbb{T}^{3}\times\mathbb{R}^{3},$
let $(X^{n}_{t,s},V^{n}_{t,s})$ as the solution of the following ordinary differential equations
  $$
   \left\{\begin{array}{l}
   \frac{d}{dt}X^{n+1}_{t,s}(x,v)=V^{n+1}_{t,s}(x,v),\\
 X^{n+1}_{s,s}(x,v)=x,
  \end{array}\right.
           $$

         \begin{align}
   \left\{\begin{array}{l}
   \frac{d}{dt}V^{n+1}_{t,s}(x,v)=E[f^{n}](t,X^{n}_{t,s}(x,v)),\\
 V^{n+1}_{s,s}(x,v)=v.
  \end{array}\right.
           \end{align}
            At the same time, we consider the corresponding  linear dynamics system as follows,
            \begin{align}
   \left\{\begin{array}{l}
  \frac{d}{dt}X^{0}_{t,s}(x,v)=V^{0}_{t,s}(x,v),\quad  \frac{d}{dt}V^{0}_{t,s}(x,v)=0,\\
 X^{0}_{s,s}(x,v)=x,\quad V^{0}_{s,s}(x,v)=v,
  \end{array}\right.
           \end{align}
             It is easy to check that
             $$\Omega^{n}_{t,s}-Id\triangleq(\delta X^{n}_{t,s},\delta V^{n}_{t,s})\circ (X^{0}_{s,t},V^{0}_{s,t})
             = ( X^{n}_{t,s}\circ(X^{0}_{s,t},V^{0}_{s,t})-Id, V^{n}_{t,s}\circ(X^{0}_{s,t},V^{0}_{s,t})-Id).$$
            Therefore, in order to estimate $( X^{n}_{t,s}\circ(X^{0}_{s,t},V^{0}_{s,t})-Id, V^{n}_{t,s}\circ(X^{0}_{s,t},V^{0}_{s,t})-Id),$
             we only need to study $(\delta X^{n}_{t,s},\delta V^{n}_{t,s})\circ (X^{0}_{s,t},V^{0}_{s,t}).$

             From Eqs.(3.3) and (3.4),
             $$
   \left\{\begin{array}{l}
   \frac{d}{dt}\delta X^{n+1}_{s,t}(x,v)=\delta V^{n+1}_{s,t}(x,v),\\
 \delta X^{n+1}_{s,s}(x,v)=0,
  \end{array}\right.
          $$
            \begin{align}
   \left\{\begin{array}{l}
   \frac{d}{dt}\delta V^{n+1}_{s,t}(x,v)=E[f^{n}](t,X^{n}_{s,t}(x,v)),\\
\delta V^{n}_{s,s}(x,v)=0.
  \end{array}\right.
           \end{align}

           Integrating (3.2) in time and $ h^{n+1}(0,x,v) =0,$ we get
            \begin{align}
            h^{n+1}(t,X^{n}_{t,0}(x,v),V^{n}_{t,0}(x,v))=\int^{t}_{0}e^{-\nu(t-s)}\Sigma^{n+1}(s,X^{n}_{s,0}(x,v),V^{n}_{s,0}(x,v))ds,
            \end{align}
           where$$\Sigma^{n+1}(t,x,v)=-E[h^{n+1}]\cdot\nabla_{v}f^{n}-E[h^{n}]\cdot\nabla_{v}h^{n}
          +\nu\rho[h^{n+1}]f^{0}.$$

           By the definition of $(X^{n}_{t,s}(x,v),V^{n}_{t,s}(x,v)),$ we have
           $$h^{n+1}(t,x,v)=\int^{t}_{0}e^{-\nu(t-s)}\Sigma^{n+1}(s,X^{n}_{s,t}(x,v),V^{n}_{s,t}(x,v))ds$$
           $$=\int^{t}_{0}e^{-\nu(t-s)}\Sigma^{n+1}(s,\delta X^{n}_{s,t}(x,v)+ X^{0}_{s,t}(x,v),\delta V^{n}_{s,t}(x,v)+ V^{0}_{s,t}(x,v))ds.$$

           Since the unknown $h^{n+1}$ appears on both sides of (3.6), we hope to get a self-consistent estimate. For this, we have little choice but to
           integrate in $v$ and get an integral equation on $\rho[h^{n+1}]=\int_{\mathbb{R}^{3}}h^{n+1}dv,$ namely

 \begin{align}
           &\rho[h^{n+1}](t,x)=\int^{t}_{0}\int_{\mathbb{R}^{3}}e^{-\nu(t-s)}(\Sigma^{n+1}\circ\Omega^{n}_{s,t}(x,v))(s, X^{0}_{s,t}(x,v),V^{0}_{s,t}(x,v))dvds\notag\\
           &=\int^{t}_{0}\int_{\mathbb{R}^{3}}-e^{-\nu(t-s)}\bigg[(\mathcal{E}^{n+1}_{s,t}\cdot G_{s,t}^{n})-(\mathcal{E}^{n}_{s,t}\cdot H^{n}_{s,t})+\nu(\rho[h^{n+1}]f^{0})\circ\Omega^{n}_{s,t}(x,v)\bigg](s, x-v(t-s),v)\notag\\
           &=I^{n+1,n}+II^{n,n}+III^{n+1,0},
           \end{align}
  where
  $$
  \left\{\begin{array}{l}
   \mathcal{E}^{n+1}_{s,t}=E[h^{n+1}]\circ\Omega^{n}_{s,t}(x,v),\quad \mathcal{E}^{n}_{s,t}=E[h^{n}]\circ\Omega^{n}_{s,t}(x,v),\\
   G_{s,t}^{n}=(\nabla_{v}f^{n})\circ\Omega^{n}_{s,t}(x,v),\quad H^{n}_{s,t}=(\nabla_{v}h^{n})\circ\Omega^{n}_{s,t}(x,v).\\
   \end{array}\right.
   $$

   \begin{center}
   \item\subsection{Inductive hypothesis}
%{\bf\large 1. \quad Introduction }
\end{center}

For n=1, from (3.1), it is known that (3.1) is a linear Vlasov equation. From section 2,  the conclusions of Theorem 0.1 hold.

 Now for any $i\leq n,i\in\mathbb{N}_{0},$ we assume that the following estimates hold,
 \begin{align}
 &\sup_{t\geq0}\|\rho[h^{i}](t,\cdot)\|_{\mathcal{F}^{\lambda_{i}t+\mu_{i}}},\notag\\
 &\sup_{0\leq s\leq t}\|h^{i}_{s}\circ\Omega^{i-1}_{t,s}\|
_{\mathcal{Z}_{s-\frac{bt}{1+b}}^{\lambda_{i}(1+b),\mu_{i};1}}\leq\delta_{i},\notag\\
\end{align}
then we have the following inequalities, denote $(\mathbf{E}^{n}):$
$$\sup_{t\geq 0}\|E[h^{i}](t,\cdot)\|_{\mathcal{F}^{\lambda_{i}t+\mu_{i}}}<\delta_{i},$$
$$\sup_{0\leq s\leq t}\|\nabla_{x}(h^{i}_{s}\circ\Omega^{i-1}_{t,s})\|
_{\mathcal{Z}_{s-\frac{bt}{1+b}}^{\lambda_{i}(1+b),\mu_{i};1}}\leq\delta_{i},\sup_{0\leq s\leq t}\|(\nabla_{x}h^{i}_{s})\circ\Omega^{i-1}_{t,s}\|
_{\mathcal{Z}_{s-\frac{bt}{1+b}}^{\lambda_{i}(1+b),\mu_{i};1}}\leq\delta_{i},$$
$$\|(\nabla_{v}+s\nabla_{x})(h^{i}_{s}\circ\Omega^{i-1}_{t,s})\|
_{\mathcal{Z}_{s-\frac{bt}{1+b}}^{\lambda_{i}(1+b),\mu_{i};1}}\leq\delta_{i},\|((\nabla_{v}+s\nabla_{x})h^{i}_{s})\circ\Omega^{i-1}_{t,s})\|
_{\mathcal{Z}_{s-\frac{bt}{1+b}}^{\lambda_{i}(1+b),\mu_{i};1}}\leq\delta_{i},$$
$$\sup_{0\leq s\leq t}\frac{1}{(1+s)^{2}}\|(\nabla\nabla h^{i}_{s})\circ\Omega^{i-1}_{t,s}\|
_{\mathcal{Z}_{s-\frac{bt}{1+b}}^{\lambda_{i}(1+b),\mu_{i};1}}\leq\delta_{i},$$
$$\sup_{0\leq s\leq t}(1+s)^{2}\|(\nabla_{v} h^{i})\circ\Omega^{i-1}_{t,s}-\nabla_{v}(h^{i}\circ\Omega^{i-1}_{t,s})\|
_{\mathcal{Z}_{s-\frac{bt}{1+b}}^{\lambda_{i}(1+b),\mu_{i};1}}\leq\delta_{i}.$$

It is easy to check  that the first inequality of $(\mathbf{E}^{n})$ holds since $E[h^{i}]$ satisfies the Poisson equation, so we only  need to show that the other inequalities of
$(\mathbf{E}^{n})$ also hold, the related proofs are found in section 4.
  \begin{center}
   \item\subsection{Local time iteration}
%{\bf\large 1. \quad Introduction }
\end{center}

Before working out the core of the proof of Theorem 0.1, we shall show a short time estimate,which will play a role as an initial data layer for the Newton scheme.
The main tool in this section is given by the following lemma.
\begin{lem} Let $f$ be an analytic function, $\lambda(t)=\lambda-Kt$ and $\mu(t)=\mu-Kt, K>0,$ let $T>0$ be so small that $\lambda(t)>0,
\mu(t)>0$ for $0\leq t\leq T.$ Then
for any $s\in[0,T]$ and any $p\geq1,$
$$\frac{d^{+}}{dt}\bigg|_{t=s}\|f\|_{\mathcal{Z}_{s}^{\lambda(t),\mu(t);p}}\leq-\frac{K}{(1+s)}\|\nabla f\|_{\mathcal{Z}_{s}^{\lambda(s),\mu(s);p}},$$
where $\frac{d^{+}}{dt}$ stands for the upper right derivative.
\end{lem}

\begin{prop} There exists some small constant $T>0,$ such that when  all  conditions of Theorem 0.1 hold, then
                  %\item[(v)] $f^{0}(\eta_{1},\eta_{2},v_{3})\leq Ce^{-\lambda_{1}|\eta_{1}|}e^{-\lambda_{1}|\eta_{2}|}e^{-\alpha_{1}|v_{3}|},$ otherwise $f^{0}(\eta_{1},\eta_{2},v_{3})<\varepsilon(v_{3}),$ for some $0<\varepsilon(v_{3})<1$ sufficient small.%$\inf_{k\in\mathbb{Z}^{3}}|\tilde{\mathcal{L}}(\omega,k)-1|>\kappa,$ for some $0<\kappa<1.$

        %if $f_{0}(x,v),f^{0}(v)$ are axisymetric, then the solution $f$ of (1.1) is also axisymmetric.
for any fixed $\eta,k,$  for all $t\in(0,T],$ $0<\lambda<\lambda_{0},$  we have
         \begin{align}
         &|\hat{f}(t,k,\eta)-\hat{f}_{0}(k,\eta)|\leq C(C_{0},\delta_{0}) e^{-2\pi\lambda|\eta+kt|}e^{-2\pi\mu |k|},\notag\\
         &|\hat{\rho}(t,k)-\hat{\rho}_{0}|\leq C(C_{0},\delta_{0}) e^{-2\pi(\lambda t+\mu)|k|},\quad |\hat{E}(t,k)|\leq C(C_{0},\delta_{0})e^{-2\pi\lambda|k|t}e^{-2\pi\mu|k| }.
         \end{align}
        where $\rho_{0}=\int_{\mathbb{T}^{3}}\int_{\mathbb{R}^{3}}f_{0}(x,v)dvdx.$
         %(这里可能不需要假设初值为轴对称在速度方向 %maybe not assume $f_{0}$ axisymmetric in $v$)
\end{prop}

\textbf{Proof.}  The first stage of the iteration,namely, $h^{1}$ was considered in $\S2.$ So we only need to care about the higher orders.
Recall that $f^{k}=f^{0}+h^{1}+\ldots+h^{k}.$ And we define
$$\lambda_{k}(t)=\lambda_{k}-2Kt\quad \textmd{and}\quad \mu_{k}(t)=\mu_{k}-Kt,$$
where $\{\lambda_{k}\}^{\infty}_{k=1}$ and $\{\mu_{k}\}^{\infty}_{k=1}$ are decreasing sequences of positive numbers.

We assume inductively that at stage $ n$ of the iteration, we have constructed $\{\lambda_{k}\}^{n}_{k=1},\{\mu_{k}\}^{n}_{k=1},\{\delta_{k}\}^{n}_{k=1}$ such that
$$\sup_{0\leq t\leq T}\|h^{k}(t,\cdot)\|_{\dot{\mathcal{Z}}^{\lambda_{k}(t),\mu_{k}(t);1}_{t}}\leq\delta_{k}\quad \textmd{for all}\quad  1\leq k\leq n$$
for some fixed $T>0.$

In the following we need to show that the induction hypothesis are satisfied at stage $n+1.$ For this, we have to construct $\lambda_{n+1},\mu_{n+1},\delta_{n+1}.$

Note that $h_{n+1}=0,$ at $t=0.$
 For $n\geq1,$ now let us solve
 $$ \partial_{t}h^{n+1}+v\cdot\nabla_{x}h^{n+1} =\widetilde{\Sigma}^{n+1},$$
 where $$\widetilde{\Sigma}^{n+1}=-E[f^{n}]\cdot\nabla_{v}h^{n+1}-E[h^{n+1}]\cdot\nabla_{v}f^{n}
-E[h^{n}]\cdot\nabla_{v}h^{n}+\nu(\rho[h^{n+1}]f^{0}-h^{n+1}).$$

 Hence, $$\|h^{n+1}\|_{\mathcal{Z}_{t}^{\lambda_{n+1}(t),\mu_{n+1}(t);1}}\leq\int^{t}_{0}e^{-\nu(t-s)}\|\overline{\Sigma}_{\tau}^{n+1}\circ S_{-(t-s)}^{0}\|
 _{\mathcal{Z}_{t}^{\lambda_{n+1}(t),\mu_{n+1}(t);1}}ds\leq\int^{t}_{0}e^{-\nu(t-s)}\|\overline{\Sigma}_{s}^{n+1}\|
 _{\mathcal{Z}_{s}^{\lambda_{n+1}(t),\mu_{n+1}(t);1}}ds,$$
  where $$\overline{\Sigma}^{n+1}=-E[f^{n}]\cdot\nabla_{v}h^{n+1}-E[h^{n+1}]\cdot\nabla_{v}f^{n}
-E[h^{n}]\cdot\nabla_{v}h^{n}+\nu\rho[h^{n+1}]f^{0}.$$
 then by Lemma 3.1,
 $$\frac{d^{+}}{dt}\|h^{n+1}\|_{\mathcal{Z}_{t}^{\lambda_{n+1}(t),\mu_{n+1}(t);1}}\leq-K\|\nabla_{x}h^{n+1}\|_{\mathcal{Z}_{t}
 ^{\lambda_{n+1},\mu_{n+1};1}}
 -K\|\nabla_{v}h^{n+1}\|_{\mathcal{Z}_{t}^{\lambda_{n+1},\mu_{n+1};1}}$$
 $$+\|E[f^{n}]\|_{\mathcal{F}^{\lambda_{n+1}t+\mu_{n+1}}}
 \|\nabla_{v}h^{n+1}\|_{\mathcal{Z}_{t}^{\lambda_{n+1},\mu_{n+1};1}}
+\|E[h^{n+1}]\|_{\mathcal{F}^{\lambda_{n+1}t+\mu_{n+1}}}
 \|\nabla_{v}f^{n}\|_{\mathcal{Z}_{t}^{\lambda_{n+1},\mu_{n+1};1}}$$
 $$+\|E[h^{n}]\|_{\mathcal{F}^{\lambda_{n+1}t+\mu_{n+1}}}
 \|\nabla_{v}h^{n}\|_{\mathcal{Z}_{t}^{\lambda_{n+1},\mu_{n+1};1}}+\nu\|\rho[h^{n+1}]\|_{\mathcal{F}^{\lambda_{n+1}t+\mu_{n+1}}}
 \|f^{0}\|_{\mathcal{Z}_{t}^{\lambda_{n+1},\mu_{n+1};1}}$$
 $$+\nu\int^{t}_{0}e^{-\nu(t-s)}\|\overline{\Sigma}_{s}^{n+1}\|
 _{\mathcal{Z}_{s}^{\lambda_{n+1}(t),\mu_{n+1}(t);1}}ds$$

 First, we easily get $\|E[h^{n}]\|_{\mathcal{F}^{\lambda_{n+1}t+\mu_{n+1}}}\leq C\|\nabla h^{n}\|_{\mathcal{Z}_{t}^{\lambda_{n+1},\mu_{n+1};1}}.$
Moreover,
 $$ \|\nabla_{v}f^{n}\|_{\mathcal{Z}_{t}^{\lambda_{n+1},\mu_{n+1};1}}\leq\sum^{n}_{i=1}\|\nabla_{v}h^{i}\|_{\mathcal{Z}_{t}^{\lambda_{n+1},\mu_{n+1};1}}
 \leq C\sum^{n}_{i=1}\frac{\|h^{i}\|_{\mathcal{Z}_{t}^{\lambda_{i+1},\mu_{i+1};1}}}{\min\{\lambda_{i}-\lambda_{n+1},\mu_{i}-\mu_{n+1}\}}.$$

 We gather the above estimates,
 $$\frac{d^{+}}{dt}\|h^{n+1}\|_{\mathcal{Z}_{t}^{\lambda_{n+1}(t),\mu_{n+1}(t);1}}\leq \bigg( C\sum^{n}_{i=1}\frac{\delta_{i}}
 {\min\{\lambda_{i}-\lambda_{n+1},\mu_{i}-\mu_{n+1}\}}-K+\nu C_{0}\bigg)\|\nabla h^{n+1}\|_{\mathcal{Z}_{t}^{\lambda_{n+1},\mu_{n+1};1}}$$
 $$
 +\frac{\delta^{2}_{n}}{\min\{\lambda_{n}-\lambda_{n+1},\mu_{n}-\mu_{n+1}\}}.$$

We may choose $$\delta_{n+1}=\frac{\delta^{2}_{n}}{\min\{\lambda_{n}-\bar{\lambda}_{n+1},\mu_{n}-\bar{\mu}_{n+1}\}},$$
if \begin{align}
C\sum^{n}_{i=1}\frac{\delta_{i}}
 {\min\{\lambda_{i}-\lambda_{n+1},\mu_{i}-\mu_{n+1}\}}\leq K-\nu C_{0}\end{align}
  holds.

 We choose $\lambda_{i}-\lambda_{i+1}=\mu_{i}-\mu_{i+1}=\frac{\Lambda}{i^{2}},$ where $\Lambda>0$ is arbitrarily small. Then for $i\leq n,\lambda_{i}-\lambda_{n+1}\geq
 \frac{\Lambda}{i^{2}},$ and $\delta_{n+1}\leq\delta^{2}_{n}n^{2}/\Lambda.$ Next we need to check that $\sum^{\infty}_{n=1}\delta_{n}n^{2}<\infty.$
 In fact, we choose $ K $ large enough and T small enough such that $\lambda_{0}-KT\geq \lambda_{\ast},\mu_{0}-KT\geq \mu_{\ast},$ and (3.9) holds, where
  $\lambda_{0}>\lambda_{\ast},\mu_{0}>\mu_{\ast}$ are fixed.

  If $\delta_{1}=\delta,$ then $\delta_{n}=n^{2}\frac{\delta^{2^{n}}}{\Lambda^{n}}(2^{2})^{2^{n-2}}(4^{2})^{2^{n-2}}\ldots((n-1)^{2})^{2}n^{2}.$
  To prove the sequence convergence for $\delta$ small enough, by induction that $\delta_{n}\leq z^{a^{n}},$ where $z$ small enough and $a\in(1,2).$ We
  claim that the conclusion holds for $n+1.$ Indeed, $\delta_{n+1}\leq\frac{z^{2a^{n}}}{\Lambda}n^{2}\leq z^{a^{n+1}}\frac{z^{(2-a)a^{n}}n^{2}}{\Lambda}.$ If $ z$ is
  so small that $z^{(2-a)a^{n}}\leq\frac{\Lambda}{n^{2}}$ for all $n\in\mathbb{N},$ then $\delta_{n+1}\leq z^{a^{n+1}},$ this concludes the local-time argument.

   \begin{center}
\item\subsection{Global time iteration }
%{\bf\large 1. \quad Introduction }
\end{center}

Based on the estimates of the local-time iteration, without loss of generality, sometimes  we only consider the case  $s\geq\frac{bt}{1+b},$ where $b$
is small enough.

First, we give deflection estimates that compare the free evolution with the true evolution for the particles trajectories.
  \begin{prop}
           Assume  for any $i\in\mathbb{N},0<i\leq n, $
           $$\sup_{t\geq 0}\|E[h^{i}](t,\cdot)\|_{\mathcal{F}^{\lambda_{i}t+\mu_{i}}}<\delta_{i}.$$

           And there exist  constants $\lambda_{\star}>0,\mu_{\star}>0$ such that $\lambda_{0}>\lambda'_{0}>\lambda_{1}
           >\lambda'_{1}>\ldots>\lambda_{i}>\lambda'_{i}>\ldots>\lambda_{\star},
           \mu_{0}>\mu_{1}>\mu'_{1}>\ldots>\mu_{i}>\mu'_{i}>\ldots>\mu_{\star}.$
           % satisfy the following conditions:

            %$ C(\lambda'_{k})\varepsilon e^{-2\pi(\lambda_{k}-\lambda'_{k})s}\min\bigg\{(s-\tau),\frac{1}{2\pi(\lambda_{k}-\lambda'_{k})}\bigg\}\leq\frac{1}{2}(\mu_{k+1}-\mu'_{k+1}),$

         % and

          %$|B_{0}|C(\lambda'_{k}) e^{-2\pi(\lambda_{k}-\lambda'_{k})s}\min\bigg\{(s-\tau),\frac{1}{2\pi(\lambda_{k}-\lambda'_{k})}\bigg\}\leq C(\lambda'_{k+1})e^{-2\pi(\lambda_{k+1}-\lambda'_{k+1})s},$   for all
          %$\tau\leq s\leq t.$

           % Meanwhile, $\mathbf{here}$  $\mathbf{assume}$  $\mathbf{that}$ $\mathbf{there}$ $\mathbf{exist}$ $\mathbf{a}$  $\mathbf{zero}$  $\mathbf{measure}$
            % $set$  $A$  and sufficiently large constant $V_{M}>0$ such that $$\sup_{\mathbb{R}^{3}\setminus A}|v|\leq V_{M}.$$
          Then we have %for any $k\in\mathbb{Z}^{3},i \in\mathbb{N}$
          $$\|\delta X^{n+1}_{t,s}\circ (X^{0}_{s,t},V^{0}_{s,t})\|_{\mathcal{Z}_{s-\frac{bt}{1+b}}^{\lambda'_{n},\mu'_{n}}}
          \leq C\sum^{n}_{i=1}\delta_{i}  e^{-\pi(\lambda_{i}-\lambda'_{i})s}\min\bigg\{\frac{(t-s)^{2}}{2},
           \frac{1}{2\pi(\lambda_{i}-\lambda'_{i})^{2}}\bigg\},$$

            $$\|\delta V^{n+1}_{t,s}\circ (X^{0}_{s,t},V^{0}_{s,t})\|_{\mathcal{Z}_{s-\frac{bt}{1+b}}^{\lambda'_{n},\mu'_{n}}}
            \leq C\sum^{n}_{i=1}\delta_{i} e^{-\pi(\lambda_{i}-\lambda'_{i})s}\min\bigg\{\frac{(t-s)}{2},
           \frac{1}{2\pi(\lambda_{i}-\lambda'_{i})}\bigg\},$$
           for $0<s<t,b=b(t,s)$ sufficiently small.
           %$$\|\delta X^{k}_{t,\tau}\|_{\mathcal{Z}_{t+\sigma}^{\lambda'_{k},\mu'_{k}}}\leq C(\lambda'_{k})\varepsilon e^{-2\pi(\lambda_{k}-\lambda'_{k})t}\min\bigg\{\frac{(t-\tau)^{2}}{2},\frac{1}{2\pi(\lambda_{k}-\lambda'_{k})^{2}}\bigg\},$$
    %$$\|\delta V^{k}_{t,\tau}\|_{\mathcal{Z}_{t+\sigma}^{\lambda'_{k},\mu'_{k}}}\leq C(\lambda'_{k})\varepsilon e^{-2\pi(\lambda_{k}-\lambda'_{k})t}\min\bigg\{(t-\tau),\frac{1}{2\pi(\lambda_{k}-\lambda'_{k})}\bigg\}.$$
           \end{prop}
           \begin{rem} From the above proposition, we know that weak collision has little impact on the trajectory of the plasma particles.
\end{rem}
\begin{prop} Under the assumptions of Proposition 3.3, then
           $$\bigg\|\nabla\Omega^{n+1}X_{t,s}-(Id,0)\bigg\|_{\mathcal{Z}_{s+\frac{bt}{1-b}}^{\lambda'_{n+1}(1-b),\mu'_{n+1}}}<\mathcal{C}_{1}^{n},
           \bigg\|\nabla\Omega^{n+1}V_{t,s}-(0,Id)\bigg\|_{\mathcal{Z}_{s+\frac{bt}{1-b}}^{\lambda'_{n+1}(1-b),\mu'_{n+1}}}<\mathcal{C}_{1}^{n}+
           \mathcal{C}_{2}^{n},$$
           where $\mathcal{C}_{1}^{n}= C\sum^{n}_{i=1}\frac{ e^{-\pi(\lambda_{i}-\lambda'_{i})s}\delta_{i}}
           {2\pi(\lambda_{i}-\lambda'_{i})^{2}}\min\bigg\{\frac{(t-s)^{2}}{2},
           1\bigg\},\mathcal{C}_{2}^{n}= C\sum^{n}_{i=1}\frac{ e^{-\pi(\lambda_{i}-\lambda'_{j})s}\delta_{i}}{2\pi(\lambda_{i}-\lambda'_{i})}
           \min\bigg\{t-s,1\bigg\}.$
           \end{prop}
           \begin{prop} Under the assumptions of Proposition 3.3, then
           $$\bigg\|\Omega^{i}X_{t,s}-\Omega^{n}X_{t,s}\bigg\|_{\mathcal{Z}_{s-\frac{bt}{1+b}}^{\lambda'_{n}(1-b),\mu'_{n}}}
           < \mathcal{C}_{1}^{i,n},\bigg\|\Omega^{i}V_{t,s}-\Omega^{n}V_{t,s}\bigg\|_{\mathcal{Z}_{s-\frac{bt}{1+b}}^{\lambda'_{n}(1-b),\mu'_{n}}}
           <\mathcal{C}_{1}^{i,n}+ \mathcal{C}_{2}^{i,n},$$
           where $\mathcal{C}_{1}^{i,n}=C\sum^{n}_{j=i+1} \frac{ e^{-\pi(\lambda_{j}-\lambda'_{j})s}\delta_{j}}
           {2\pi(\lambda_{j}-\lambda'_{j})^{2}}\min\bigg\{\frac{(t-s)^{2}}{2},
           1\bigg\},\mathcal{C}_{2}^{i,n}=C\sum^{n}_{j=i+1}
           \frac{ e^{-\pi(\lambda_{j}-\lambda'_{j})s}\delta_{j}}{2\pi(\lambda_{j}-\lambda'_{j})}
           \min\bigg\{t-s,1\bigg\}.$
           \end{prop}
           \begin{rem} Note that $\mathcal{C}_{1}^{i,n},\mathcal{C}_{2}^{i,n}$  decay fast as $s\rightarrow\infty,i\rightarrow\infty,$ and uniformly in
           $n\geq i,$ since the sequence $\{\delta_{n}\}^{\infty}_{n=1}$ has fast convergence. Hence, if $r\in\mathbb{N}$  given, we shall have
           \begin{align}
           \mathcal{C}_{1}^{i,n}\leq\omega_{i,n}^{r,1},\quad \textmd{and}\quad \mathcal{C}_{2}^{i,n}\leq\omega_{i,n}^{r,2},\quad \textmd{all}\quad r\geq1,
           \end{align}
           with $\omega_{i,n}^{r,1}=C^{r}_{\omega}\sum^{n}_{j=i+1} \frac{ \delta_{j}}
           {2\pi(\lambda_{j}-\lambda'_{j})^{2+r}}\frac{\min\{\frac{(t-s)^{2}}{2},
           1\}}{(1+s)^{r}}$ and $\omega_{i,n}^{r,2}=C^{r}_{\omega}\sum^{n}_{j=i+1} \frac{ \delta_{j}}
           {2\pi(\lambda_{j}-\lambda'_{j})^{1+r}}\frac{\min\{\frac{(t-s)^{2}}{2},
           1\}}{(1+s)^{r}},$ for some absolute constant $C^{r}_{\omega}$ depending only on $ r.$
           \end{rem}
           \begin{prop} Under the assumptions of Proposition 3.3, then
           $$\bigg\|(\Omega^{i}_{t,s})^{-1}\circ\Omega_{t,s}^{n}-Id\bigg\|_{\mathcal{Z}_{s-\frac{bt}{1+b}}^{\lambda'_{n}(1-b),\mu'_{n}}}
            <\mathcal{C}_{1}^{i,n}+ \mathcal{C}_{2}^{i,n}.$$
          \end{prop}

            To give a self-consistent estimate, we have to control each term of Eq.(3.7): I,II,III. And the most difficult term is $I,$
            because there is some resonance phenomena occurring in this term that makes the  propagated wave away from equilibrium.

            Let us first consider the first  term I.
 \begin{align}
 &I^{n+1,n}(t,x)=\int^{t}_{0}\int_{\mathbb{R}^{3}}-e^{-\nu(t-s)}(\mathcal{E}^{n+1}_{s,t}\cdot G_{s,t}^{n})(s, x-v(t-s),v)dvds.
 \end{align}

 To handle this term, we start by introducing
\begin{align}
&\bar{G}^{n}_{s,t}=\nabla_{v}f^{0}+\sum^{n}_{i=1}\nabla_{v}(h^{i}\circ\Omega^{i-1}_{s,t}),
\end{align}
and the error terms $\mathcal{R}_{0},\tilde{\mathcal{R}}_{0}$ are defined by

\begin{align}
\mathcal{R}_{0}=\int^{t}_{0}\int_{\mathbb{R}^{3}}e^{-\nu(t-s)}((E[h^{n+1}]\circ\Omega^{n}_{s,t}(x,v)-E[h^{n+1}])\cdot G_{s,t}^{n})(s, x-v(t-s),v)dvds,
\end{align}
\begin{align}
\tilde{\mathcal{R}}_{0}=\int^{t}_{0}\int_{\mathbb{R}^{3}}e^{-\nu(t-s)}(E[h^{n+1}]\cdot (G_{s,t}^{n}-\bar{G}_{s,t}^{n}))(s, x-v(t-s),v)dvds,
\end{align}
then we can decompose $$I^{n+1,n}=\bar{I}^{n+1,n}+\mathcal{R}_{0}+\tilde{\mathcal{R}}_{0}.$$
We decompose as
$$\bar{I}^{n+1,n}=\bar{I}_{0}^{n+1,n}+\sum^{n}_{i=1}\bar{I}_{i}^{n+1,n},$$
where
$$\bar{I}_{0}^{n+1,n}(t,x)=\int^{t}_{0}\int_{\mathbb{R}^{3}}E[h^{n+1}](s, x-v(t-s))\cdot\nabla_{v}f^{0}(v)dvd\tau,$$
$$\bar{I}_{i}^{n+1,n}(t,x)=\int^{t}_{0}\int_{\mathbb{R}^{3}}E[h^{n+1}](s, x-v(t-s))\cdot(\nabla_{v}h_{s}^{i}\circ\Omega^{i-1}_{t,s})(s, x-v(t-s),v))dvds.$$

\begin{prop} ($\textbf{M}$ain term $I$)
Assume %$\lambda^{\star}>\lambda_{1}>\lambda_{2}>\ldots>\lambda_{n}>\ldots>0,\mu^{\star}>\mu_{1}>\mu_{2}>\ldots>\mu_{n}>\ldots>0,\mu'_{n+1}=\mu_{n+1}+\eta(\frac{t-\tau}{1+t}),
%\nu'_{n+1}=\lambda_{n+1}(1+b)|\tau-\frac{bt}{1+b}|+\mu'_{n+1},$ where
 $b(t,s)\geq0$ small. And there exist  constants $\lambda_{\star}>0,\mu_{\star}>0$ such that $\lambda_{0}>\lambda'_{0}
 >\lambda_{1}>\lambda'_{1}>\ldots>\lambda_{i}>\lambda'_{i}>\ldots>\lambda_{\star},
 \mu_{0}>\mu_{1}>\mu''_{1}>\mu'_{1}>\ldots>\mu_{i}>\mu''_{i}>\mu'_{i}>\ldots>\mu_{\star}.$

We have
$$\| \bar{I}_{i}^{n+1,n}(t,\cdot)\|_{\mathcal{F}^{\lambda'_{n} t+\mu'_{n}}}$$
 $$\leq C\int^{t}_{0}e^{-(t-s)\nu} K^{n+1}_{1}(t,s)\| \nabla_{v}
(h^{i}_{s}\circ\Omega^{i-1}_{t,s})-\langle\nabla_{v}
(h^{i}_{s}\circ\Omega^{i-1}_{t,s})\rangle\|_{\mathcal{Z}^{\lambda'_{i}(1+b),
\mu'_{i};1}_{s-\frac{bt}{1+b}}}$$
$$
\| E[h^{n+1}]\|_{_{\mathcal{F}^{\bar{\nu}}}}ds+\int^{t}_{0}e^{-(t-s)\nu} K^{n+1}_{0}(t,s)\| \langle\nabla_{v}
(h^{i}_{s}\circ\Omega^{i-1}_{t,s})\rangle\|_{\mathcal{C}^{\lambda_{i}'(1+b);1}}\cdot\| E[h^{n+1}]\|_{_{\mathcal{F}^{\bar{\nu}}}}ds,$$
 where
$$\bar{\nu}=\max\bigg\{\lambda_{n}'s+\mu''_{n}-\frac{1}{2}\lambda_{n}'b(t-s),0\bigg\},$$
$$K^{n}_{0}(t,s)=e^{-\pi(\lambda'_{i}-\lambda'_{n})(t-s)},$$
$$K^{n+1}_{1}(t,s)=\sup_{k,l\in\mathbb{Z}^{3}}e^{-2\pi(\mu'_{i}-\mu'_{n})|l|}e^{-\pi(\lambda'_{i}-\lambda'_{n})|k(t-s)+ls|}
  e^{-2\pi(\frac{\lambda'_{n}}{2}(s-s')+\mu''_{n}-\mu'_{n})|k-l|}.$$
\end{prop}

%$$\leq\int^{\tau}_{0}e^{-2\pi(\nu''_{n+1}\tau-\nu'_{n+1}\tau)}\| \rho[h^{n+1}]\|_{\mathcal{\dot{F}}^{\nu''_{n+1}}}d\tau,$$
%where we set $\nu''_{n+1}=\lambda_{n+1}(1+b)|\tau-\frac{bt}{1+b}|+\mu''_{n+1},\mu''_{n+1}=\mu_{n+1}+2\eta(\frac{t-\tau}{1+t}).$
%$\mathbf{maybe}$ $\mathbf{we}$ $\mathbf{need}$ $\mathbf{take}$ $\mathcal{F}$ $\mathbf{norm}$ on $\mathbf{x_{\perp}}$ in order to remove $\nabla_{x_{\perp}}.$

%Indeed, leaving apart the small-time case, we assume $\tau\geq\frac{bt}{1+b},$ then
%$$\nu''_{n+1}=\lambda_{n+1}(1+b)|\tau-\frac{bt}{1+b}|+\mu_{n+1}+2\eta(\frac{t-\tau}{1+t})$$
%$$=\lambda_{n+1}\tau+\mu_{n+1}+(\frac{2\eta}{1+t}-\lambda_{n+1}b)(t-\tau),$$
%which is bounded by $\lambda_{n+1}\tau+\mu_{n+1}$ if we choose $\eta\leq\frac{1}{2}\lambda^{\star}(1+t).$

%So we obtain $$\|B[h^{n+1}]\|_{\mathcal{F}^{\nu'_{n+1}}}\leq\int^{t}_{0}e^{-2\pi(\nu''_{n+1}\tau-\nu'_{n+1}\tau)}\| \rho[h^{n+1}]\|_{\mathcal{\dot{F}}^{\lambda_{n+1}\tau+\mu_{n+1}}}d\tau.$$
\begin{cor}
From the above statement, we have
$$
\|\bar{I}_{i}^{n+1,n}(t,\cdot)\|_{\mathcal{F}^{\lambda'_{n} t+\mu'_{n}}}
\leq\int^{t}_{0}e^{-(t-s)\nu}K^{n+1}_{0}(t,s)\delta_{i}\| \rho[h^{n+1}]\|_{\mathcal{F}^{\lambda'_{n}s+\mu'_{n}}}ds$$
$$+\int^{t}_{0}e^{-(t-s)\nu}K^{n+1}_{1}(t,s)(1+s)\delta_{i}\|\rho[h^{n+1}]\|_{\mathcal{F}^{\lambda'_{n}s+\mu'_{n}}}ds,
$$
where $K^{n}_{0}(t,s)=e^{-\pi(\lambda'_{i}-\lambda'_{n})(t-s)},$
and
$$K^{n+1}_{1}(t,s)=e^{-2\pi(\mu'_{i}-\mu'_{n})|l|}e^{-\pi(\lambda'_{i}-\lambda'_{n})|k(t-s)+ls|}
  e^{-2\pi(\frac{\lambda'_{n}}{2}(s-s')+\mu''_{n}-\mu'_{n})|k-l|}.$$
\end{cor}
\begin{prop}($\textbf{E}$rror term $\mathcal{R}_{0}$)
 $$\|\mathcal{R}_{0}(t,\cdot)\|_{\mathcal{F}^{\lambda'_{n} t+\mu'_{n}}}
 \leq C\bigg(C'_{0}+\sum^{n}_{i=1}\delta_{i}\bigg)\bigg(\sum^{n}_{i=1}\frac{\delta_{i}}{(\lambda_{i}-\lambda'_{i})^{5}}\bigg)
\int^{t}_{0}\|\rho[h^{n+1}]\|_{\mathcal{F}^{\lambda'_{n}s+\mu'_{n}}}\frac{ds}{(1+s)^{2}}.$$
 \end{prop}
 \begin{prop}($\textbf{E}$rror term $\tilde{\mathcal{R}}_{0}$)
$$
\|\tilde{\mathcal{R}}_{0}(t,\cdot)\|_{\mathcal{F}^{\lambda'_{n}s+\mu'_{n}}}
\leq\bigg(C^{4}_{\omega}\bigg(C'_{0}+\sum^{n}_{i=1}\delta_{i}\bigg)\bigg(\sum^{n}_{j=1}\frac{\delta_{j}}{2\pi(\lambda_{j}-\lambda'_{j})^{6}}\bigg)
+\sum^{n}_{i=1}\delta_{i}\bigg)
\int^{t}_{0}\|\rho\|_{\mathcal{F}^{\lambda'_{n}s+\mu'_{n}}}\frac{1}{(1+s)^{2}}ds$$
$$=\int^{t}_{0}\tilde{K}_{1}^{n+1}\|\rho\|_{\mathcal{F}^{\lambda'_{n}s+\mu'_{n}}}\frac{1}{(1+s)^{2}}ds.$$
\end{prop}

   \begin{center}
\item\subsection{  The proof of main theorem}
%{\bf\large 1. \quad Introduction }
\end{center}

%$$K^{n}_{1,2}(t,\tau)=\sup_{k_{1,2},l_{1,2}} e^{-2\pi(\mu_{i} -\mu_{n+1})|l_{1,2}|}$$
 %$$\cdot e^{-2\pi(\lambda_{i}-\lambda_{n+1})|k_{1,2}\frac{\sin\Omega(t-\tau)}
%{\Omega}+l_{1,2}\frac{\sin\Omega\tau}{\Omega}|}e^{-2\pi(\frac{\lambda_{n+1}}{2}(|\frac{\sin\Omega\tau}{\Omega}|
%-|\frac{\sin\Omega\tau'}{\Omega}|)+\mu'_{n+1}-\mu_{n+1})|k_{1,2}-l_{1,2}|},$$
  %$$K^{n}_{3}(t,\tau)=\sup_{k_{3},l_{3}}e^{-2\pi(\mu_{i}-\mu_{n+1})|l_{3}|}e^{-\pi(\lambda_{i}-\lambda_{n+1})|k_{3}(t-\tau)+l_{3}\tau|}
 % e^{-2\pi(\frac{\lambda_{n+1}}{2}(\tau-\tau')+\mu'_{n+1}-\mu_{n+1})|k_{3}-l_{3}|}.$$
$$$$
$\mathbf{Step }$ $\mathbf{2}.$
Note from the definition  of $\delta_{n+1}$ in (7.17),  more smaller $\nu$  is, more lager the coefficient of  $\delta_{n}^{2}$  is. Therefore, without loss of generality, we assume that  $\nu$ is small enough,  up to slightly lowering $\lambda_{1},$ we may choose all parameters in such a way that
$\lambda_{k},\lambda'_{k}\rightarrow\lambda_{\infty}>\underline{\lambda}\quad \textmd{and}\quad\mu_{k},
\mu'_{k}\rightarrow\mu_{\infty}>\underline{\mu},\quad \textmd{as}
\quad k\rightarrow\infty;$
then we pick up $B>0$ such that
$\mu_{\infty}-\lambda_{\infty}(1+B)B\geq\mu'_{\infty}>\underline{\mu},$
and we let $b(t)=\frac{B}{1+t}.$
From the iteration, we have,  for all $k\geq2,$
\begin{align}
\sup_{0\leq s\leq t}\|h^{k}_{s}\circ\Omega^{k-1}_{t,s}\|_{\mathcal{Z}^{\lambda_{\infty}(1+b),\mu_{\infty};1}_{t-\frac{bt}{1+b}}}\leq\delta_{k},
\end{align}
where $\sum^{\infty}_{k=2}\delta_{k}\leq C\delta.$
Choosing $t=s$ in (3.15) yields
$\sup_{0\leq s\leq t}\|h^{k}_{s}\|_{\mathcal{Z}^{\lambda_{\infty}(1+B),\mu_{\infty};1}_{t-\frac{Bt}{1+B+t}}}\leq\delta_{k}.$
This implies that
$\sup_{ t\geq0}\|h^{k}_{t}\|_{\mathcal{Z}^{\lambda_{\infty}(1+B),\mu_{\infty}-\lambda_{\infty}(1+B)B;1}_{t}}\leq\delta_{k}.$
In particular, we have a uniform estimate on $h^{k}_{t}$ in $\mathcal{Z}^{\lambda_{\infty},\mu'_{\infty};1}_{t}.$ Summing up over $ k $ yields for
$f=f^{0}+\sum^{\infty}_{k=1}h^{k},$ the estimate
\begin{align}
\sup_{t\geq0}\|f(t,\cdot)-f^{0}\|_{\mathcal{Z}^{\lambda_{\infty},\mu'_{\infty};1}_{t}}\leq C\delta.
\end{align}
From (viiii) of Proposition 1.3, we can deduce from (3.21) that
\begin{align}
\sup_{t\geq0}\|f(t,\cdot)-f^{0}\|_{\mathcal{Y}^{\underline{\lambda},\underline{\mu}}_{t}}\leq C\delta.
\end{align}
Moreover, $\rho=\int_{\mathbb{R}^{3}}fdv$ satisfies similarly
$\sup_{t\geq0}\|\rho(t,\cdot)\|_{\mathcal{F}^{\lambda_{\infty}t+\mu_{\infty}}}\leq C\delta.$
It follows that $|\hat{\rho}(t,k)|\leq C\delta $
$e^{-2\pi\lambda_{\infty}|k|t}e^{-2\pi\mu_{\infty}|k|}$ for any $k\neq 0.$ On the one hand, by Sobolev embedding,
we deduce that for any $r\in\mathbb{N},$
$$\|\rho(t,\cdot)-\langle\rho\rangle\|_{C^{r}(\mathbb{T}^{3})}\leq C_{r}\delta e^{-2\pi\lambda't};$$
on the other hand, multiplying $\hat{\rho}$ by the Fourier transform of $W,$  we see that the electric  field
$E$ satisfies
\begin{align}
|\hat{E}(t,k)|\leq C\delta e^{-2\pi\lambda'|k|t}e^{-2\pi\mu'|k|};
\end{align}
for some $\lambda_{0}>\lambda'>\underline{\lambda},\mu_{0}>\mu'>\underline{\mu}.$

Now, from (3.15), we have, for any $(k,\eta)\in\mathbb{Z}^{3}\times\mathbb{R}^{3}$ and any $t\geq0,$
\begin{align}
|\hat{f}(t,k,\eta+kt)-\hat{f}^{0}(\eta)|\leq C\delta e^{-2\pi\mu'|k|}e^{-2\pi\lambda'|\eta|},
\end{align}
this finishes the proof of Theorem 0.1.

         $$$$
         \begin{center}
\item\section{ Dynamics of the particles' trajectory}
%{\bf\large 1. \quad Introduction }
\end{center}

Because the proof of Proposition 3.3 can be found in [26,28] here we sketch the key steps in the proof.
          To prove Proposition 3.3,  the idea is to use the classical Picard iteration, we only
          need to  consider the following  equations
           \begin{align}
   \left\{\begin{array}{l}
   \frac{d}{dt}\delta X^{n+1}_{t,s}(x,v)=\delta V^{n+1}_{t,s}(x,v),\\
  \frac{d}{dt}\delta V^{n+1}_{t,s}(x,v)=E[f^{n}](t,\delta X^{n}_{t,s}(x,v)+X^{0}_{t,s}(x,v)),\\
 \delta X^{n+1}_{s,s}(x,v)=0,\delta V^{n+1}_{s,s}(x,v)=0.
  \end{array}\right.
           \end{align}
            It is easy to check that $$\Omega^{n+1}_{t,s}-Id\triangleq(\delta X^{n+1}_{t,s},\delta V^{n+1}_{t,s})\circ (X^{0}_{s,t},V^{0}_{s,t})
            = ( X^{n+1}_{t,s}\circ(X^{0}_{s,t},V^{0}_{s,t})-Id, V^{n+1}_{t,s}\circ(X^{0}_{s,t},V^{0}_{s,t})-Id).$$
            Therefore, in order to estimate $( X^{n+1}_{t,s}\circ(X^{0}_{s,t},V^{0}_{s,t})-Id, V^{n+1}_{t,s}\circ(X^{0}_{s,t},V^{0}_{s,t})-Id),$
             we only need to study $(\delta X^{n+1}_{t,s},\delta V^{n+1}_{t,s})\circ (X^{0}_{s,t},V^{0}_{s,t}).$

   $$$$

Note that in the proof, in order to obtain $$\|\delta X^{n}_{t,s}\circ (X^{0}_{t,s},V^{0}_{t,s})\|
_{\mathcal{Z}^{\lambda_{n}'(1+b),\mu_{n}'}_{s-\frac{bt}{1+b}}}\leq \mathcal{C}^{n}_{1},$$ we need the following assumptions:

 If $s\geq\frac{bt}{1+b},$ then
 $$\nu'_{n}\leq\lambda'_{n}s+\mu'_{n}+\mathcal{C}^{n}_{1}\leq \lambda_{i}s+\mu_{i}-(\lambda_{i}-\lambda'_{n})s$$ as soon as
 $$\mathcal{C}^{n}_{1}\leq\frac{\lambda_{i}b(t-s)}{2}\quad\quad(I);$$
If $s\leq\frac{bt}{1+b},$ then $$\nu'_{n}\leq\lambda'_{n}bt+\mu'_{n}-\lambda'_{n}(1+b)s+\mathcal{C}^{n}_{1}
\leq\lambda'_{n}B+\mu'_{n}-(\lambda_{i}-\lambda'_{n})s+\mathcal{C}^{n}_{1}\leq \mu_{0}-(\lambda_{i}-\lambda'_{n})s$$ as soon as
$$\mathcal{C}^{n}_{1}\leq\frac{\mu_{0}-\mu'_{n}}{2}\quad\quad(II).$$ In order to the feasibility of the conditions $(I)$ and $(II),$ we only need to check that the following
assumption $(\mathbf{I})$ holds
$$2C^{1}_{\omega}\bigg(\sum^{n}_{i=1}\frac{\delta_{i}}{(2\pi(\lambda_{i}-\lambda'_{n}))^{3}}\bigg)\leq\min\bigg\{\frac{\lambda_{i}b(t-s)}{6},
\frac{\mu_{0}-\mu'_{n}}{2}\bigg\},\quad\quad (\mathbf{I})$$
since $\mathcal{C}^{n}_{1}\leq\omega^{1,2}_{0,n}=2C^{1}_{\omega}\bigg(\sum^{n}_{i=1}\frac{\delta_{i}}{(2\pi(\lambda_{i}-\lambda'_{n}))^{3}}\bigg)
\frac{\min\{\frac{1}{2}(t-s)^{2},1\}}{1+s}.$

          %$$ \leq\sum_{l}e^{2\pi|l_{12}|(\lambda'_{i}|\frac{\sin\Omega t}{\Omega}|+\mu'_{i})} e^{2\pi|l_{3}|(\lambda'_{i}t+\mu'_{i})}|\hat{E[f^{i}]}(t,l)|
            %$$
           %$$\cdot e^{\varepsilon_{i-1} e^{-\pi|l_{3}-k_{3}|(\lambda_{i-1}-\lambda'_{i})t}(1+|V^{0}_{t,\tau}|)}.$$
         % ( $\mathbf{how}\quad \mathbf{to}\quad\mathbf{get}? $ $If $ $we$ $prove $ $e^{2\pi|k_{12}-l_{12}|(\lambda'_{i}|\frac{\sin\Omega t}{\Omega}|+\mu'_{i})} e^{2\pi|k_{3}-l_{3}|(\lambda'_{i}t+\mu'_{i})}|(\Omega^{i}_{\tau,t})^{-1}-Id|^{\wedge}(t,l-k)\leq 2e^{2\pi|k_{12}-l_{12}|(\lambda'_{i}|\frac{\sin\Omega t}{\Omega}|+\mu'_{i})} e^{2\pi|k_{3}-l_{3}|(\lambda'_{i}t+\mu'_{i})}|\Omega^{i}_{\tau,t}-Id|^{\wedge}(t,l-k),$)

         % $$e^{2\pi|k_{12}|(\lambda_{i}'|\frac{\sin\Omega t}{\Omega}|+\mu_{i}')}e^{2\pi|k_{3}|(\lambda_{i}'t+\mu_{i}')}|[(\Omega^{n+1}_{t,\tau})^{-1}-Id]^{\wedge}(k,v)|\leq 2\kappa_{i},\kappa_{i}=(\kappa_{xi},\kappa_{vi}),$$
         % and $$e^{2\pi|k_{12}-l_{12}|(\lambda'_{i}|\frac{\sin\Omega t}{\Omega}|+\mu'_{i})}\cdot  e^{2\pi|k_{3}-l_{3}|(\lambda'_{i}t+\mu'_{i})}|[e^{2i\pi k\cdot(\Omega^{i}_{\tau,t}-Id)}]^{\wedge}(l-k)|$$
         % $$\leq e^{2\pi|k_{12}-l_{12}|(\lambda'_{i}|\frac{\sin\Omega t}{\Omega}|+\mu'_{i})}\cdot  e^{2\pi|k_{3}-l_{3}|(\lambda'_{i}t+\mu'_{i})}|[e^{2i\pi k\cdot(\Omega^{i}_{\tau,t}-Id)^{-1}(l-k)}]||(\Omega^{i}_{\tau,t})^{-1}-Id)|$$
       Combining (4.1) and $(\mathbf{I}),$   we can obtain the following conclusion
            $$\|\delta V^{n+1}_{t,s}\circ (X^{0}_{s,t},V^{0}_{s,t})\|_{\mathcal{Z}^{\lambda_{n}'(1+b),\mu_{n}'}
            _{s-\frac{bt}{1+b}}}\leq C\sum^{n}_{i=1}\delta_{i}\int^{t}_{s}e^{-2\pi(\lambda_{i}-\lambda'_{i})s}ds
            \leq C\sum^{n}_{i=1}\delta_{i} e^{-2\pi(\lambda_{i}-\lambda'_{i})s}
             \min\bigg\{\frac{(t-s)}{2},\frac{1}{2\pi(\lambda_{i}-\lambda'_{i})}\bigg\},$$

then we have  $$\|\delta X^{n+1}_{t,s}\circ (X^{0}_{s,t},V^{0}_{s,t})\|_{\mathcal{Z}^{\lambda_{n}'(1+b),\mu_{n}'}_{s-\frac{bt}{1+b}}}
\leq C\sum^{n}_{i=1}\delta_{i}\int^{t}_{s}e^{-2\pi(\lambda_{i}-\lambda'_{i})s}ds
\leq C\sum^{n}_{i=1}\delta_{i} e^{-2\pi(\lambda_{i}-\lambda'_{i})s}\min\bigg\{\frac{(t-s)^{2}}{2},\frac{1}{2\pi(\lambda_{i}
             -\lambda'_{i})^{2}}\bigg\}.$$
             We finish the proof of  Proposition 3.3.

            In the following we  estimate $\nabla\Omega_{t,s}^{n}-Id.$
In fact, we  write
$(\Omega_{t,s}^{n}-Id)(x,v)=(\delta X^{n}_{t,s},\delta V^{n}_{t,s})\circ (X^{0}_{s,t},V^{0}_{s,t}),$
and get by differentiation
$\nabla_{x}\Omega_{t,s}^{n+1}-(I,0)=\nabla_{x}(\delta X^{n}_{t,s}\circ (X^{0}_{s,t},V^{0}_{s,t}),
\delta V^{n}_{t,s}\circ (X^{0}_{s,t},V^{0}_{s,t})),$
$\nabla_{v}\Omega_{t,s}^{n}-(0,I)=(\nabla_{v}+(t-s)\nabla_{x})(\delta X^{n}_{t,s}\circ (X^{0}_{s,t},V^{0}_{s,t}),
\delta V^{n}_{t,s}\circ (X^{0}_{s,t},V^{0}_{s,t})).$
 \begin{align}
   \left\{\begin{array}{l}
   \frac{d}{dt}\nabla_{x}\delta X^{i}_{t,s}(x,v)=\nabla_{x}\delta V^{i}_{t,s}(x,v),\\
  \frac{d}{dt}\nabla_{x}\delta V^{i}_{t,s}(x,v)=
  \nabla_{x}E[f^{i}](t,\delta X^{i}_{t,s}(x,v) +X^{0}_{t,s}(x,v)),\\
 \delta X^{i}_{s,s}(x,v)=0,\quad\delta V^{i}_{s,s}(x,v)=0.
  \end{array}\right.
           \end{align}

           Using the same process in the proof of Proposition 3.3, we can obtain Proposition 3.5.

           To establish a control of $\Omega^{i}_{t,s}-\Omega^{n}_{t,s}$ in norm $\mathcal{Z}^{\lambda'_{n}(1+b),\mu'_{n}}_{s-\frac{bt}{1+b}},$
           we start again from the differential equation satisfied by $\delta V^{i}_{t,s}$ and $\delta V^{n}_{t,s}:$
          \begin{align}
   \left\{\begin{array}{l}
   \frac{d}{dt}(\delta X^{i}_{t,s}-\delta X^{n}_{t,s})(x,v)=\delta V^{i}_{t,s}(x,v)-\delta V^{n}_{t,s}(x,v),\\
  \frac{d}{dt}(\delta V^{i}_{t,s}-\delta V^{n}_{t,s})(x,v)=
  E[f^{i-1}](t,\delta X^{i-1}_{t,s}(x,v) +X^{0}_{t,\tau}(x,v))\\
  -E[f^{n-1}](t,\delta X^{n-1}_{t,s}(x,v) +X^{0}_{t,s}(x,v)),\\
 (\delta X^{i}_{t,s}-\delta X^{n}_{t,s})(x,v)=0,\quad (\delta V^{i}_{t,s}-\delta V^{n}_{t,s})(x,v)=0.
  \end{array}\right.
           \end{align}

           So from (4.3), $\delta V^{i}_{t,s}-\delta V^{n}_{t,s}$ satisfies the equation:
           $$\frac{d}{dt}(\delta V^{i}_{t,s}-\delta V^{n}_{t,s})(x,v)=
  E[f^{i-1}](t,\delta X^{i-1}_{t,s}(x,v) +X^{0}_{t,s}(x,v))$$
 $$-E[f^{n-1}](t,\delta X^{i-1}_{t,s}(x,v) +X^{0}_{t,s}(x,v))+E[f^{n-1}](t,\delta X^{i-1}_{t,s}(x,v) +X^{0}_{t,s}(x,v))
  -E[f^{n-1}](t,\delta X^{n-1}_{t,s}(x,v) +X^{0}_{t,s}(x,v)).$$

  Under the assumption $(\mathbf{I}),$  we can use the similar proof of Proposition 3.3 to finish Proposition 3.6.

  Let $\varepsilon$ be the small constant appearing in Lemma 1.7. If
  $$3\mathcal{C}_{1}^{i}+\mathcal{C}_{2}^{i}\leq\varepsilon,\quad \textmd{for}\quad \textmd{all}\quad i\geq1,\quad (\mathbf{II});$$  if in addition
  $$2(1+s)(1+B)(3\mathcal{C}_{1}^{i,n}+\mathcal{C}_{2}^{i,n})(s,t)\leq\max\{\lambda'_{i}-\lambda'_{n},\mu'_{i}-\mu'_{n}\}\quad (\mathbf{III})$$
  for all $i\in\{1,\ldots,n-1\}$ and all $t\geq s,$
  then  \begin{align}
   \left\{\begin{array}{l}
   \lambda'_{n}(1+b)+2\|\Omega^{n}-\Omega^{i}\|_{\mathcal{Z}^{\lambda'_{n}(1+b),\mu'_{n}}_{s-\frac{bt}{1+b}}}\leq \lambda'_{i}(1+b),\\
  \mu'_{n}+2(1+|s-\frac{bt}{1+b}|)\|\Omega^{n}-\Omega^{i}\|_{\mathcal{Z}^{\lambda'_{n}(1+b),\mu'_{n}}_{s-\frac{bt}{1+b}}}\leq\mu'_{i}.
 \end{array}\right.
           \end{align}

Then Lemma 1.7 and (4.6) yield Proposition 3.7.

As  a corollary of Proposition 3.7 and Proposition 1.8, under the assumption $(\mathbf{IV}):
$
$$4(1+s)(\mathcal{C}_{1}^{i,n}+\mathcal{C}_{2}^{i,n})\leq\min\{\lambda_{i}-\lambda'_{n},\mu_{i}-\mu'_{n}\}$$
for all $i\in\{1,\ldots,n\}$ and all $s\in[0,t],$
we have
\begin{cor} under the assumption (3.8), we have
$$\|h^{i}_{s}\circ\Omega^{n}_{t,s}\|
_{\mathcal{Z}_{s-\frac{bt}{1+b}}^{\lambda'_{n}(1+b),\mu'_{n};1}}\leq\delta_{i},$$
$$\|((\nabla_{v}+s\nabla_{x})h^{i}_{s})\circ\Omega^{n}_{t,s}\|
_{\mathcal{Z}_{s-\frac{bt}{1+b}}^{\lambda'_{n}(1+b),\mu'_{n};1}}\leq\delta_{i}.$$
\end{cor}

$$$$
\begin{center}
\item\section{The estimates of main terms}
%{\bf\large 1. \quad Introduction }
\end{center}

In order to estimate  that term $I,$ we have to make good understanding of plasma echo. First of all, one of the key steps is that  we need to translate the physical phenomenon into the mathematical language. In the following we give the corresponding mathematical analysis.
\begin{center}
\item\subsection{Plasma echoes: Mathematical expression}
%{\bf\large 1. \quad Introduction }
\end{center}

From the above physical point of view, under the assumption of the stability condition, we have known that,  echoes occurring
 at distinct frequencies are asymptotically well separated. In the following, through  complicate computation, we give a detailed description by
  using mathematical tool. And the proof is simple, so we omit.

 \begin{thm} Let $\lambda,\bar{\lambda},\mu,\bar{\mu},\mu'$ be such that
$2\lambda\geq\bar{\lambda}>\lambda>0, \bar{\mu}\geq\mu'>\mu>0,$ and let $b=b(t,s)>0,$
  $ R=R(t,x),G=G(t,x,v)$ and assume $\widehat{G}(t,0,v)=0,$ we have, if
$$\sigma(t,x,v)=\int^{t}_{0}R(s,x+(t-s)v)G(s,x+(t-s)v,v)ds,$$
 then
\begin{align}
 & \|\sigma(t,\cdot)\|_{\mathcal{Z}^{\lambda,\mu;1}_{t}}\leq\int^{t}_{0} K(t,s)\|R\|_{\mathcal{F}^{\lambda s+\mu'-\lambda b(t-s)}}\frac{\|G\|_{\mathcal{Z}^{\bar{\lambda}(1+b),\bar{\mu};1}_{s-\frac{bt}{1+b}}}}{1+s}ds,
\end{align}
where $K(t,s)=(1+s)\sup_{k,l\in\mathbb{Z}^{3}_{\ast}}e^{-\pi(\bar{\mu}-\mu)|l|}
e^{-\pi(\bar{\lambda}-\lambda)|k(t-s)+ls|}e^{-2\pi[\mu'-\mu+\lambda b(t-s)]|k-l|}.$

\end{thm}

Now we try to explain this above theorem from the two aspects: mathematical and physical,respectively. First, the inequality (5.1) is vital for the Vlasov-Poisson equations (0.11). Now we assume that the function $G(t,x,v)$ is known and is good enough, then  in some sense, if the kernel is ``good", (5.1) is
considered as    a  inverse ``H$\ddot{\textmd{o}}$lder" inequality on the function $R(t,x);$ on the other hand,  from the energy viewpoint, for (0.11), (5.1) can be also regarded as a ``monotone" energy formula.  However, in order to prove the inverse H$\ddot{\textmd{o}}$lder inequality or the ``monotone" energy formula holds, we have to check whether  the kernel $K_{t,s}$  is good or not.

Indeed, let $\mu'-\mu=\sigma,$ assuming that $b=\frac{B}{1+t}$ with $ B $ so small that $(\mu'-\mu)-\lambda b(t-s)\geq\frac{\sigma }{2},$ then $K$ is bounded by
$$K^{\alpha}(t,s)=(1+s)\sup_{k,l\in\mathbb{Z}^{3}_{\ast}}e^{-\pi\alpha|l|}
e^{-\pi\alpha|k(t-s)+ls|}e^{-2\pi\alpha|k-l|},$$
where $\alpha=\frac{1}{2}\min\{\bar{\lambda}-\lambda,\bar{\mu}-\mu,\sigma\}.$ Through  simple computation, it is easy to find that when $s\leq\frac{t}{2},$ $K^{\alpha}(t,s)$ is ``good"; however, when $\frac{1}{2}t< s\leq t,$  whenever $s/t$ is a rational number distinct from 0 or 1, there are $k,l\in\mathbb{Z}^{3}$ such that $|k(t-s)+ls|=0,$ $K^{\alpha}(t,s)$ is ``bad", that is, as $t\rightarrow\infty,$ $K^{\alpha}_{t,s}$ maybe cannot be controlled.
From the physical point of view, one can consider $l,k-l$ as frequencies of two different waves, and start a wave at frequency $l$ at time 0, and force it at time $ s $ by a wave of frequency $k-l,$ a strong response is obtained at time $ t $ and frequency $ k$ such that $k(t-s)+ls=0.$ And the corresponding strong response is called plasma echo in plasma physics. It is worthy mentioned that the condition $x\in\mathbb{T}^{3}$ guarantees the asymptotically well separated behavior of echoes occurring at distinct frequencies, which was discovered by Mouhot and Villani. The detailed computation is found in the following section 7.2 ( also see the paper [26,28] ).

%$$\|(RG)(s,\cdot)\|_{\mathcal{Z}^{\lambda,\mu;1}_{s}}\leq \sup_{k,l\in\mathbb{Z}^{3},k\neq l} e^{2\pi\lambda b|k(t-s)+ls|}e^{-2\pi(\bar{\lambda}-\lambda)|k-l|s}\|R(s,\cdot)\|_{\mathcal{F}^{\bar{\lambda}s+\mu}}
%\|G(s,\cdot)\|_{\mathcal{Z}^{\lambda(1-b),\mu;1}_{s}},$$

%Note that
%\begin{align}
%|k(t-s)+ls|\sim
%\left\{
%\begin{array}
%c|k-l|t,\quad if \quad |k|\gg|l|\quad or\quad |l|\gg|k| \\
%|l|t,\quad \quad if\quad\quad \quad |k|\sim|l| ,\\
%\end{array}
%\right.
%\end{align}

\begin{center}
\item\subsection{ Estimates of main terms}
%{\bf\large 1. \quad Introduction }
\end{center}

% $ \mathbf{ Note}$ that,
  % \begin{align}
  % &\int_{\mathbb{T}^{3}}\int_{\mathbb{R}^{3}}-(\mathcal{E}^{n+1}_{t,s}\cdot G_{t,s}^{n})(s, X^{0}_{t,s}(x,v),V^{0}_{t,s}(x,v))dxdv\neq0\notag\\
   %&\int_{\mathbb{T}^{3}}\int_{\mathbb{R}^{3}}-(F^{n+1}_{t,s}\cdot G_{t,s}^{n,v})(s, X^{0}_{t,s}(x,v),V^{0}_{t,s}(x,v)))dxdv\neq0\notag\\
   %&\int_{\mathbb{T}^{3}}\int_{\mathbb{R}^{3}}-(\mathcal{E}^{n}_{t,s}\cdot H^{n}_{t,s})(s, X^{0}_{t,s}(x,v),V^{0}_{t,s}(x,v)))dxdv\neq0\notag\\
 % &\int_{\mathbb{T}^{3}}\int_{\mathbb{R}^{3}}-(F^{n}_{t,s}\cdot H_{t,s}^{n,v})(s, X^{0}_{t,s}(x,v),V^{0}_{t,s}(x,v)))dxdv\neq0\notag\\
  % \end{align}

   %$\mathbf{Problem:}$ how to add condition such that $\hat{\rho}^{n+1}(t,k_{\perp},0)=0$?

%Since $\hat{W}(k_{1},k_{2},0)=0,$ it is trivial to check that $\hat{\rho}^{n+1}(t,k_{\perp},0)=0.$

In the following we estimate $\bar{I}_{i}^{n+1,n}(t,x).$ Note that their zero modes vanish. For any $n\geq i\geq1,$
$$\widehat{\bar{I}_{i}^{n+1,n}}(t,k)
=\int^{t}_{0}\int_{\mathbb{T}^{3}}\int_{\mathbb{R}^{3}}e^{-(t-s)\nu}e^{-2\pi ik\cdot x} \bigg(E[h^{n+1}]
\cdot (\nabla_{v}(h^{i}_{s}\circ\Omega^{i-1}_{t,\tau})\bigg)(s,x-v(t-s),v)dvdxds,$$
\begin{align}
&|\widehat{\bar{I}_{i}^{n+1,n}}(t,k)|\leq\int^{t}_{0}\bigg(\sum_{l\in\mathbb{Z}^{3}_{\ast}}\bigg|\int_{\mathbb{R}^{3}}e^{-(t-s)\nu}e^{-2\pi ik\cdot(v(t-s))} (\widehat{\nabla_{v}(h^{i}_{s}\circ\Omega^{i-1}_{t,s}))}(s,l,v)dv\bigg||\widehat{E[h^{n+1}]}(s,k-l)|\bigg)ds\notag\\
&=\int^{t}_{0}\bigg(\sum_{l\in\mathbb{Z}^{3}_{\ast}}\bigg|\int_{\mathbb{R}}e^{-(t-s)\nu}e^{-2\pi ik\cdot(v(t-s))} (\widehat{\nabla_{v}(h^{i}_{\tau}\circ\Omega^{i-1}_{t,s}))}(s,l,v)dv\bigg||\widehat{E[h^{n+1}]}(s,k-l)|\bigg)ds\notag\\
\end{align}

From (5.1) of Theorem 5.1 and (5.2), we can get Proposition 3.9.

$$$$

\begin{center}
\item\section{ Estimates of error terms}
%{\bf\large 1. \quad Introduction }
\end{center}

In the following we estimate one of the error terms $\mathcal{R}_{0}.$

 Recall
$$\mathcal{R}_{0}(t,x)=\int^{t}_{0}\int_{\mathbb{R}^{3}}e^{-(t-s)\nu}\bigg(\bigg(E[h^{n+1}]\circ\Omega^{n}_{t,s}(x,v)-E[h^{n+1}]\bigg)\cdot G_{t,s}^{n}\bigg)(s, X^{0}_{t,s}(x,v),V^{0}_{t,s}(x,v)))dvds.$$
First,
$$\|\mathcal{R}_{0}(t,\cdot)\|_{\mathcal{F}^{\lambda'_{n} t+\mu'_{n}}}$$
$$\leq\int^{t}_{0}e^{-(t-s)\nu}\|E[h^{n+1}]\circ\Omega^{n}_{t,s}(x,v)-E[h^{n+1}]\|_{\mathcal{Z}_{s-\frac{bt}{1+b}}^{\lambda'_{n}(1+b),\mu'_{n}}}
\|G_{t,s}^{n}\|_{\mathcal{Z}_{s-\frac{bt}{1+b}}^{\lambda'_{n}(1+b),\mu'_{n};1}}ds.$$
Next,
\begin{align}
&\|(E[h^{n+1}]\circ\Omega^{n}_{t,s}(x,v)-E[h^{n+1}]\|_{\mathcal{Z}_{s-\frac{bt}{1+b}}^{\lambda'_{n}(1+b),\mu'_{n}}}\notag\\
&\leq\|\Omega^{n}_{t,s}-Id\|_{\mathcal{Z}_{s-\frac{bt}{1+b}}^{\lambda'_{n}(1+b),\mu'_{n}}}
\int^{1}_{0}\|\nabla E[h^{n+1}]((1-\theta)Id+\theta\Omega^{n}_{t,s})\|_{\mathcal{Z}_{s-\frac{bt}{1+b}}^{\lambda'_{n}(1+b),\mu'_{n}}}d\theta\notag\\
&\leq\|\nabla E[h^{n+1}]\|_{\mathcal{F}^{\nu'_{n}}}
\|\Omega^{n}_{t,s}-Id\|_{\mathcal{Z}_{s-\frac{bt}{1+b}}^{\lambda'_{n}(1+b),\mu'_{n}}},\notag\\
\end{align}
where
$\nu'_{n}=\lambda'_{n}(1+b)\bigg|s-\frac{bt}{1+b}\bigg|+\mu'_{n}+\|\Omega^{n}X_{t,s}-Id\|_{\mathcal{Z}_{s-\frac{bt}{1+b}}
^{\lambda'_{n}(1+b),\mu'_{n}}}.$

Here we only focus on the case $s\geq\frac{bt}{1+b},$ then we need to $ \mathbf{show }$
$\|\nabla E[h^{n+1}]\|_{\mathcal{F}^{\nu'_{n}}}\leq\|\rho[h^{n+1}]\|_{\mathcal{F}^{\lambda'_{n}s+\mu'_{n}}}.$
For that, we need to use the fact $E[h^{n+1}]=\nabla_{x}W(x)\ast\rho([h^{n+1}])$  and prove
$\nu'_{n}<\lambda'_{n}s+\mu'_{n}-\iota,$
for some constant $\iota>0$ sufficiently small.

Indeed,$$\nu'_{n}\leq\lambda'_{n}s+\mu'_{n}-\lambda'_{n}b(t-s)+C\sum^{n}_{i=1}
\delta_{i}e^{-\pi|k|(\lambda_{i}-(\lambda'_{i}))t}\cdot\min\bigg\{\frac{(t-s)^{2}}{2},\frac{1}{2\pi(\lambda_{i}-\lambda'_{i})^{2}}\bigg\}$$
$$\leq\lambda'_{n}s+\mu'_{n}-\lambda'_{n}\frac{B(t-s)}{1+t}+C\bigg(\sum^{n}_{i=1}\frac{\delta_{i}}{(\lambda_{i}
-\lambda'_{i})^{3}}\bigg)
\frac{\min\{t-s,1\}}{1+s}.$$
Note that
$\frac{\min\{t-s,1\}}{1+s}\leq3\frac{t-s}{1+t}.$
In the following we also need to show that
$$
C\sum^{n}_{i=1}\frac{\delta_{i}}{(\lambda_{i}-\lambda'_{i})^{3}}\leq\frac{\lambda^{\ast}B}{3}-\iota, \quad (\mathbf{VI})
$$

From Proposition 3.3, then
$$\|E[h^{n+1}]\circ\Omega^{n}_{t,s}(x,v)-E[h^{n+1}]\|_{\mathcal{Z}_{s-\frac{bt}{1+b}}^{\lambda'_{n}(1+b),\mu'_{n}}}
\leq C\bigg(\sum^{n}_{i=1}\frac{\delta_{i}}{(\lambda_{i}-\lambda'_{i})^{5}}\bigg)
\frac{1}{(1+s)^{3}}\|\rho[h^{n+1}]\|_{\mathcal{F}^{\lambda'_{n}s+\mu'_{n}}}$$

 Since $G_{t,s}^{n}=\nabla_{v}f^{n}\circ\Omega^{n}_{t,s},$
 $$\|G_{t,s}^{n}\|_{\mathcal{Z}_{s-\frac{bt}{1+b}}^{\lambda'_{n}(1+b),\mu'_{n};1}}
 \leq\|(\nabla_{v}f^{0})\circ\Omega^{n}_{t,s}\|
 _{\mathcal{Z}_{s-\frac{bt}{1+b}}^{\lambda'_{n}(1+b),\mu'_{n};1}}$$
$$+\sum^{n}_{i=1}\|\nabla_{v}h_{s}^{i}\circ\Omega^{n}_{t,s}\|_{\mathcal{Z}_{s-\frac{bt}{1+b}}^{\lambda'_{n}(1+b),\mu'_{n};1}}$$
$$\leq C'_{0}+\bigg(\sum^{n}_{i=1}\delta_{i}\bigg)(1+s).$$

  We can conclude
$$\|\mathcal{R}_{0}(t,\cdot)\|_{\mathcal{F}^{\lambda'_{n} t+\mu'_{n}}}\leq C\bigg(C'_{0}+\sum^{n}_{i=1}\delta_{i}\bigg)\bigg(\sum^{n}_{i=1}\frac{\delta_{i}}{(\lambda_{i}-\lambda'_{i})^{5}}\bigg)
\int^{t}_{0}e^{-(t-s)\nu}\|\rho[h^{n+1}]\|_{\mathcal{F}^{\lambda'_{n}s+\mu'_{n}}}\frac{ds}{(1+s)^{2}}.$$

In order to finish the control of $\tilde{\mathcal{R}}_{0},$
 we still need the estimate of $\|G_{t,s}^{n}-\bar{G}_{t,s}^{n}\|_{\mathcal{Z}_{s-\frac{bt}{1+b}}^{\lambda'_{n}
 (1+b),\mu'_{n};1}}.$

  In fact,
$$\|G_{t,s}^{n}-\bar{G}_{t,s}^{n}\|_{\mathcal{Z}_{s-\frac{bt}{1+b}}^{\lambda'_{n}(1+b),\mu'_{n};1}}$$
$$\leq
\|\nabla_{v}f^{0}\circ\Omega^{n}_{t,s}-\nabla_{v}f^{0}\|_{\mathcal{Z}_{s-\frac{bt}{1+b}}^{\lambda'_{n}(1+b),\mu'_{n};1}}$$
$$+\sum^{n}_{i=1}\|\nabla_{v}h^{i}\circ\Omega^{n}_{t,s}-\nabla_{v}h^{i}\circ\Omega^{i-1}_{t,s}\|
_{\mathcal{Z}_{s-\frac{bt}{1+b}}^{\lambda'_{n}(1+b),\mu'_{n};1}}$$
$$+\sum^{n}_{i=1}\|(\nabla_{v}h^{i})\circ\Omega^{i-1}_{t,s}-\nabla_{v}(h^{i}\circ\Omega^{i-1}_{t,s})\|
_{\mathcal{Z}_{s-\frac{bt}{1+b}}^{\lambda'_{n}(1+b),\mu'_{n};1}}.$$

  Now on the one hand, we treat the second term
 $$\sum^{n}_{i=1}\|(\nabla_{v}h^{i})\circ\Omega^{n}_{t,s}-(\nabla_{v}h^{i})\circ\Omega^{i-1}_{t,s}\|_{\mathcal{Z}_{s-\frac{bt}{1+b}}
 ^{\lambda'_{n}(1+b),\mu'_{n};1}}$$
 $$\leq\int^{1}_{0}\|\nabla\nabla_{v}h^{i}_{s}((1-\theta)\Omega^{n}_{t,s}+\theta\Omega^{i-1}_{t,s})\|_{\mathcal{Z}
 _{s-\frac{bt}{1+b}}^{\lambda'_{n}(1+b),\mu'_{n};1}}
 \cdot\|\Omega^{n}_{t,s}-\Omega^{i-1}_{t,s}\|
 _{\mathcal{Z}_{s-\frac{bt}{1+b}}^{\lambda'_{n}(1+b),\mu'_{n}}}d\theta$$
 $$\leq2\|\nabla\nabla_{v}h^{i}_{s}\circ\Omega^{i-1}_{t,s}\|_{\mathcal{Z}
 _{s-\frac{bt}{1+b}}^{\lambda'_{i}(1+b),\mu'_{i};1}}\|\Omega^{n}_{t,s}-\Omega^{i-1}_{t,s}\|
 _{\mathcal{Z}_{s-\frac{bt}{1+b}}^{\lambda'_{n}(1+b),\mu'_{n}}}$$
 $$\leq4C\delta_{i}\bigg(\sum^{n}_{j=i}\frac{\delta_{j}}{(\lambda_{j}-\lambda'_{j})^{6}}\bigg)\frac{1}{(1+s)^{2}},$$
 where we need to prove
 $$\|\Omega^{n}X_{t,s}-\Omega^{i-1}X_{t,s}\|
 _{\mathcal{Z}_{s-\frac{bt}{1+b}}^{\lambda'_{n}(1+b),\mu'_{n}}}\leq2\mathcal{R}^{i-1,n}_{2}(t,s),$$
 $$\|\Omega^{n}V_{t,s}-\Omega^{i-1}V_{t,s}\|
 _{\mathcal{Z}_{s-\frac{bt}{1+b}}^{\lambda'_{n}(1+b),\mu'_{n}}}
 \leq\mathcal{R}^{i-1,n}_{1}(t,s)+\mathcal{R}^{i-1,n}_{2}(t,s)$$
 with $$\mathcal{R}^{i-1,n}_{1}(t,s)=\bigg(\sum^{n}_{j=i}\frac{\delta_{j}e^{-2\pi(\lambda_{j}-\lambda'_{j})s}}{2\pi(\lambda_{j}-\lambda'_{j})}\bigg)
 \min\{t-s,1\},\mathcal{R}^{i-1,n}_{2}(t,s)=\bigg(\sum^{n}_{j=i}\frac{\delta_{j}e^{-2\pi(\lambda_{j}-\lambda'_{j})s}}
 {(2\pi(\lambda_{j}-\lambda'_{j}))^{2}}\bigg)
 \min\bigg\{\frac{(t-s)^{2}}{2},1\bigg\}.$$
On the other hand, by the  induction hypothesis, since  $\mathcal{Z}^{\lambda,\mu}_{s}$ norm is increasing as a function of $\lambda$ and $\mu,$ then
$$\sum^{n}_{i=1}\|(\nabla_{v}h^{i}_{s})\circ\Omega^{i-1}_{t,s}-\nabla_{v}(h^{i}_{s}\circ\Omega^{i-1}_{t,s})\|
_{\mathcal{Z}_{s-\frac{bt}{1+b}}^{\lambda'_{n}(1+b),\mu'_{n}}}
\leq\bigg(\sum^{n}_{i=1}\delta_{i}\bigg)\frac{1}{(1+s)^{2}}.$$
%The proof of $\|\Omega^{n}_{t,\tau}-\Omega^{i-1}_{t,\tau}\|
% _{\mathcal{Z}^{\lambda_{i}(1+b),
%\mu_{i}}_{\tau-\frac{bt}{1+b}}}$ is similar to the proof of $Proposition$ 3.2, so here we omit the details  of the proof.
%$$\|\Omega^{n}_{t,\tau}-\Omega^{i-1}_{t,\tau}\|
% _{\mathcal{Z}^{\lambda_{i}(1+b),
%\mu_{i}}_{\tau-\frac{bt}{1+b}}}$$
%$$\leq\|\Omega^{i-1}_{t,\tau}\| _{\mathcal{Z}^{\lambda_{i}(1+b),
%\mu_{i}}_{\tau-\frac{bt}{1+b}}}\|(\Omega^{i-1}_{t,\tau})^{-1}\Omega^{n}_{t,\tau}-Id\|
% _{\mathcal{Z}^{\lambda_{i}(1+b),
%\mu_{i}}_{\tau-\frac{bt}{1+b}}}$$
%$$\leq\sum^{n}_{j=i}\frac{\delta_{j}}{(2\pi(\lambda_{j}-\lambda'_{j})^{6})}\frac{1}{(1+\tau)^{5}}$$
%$\mathbf{Assumption:}$ $$\|\nabla\nabla_{v}h^{i}_{\tau}(\Omega^{i-1}_{t,\tau})\|_{\mathcal{Z}^{\lambda_{i}(1+b),
%\mu_{i};1}_{\tau-\frac{bt}{1+b}}}\leq (1+\tau)^{2}\delta_{i}.$$
%So we have
%$$\|\nabla\nabla_{v}h^{i}_{\tau}((1-\theta)\Omega^{n}_{t,\tau}+\theta\Omega^{i-1}_{t,\tau})\|_{\mathcal{Z}^{\lambda_{i}(1+b),
%\mu_{i};1}_{\tau-\frac{bt}{1+b}}}\leq2\|\nabla\nabla_{v}h^{i}_{\tau}(\Omega^{i-1}_{t,\tau})\|_{\mathcal{Z}^{\lambda_{i}(1+b),
%\mu_{i};1}_{\tau-\frac{bt}{1+b}}}\leq (1+\tau)^{2}\delta_{i}.$$

So we have
\begin{align}
&\|\tilde{\mathcal{R}}_{0}(t,\cdot)\|_{\mathcal{F}^{\lambda'_{n} s+\mu'_{n}}}
\leq\bigg(C^{4}_{\omega}\bigg(C'_{0}+\sum^{n}_{i=1}\delta_{i}\bigg)\bigg(\sum^{n}_{j=1}\frac{\delta_{j}}{2\pi(\lambda_{j}-\lambda'_{j})^{6}}\bigg)
+\sum^{n}_{i=1}\delta_{i}\bigg)
\int^{t}_{0}e^{-(t-s)\nu}\|\rho\|_{\mathcal{F}^{\lambda'_{n} s+\mu'_{n}}}\frac{1}{(1+s)^{2}}ds\notag\\
&=\int^{t}_{0}e^{-(t-s)\nu}\tilde{K}_{1}^{n+1}\|\rho[h^{n+1}]\|_{\mathcal{F}^{\lambda'_{n} s+\mu'_{n}}}\frac{1}{(1+s)^{2}}ds.\notag\\
\end{align}
Up to now we finish the estimates of error terms.
%In the following we continue to estimate $\tilde{\mathcal{R}}_{0}$ to finish the estimate of $II^{n+1,n}.$ Now we recall
%$$\tilde{\mathcal{R}}_{0}=\int^{t}_{0}\int_{\mathbb{R}^{3}}(B[h^{n+1}]\cdot (G_{t,s}^{n,v}-\bar{G}_{t,s}^{n,v}))(s, X^{0}_{t,s}(x,v),V^{0}_{t,s}(x,v))dvds,$$
$$$$

\begin{center}
\item\section{ Iteration}
%{\bf\large 1. \quad Introduction }
\end{center}

Now let us first deal with the source term
\begin{align}
 &\widehat{II^{n,n}}(t,k)=-\int^{t}_{0}\int_{\mathbb{T}^{3}}\int_{\mathbb{R}^{3}}e^{-(t-s)\nu}e^{-2\pi ik\cdot x}(\mathcal{E}^{n}_{t,s}\cdot H^{n}_{t,s})(s, X^{0}_{t,s}(x,v),V^{0}_{t,s}(x,v))dvdxds,\notag\\
          % &=-\int^{t}_{0}\int_{\mathbb{T}^{3}}\int_{\mathbb{R}^{3}}\exp(-2\pi ik\cdot(x-v\circ\mathcal{I}_{0}))(\mathcal{E}^{n}_{t,s}\cdot H^{n}_{t,s})(s, x,V^{0}_{t,s}(x,v))dvdxds\notag\\
            % &-\int^{t}_{0}\int_{\mathbb{T}^{3}}\int_{\mathbb{R}^{3}}e^{-2\pi ik\cdot (x-v\circ\mathcal{I}_{0})}(F^{n}_{t,s}\cdot H^{n,v}_{t,s})(s, x,V^{0}_{t,s}(x,v))dvdxds\notag\\
 \end{align}
 %where $X^{0}_{t,s}(x,v)=x+v\circ\mathcal{I}_{0},V^{0}_{t,s}(x,v)=v\circ\mathcal{I}_{1}$
 %$$(v\circ\mathcal{I}_{0})(t,\tau)=(-\frac{v_{\perp}}{\Omega}[\sin(\theta+\Omega(t-\tau))-\sin\theta],\frac{v_{\perp}}{\Omega}[\cos(\theta+\Omega(t-\tau))-\cos\theta],-v_{z}(t-\tau))$$
%$$(v\circ\mathcal{I}_{1})(t,\tau)=(v_{\perp}\cos(\theta+\Omega(t-\tau)),v_{\perp}\sin(\theta+\Omega(t-\tau)),v_{z}),$$
then
%\begin{align}
 %&\hat{III}^{n,n}(t,k)+\hat{IV}^{n,n}(t,k)\notag\\
% &=-\int^{t}_{0}\sum_{l}\int_{\mathbb{R}^{3}}\exp(-2\pi ik\cdot(-v\circ\mathcal{I}_{0}))\mathcal{E}^{n}_{t,s}(s,k-l) H^{n}_{t,s}(s, l,V^{0}_{t,s}(x,v))dvds\notag\\
             %&-\int^{t}_{0}\sum_{l}\int_{\mathbb{R}^{3}}e^{-2\pi ik\cdot (-v\circ\mathcal{I}_{0})}(F^{n}_{t,s}(s,k-l) H^{n,v}_{t,s}(s, l,V^{0}_{t,s}(x,v))dvds\notag\\
             %&=-\int^{t}_{0}\sum_{l}\int_{\mathbb{R}^{3}}\exp(-2\pi ik\cdot(-v\circ\mathcal{I}_{0}))\mathcal{E}^{n}_{t,s}(s,k-l) H^{n}_{t,s}(s, l,V^{0}_{t,s}(x,v))dvds\notag\\
 %\end{align}

 $$\|II(t,\cdot)\|_{\mathcal{F}^{\lambda'_{n} t+\mu'_{n}}}\leq\int^{t}_{0}e^{-(t-s)\nu}\|\mathcal{E}^{n}_{t,s}\|_{\mathcal{Z}^{\lambda'_{n}(1+b),
\mu'_{n}}_{s-\frac{bt}{1+b}}}\|H^{n}_{\tau,s}\|_{\mathcal{Z}^{\lambda'_{n}(1+b),
\mu'_{n};1}_{s-\frac{bt}{1+b}}}ds$$
\begin{align}
&\leq\int^{t}_{0}e^{-(t-s)\nu}\|\rho^{n}_{t,s}\|_{\mathcal{F}^{\nu_{n+1}}}\|H^{n}_{t,s}\|_{\mathcal{Z}^{\lambda'_{n}(1+b),
\mu'_{n};1}_{s-\frac{bt}{1+b}}}ds\leq\int^{t}_{0}e^{-(t-s)\nu}\|\rho^{n}_{t,s}\|_{\mathcal{F}^{\nu'_{n}}}(1+s)\delta_{n}ds\notag\\
&\leq\int^{t}_{0}\|\rho^{n}_{t,s}\|
_{\mathcal{F}^{\lambda'_{n} s+\mu'_{n}}}e^{-(t-s)\nu}e^{-2\pi s(\lambda_{n} -\lambda'_{n} )}(1+s)\delta_{n}d\tau
\leq\frac{C\delta^{2}_{n}}{(\lambda_{n} -\lambda'_{n}-\frac{\nu}{2\pi})^{2}},\notag\\
\end{align}
and \begin{align}
\|III^{n+1,0}(t,\cdot)\|_{\mathcal{F}^{\lambda'_{n}t+\mu'_{n}}}\leq\nu\int^{t}_{0}e^{-(t-s)\nu}\|\rho_{s}[h^{n+1}]\|_{\mathcal{F}^{\lambda'_{n}s+\mu'_{n}}}
\|f^{0}\|_{\mathcal{F}^{\lambda'_{n}s+\mu'_{n}}}ds
\end{align}

From Propositions 3.5-3.11, combining (3.10), we conclude
$$\|\rho[h^{n+1}](t,\cdot)\|_{\mathcal{F}^{\lambda'_{n} t+\mu'_{n}}}$$
$$\leq\frac{C\delta^{2}_{n}}{(\lambda_{n} -\lambda'_{n}-\frac{\nu}{2\pi})^{2}}+\int^{t}_{0}e^{-(t-s)\nu}|K^{n}_{1}(t,s)|(1+s)\sum^{n}_{i=1}\delta_{i}\| \rho[h^{n+1}]\|
_{\mathcal{F}^{\lambda'_{n}s+\mu'_{n}}}ds+2\int^{t}_{0}e^{-(t-s)\nu}|K^{n}_{0}(t,s)|$$
$$\cdot\sum^{n}_{i=1}\delta_{i}\| \rho[h^{n+1}]\|_{\mathcal{F}^{\lambda'_{n}s+\mu'_{n}}}ds+\int^{t}_{0}e^{-(t-s)\nu}
\bigg[\frac{(\tilde{K}_{0}^{n+1}+\tilde{K}_{1}^{n+1})}{(1+s)^{2}}+\nu C_{0}\bigg]\|\rho[h^{n+1}]\|_{\mathcal{F}^{\lambda'_{n}s+\mu'_{n}}}ds,$$
where $K^{n}_{0}(t,s),K^{n}_{1}(t,s)$ are defined in Proposition 3.4, and
 $$\tilde{K}_{0}^{n+1}\triangleq 2C\bigg(C_{0}+\sum^{n}_{i=1}\delta_{i}\bigg)\bigg(\sum^{n}_{i=1}\frac{\delta_{i}}{(2\pi(\lambda_{i}-\lambda'_{i}))^{5}}\bigg),$$
$$\tilde{K}_{1}^{n+1}\triangleq\bigg(C^{4}_{\omega}\bigg(C'_{0}+\sum^{n}_{i=1}\delta_{i}\bigg)\bigg(\sum^{n}_{j=1}\frac{\delta_{j}}
{2\pi(\lambda_{j}-\lambda'_{j})^{6}}\bigg)
+\sum^{n}_{i=1}\delta_{i}\bigg).
$$
\begin{prop}
From the above inequality, we obtain the following integral inequality:
\begin{align}
&\bigg\|\rho[h^{n+1}](t,x)-\int^{t}_{0}\int_{\mathbb{R}^{3}}e^{-(t-s)\nu}E[h^{n+1}](s,x+(t-s)v)\cdot\nabla_{v}f^{0}dvd\tau\bigg\|
_{\mathcal{F}^{\lambda'_{n}t+\mu'_{n}}}\notag\\
&\leq\frac{C \delta^{2}_{n}}{(\lambda_{n}-\lambda'_{n}-\frac{\nu}{2\pi})^{2}}+\int^{t}_{0}e^{-(t-s)\nu}\bigg(K_{1}^{'n}(t,s)+K_{0}^{'n}(t,s)
+\frac{c^{n}_{0}}{(1+s)^{2}}+C_{0}\nu\bigg)
\|\rho[h^{n+1}](s,\cdot)\|
_{\mathcal{F}^{\lambda'_{n}s+\mu'_{n}}}ds,
\end{align}
where $K_{1}^{'n}(t,s)=|K^{n}_{1}(t,s)|(1+s)\sum^{n}_{i=1}\delta_{i},\quad K_{0}^{'n}(t,s)=|K^{n}_{0}(t,s)|\sum^{n}_{i=1}\delta_{i},$
$c^{n}_{0}=\tilde{K}_{0}^{n+1}+\tilde{K}_{1}^{n+1}.$
\end{prop}
%$  \mathbf{Step}$ $\mathbf{6}.$ $(\mathbf{Iteration})$
%$  \mathbf{Step}$ $\mathbf{7}.$ $(\mathbf{Growth}$ $\mathbf{control})$

\begin{center}
\item\subsection{ Exponential moments of the kernel}
%{\bf\large 1. \quad Introduction }
\end{center}

First we analyzed  the influence of plasma echoes in section 5.1 when studying the stability of plasma particles as $t\rightarrow\infty.$ This is similar to the Diophantus condition of KAM theory. However, in this paper we also consider the collision among the particles. Because the collision is often considered as a energy dispersive process, then if the effect of  collision is very strong, maybe plasma echoes will not appear. For this,  we have to compare the influence on between  the resonances with the  collisions, although the collision is weak. We will give two important theorems, and the corresponding  proofs  can be also found in  [26,28].
\begin{prop}(Exponential moments of the kernel) Let $\gamma\in(1,\infty)$ be given. For any $\alpha\in(0,1),$ let $K^{(\alpha),\gamma}$ be defined
$$K^{(\alpha),\gamma}(t,s)=(1+s)\sup_{k,l\in\mathbb{Z}_{\ast}}\frac{e^{-\alpha|l|}e^{-\alpha(t-s)\frac{|k-l|}{t}}e^{-\alpha|k(t-s)+ls|}}{1+|k-l|^{\gamma}}.$$
Then for any $\gamma<\infty,$ there is $\bar{\alpha}=\bar{\alpha}(\gamma)>0$ such that if $\alpha\leq\bar{\alpha}$ and $\nu\in(0,\nu_{0}),$ then for any
$t>0,$
$$e^{-\nu t}\int^{t}_{0}K^{(\alpha),\gamma}(t,s)e^{\nu s}ds
\leq C\bigg(\frac{1}{\alpha\nu^{\gamma}t^{\gamma-1}}+\frac{1}{\alpha\nu^{\gamma}t^{\gamma}}\log\frac{1}{\alpha}
+\frac{1}{\alpha^{2}\nu^{1+\gamma}t^{1+\gamma}}+\bigg(\frac{1}{\alpha^{3}}+\frac{1}{\alpha^{2}\nu}\log\frac{1}{\alpha}\bigg)
e^{-\frac{\nu t}{4}}+\frac{e^{-\frac{\alpha t}{2}}}{\alpha^{3}}\bigg),$$
where $C=C(\gamma).$

In particular, if $\nu\leq\alpha,$ then
$e^{-\nu t}\int^{t}_{0}K^{(\alpha),\gamma}(t,s)e^{\nu s}ds\leq\frac{C(\gamma)}{\alpha^{3}\nu^{1+\gamma}t^{\gamma-1}}.$
\end{prop}

$Proof.$ Without loss of generality, we shall set $d=1$ and first consider $s\leq\frac{1}{2}t.$ We can write
$$K^{(\alpha)}(t,s)\leq(1+s)\sup_{l\in\mathbb{Z},k\in\mathbb{Z}}e^{-\alpha|l|}e^{-\alpha|k-l|/2}e^{-\alpha|k(t-s)+ls|}.$$
By symmetry we may also assume that $k>0.$

Explicit computations yield
\begin{align}
\int^{\frac{t}{2}}_{0}e^{-\alpha|k(t-s)+ls|}(1+s)ds\leq
\left\{\begin {array}{l}
\frac{1}{\alpha(l-k)}+\frac{1}{\alpha^{2}(l-k)^{2}},\quad \textmd{if}\quad l>k,\\
e^{-\alpha kt}(\frac{t}{2}+\frac{t^{2}}{8}), \quad \textmd{if}\quad l=k,\\
\frac{e^{-\alpha(k+l)t/2}}{\alpha|k-l|}(1+\frac{t}{2}), \quad \textmd{if}\quad -k\leq l<k,\\
\frac{2}{\alpha|k-l|}+\frac{2kt}{\alpha|k-l|^{2}}+\frac{1}{\alpha^{2}|k-l|^{2}}, \quad \textmd{if}\quad l<-k.\\
\end{array}\right.
\end{align}

So from (7.4), we have
$$e^{-\nu t}\int^{\frac{t}{2}}_{0}e^{-\alpha|k(t-s)+ls|}(1+s)e^{\nu s}ds$$
$$\leq e^{-\frac{\nu t}{4}}\bigg(\frac{3}{\alpha|k-l|}+\frac{1}{\alpha^{2}|k-l|^{2}}+\frac{8z}{\alpha\nu|k-l|}\bigg)
1_{k\neq l}+e^{-\frac{t\alpha}{2}}\bigg(\frac{z}{\alpha}+\frac{8z^{2}}{\alpha^{2}}\bigg)1_{l=k},$$
where $z=\sup_{x}xe^{-x}=e^{-1}.$

Using the bounds (for $\alpha\sim 0^{+}$)
$$\sum_{l\in\mathbb{Z}}e^{-\alpha l}=O\bigg(\frac{1}{\alpha}\bigg),\quad \sum_{l\in\mathbb{Z}}\frac{e^{-\alpha l}}{l}=O\bigg(\log\frac{1}{\alpha}\bigg),\quad
\sum_{l\in\mathbb{Z}}\frac{e^{-\alpha l}}{l^{2}}=O(1),$$
we end up, for $\alpha\leq\frac{1}{4},$ with a bound like
$$e^{-\nu t}\int^{\frac{t}{2}}_{0}K^{(\alpha)}(t,s)e^{\nu s}ds
\leq C\bigg[e^{-\frac{\nu t}{4}}\bigg(\frac{1}{\alpha^{3}}+\frac{1}{\alpha^{2}\nu}\bigg)
+\frac{e^{-\alpha t/2}}{\alpha^{3}}\bigg].$$

Next we turn to the more delicate contribution of $s\geq\frac{1}{2}t.$ We write
\begin{align}
K^{(\alpha)}(t,s)\leq(1+s)\sup_{l\in\mathbb{Z}_{\ast}}e^{-\alpha|l|}\sup_{k\in\mathbb{Z}}\frac{e^{-\alpha|k(t-s)+ls|}}{1+|k-l|^{\gamma}}.
\end{align}

Without loss of generality, we restrict the supremum $l>0.$ The function $x\rightarrow(1+|x-l|^{\gamma})^{-1}e^{-\alpha|x(t-s)+ls|}$
is decreasing for $x\geq l,$ increasing for $x\leq-ls/(t-s);$ and on the interval $[-ls/(t-s),l],$ its logarithmic derivative goes
from $\bigg(-\alpha+\frac{\gamma/lt}{1+((t-s)/lt)^{\gamma}}\bigg)(t-s)$ to $-\alpha(t-s).$ It is easy to check that a given integer $k$
occurs in the supremum only for some times $s$ satisfying $k-1<-ls/(t-s)<k+1.$ We can assume $k\geq0,$ then $k-1<\frac{ls}{t-s}<k+1$ holds, and it
is equivalent to $\frac{k-1}{k+l-1}t<s<\frac{k+1}{k+l+1}t.$ More importantly, $\tau>\frac{1}{2}t$ implies that $k\geq l.$ Thus, for $t\geq\frac{\gamma}{\alpha},$
we have \begin{align}
e^{-\nu t}\int^{t}_{\frac{t}{2}}K^{(\alpha)}(t,s)e^{\nu s}ds
\leq e^{-\nu t}\sum^{\infty}_{l=1}e^{-\alpha l}\sum^{\infty}_{k=l}\int^{\frac{(k+1)t}{k+l+1}}_{\frac{(k-1)t}{k+l-1}}(1+s)
\frac{e^{\alpha|k(t-s)-ls|}e^{\nu s}}{1+|k+l|^{\gamma}}ds.
\end{align}
 For $t\leq\frac{\gamma}{\alpha},$ it is easy to check that  $e^{-\nu t}\int^{t}_{\frac{t}{2}}K^{(\alpha)}(t,s)
e^{\nu s}ds\leq\frac{\gamma}{2\alpha}$ holds. Next we shall focus on $(7.6).$
According to $s$ smaller or larger than $kt/(k+l),$ we separate the integral in the right-hand side of (7.6) into two parts, and by simple computation, we get
the explicit bounds
$$ e^{-\nu t}\int^{kt/(k+l)}_{(k-l)t/(k+l-1)}(1+s)e^{-\alpha|k(t-s)-ls|}e^{\nu s}ds\leq e^{-\frac{\nu lt}{k+l}}
\bigg(\frac{1}{\alpha(k+l)}+\frac{kt}{\alpha(k+l)^{2}}\bigg),$$
$$e^{-\nu t}\int^{\frac{(k+1)t}{k+l+1}}_{\frac{kt}{k+l}}(1+s)e^{-\alpha|k(t-s)-ls|}e^{\nu s}ds
\leq e^{-\frac{\nu lt}{k+l+1}}
\bigg(\frac{1}{\alpha(k+l)}+\frac{kt}{\alpha(k+l)^{2}}+\frac{1}{\alpha^{2}(k+l)^{2}}\bigg).$$
Hence, (7.6) is bounded above by
\begin{align}
C\sum^{\infty}_{l=1}e^{-\alpha l}\sum^{\infty}_{k=l}\bigg(\frac{1}{\alpha^{2}(k+l)^{2+\gamma}}+\frac{1}{\alpha(k+l)^{1+\gamma}}
+\frac{kt}{\alpha(k+l)^{2+\gamma}}\bigg)e^{-\frac{\nu lt}{k+l}}.
\end{align}

We consider the first term $I(t)$ of (7.7) and change variables $(x,y)\mapsto (x,u),$ where $u(x,y)=\frac{\nu xt}{x+y},$ then we can find that
$$I(t)=\frac{1}{\alpha^{2}\nu^{1+\gamma}t^{1+\gamma}}\int^{\infty}_{1}\frac{e^{-\alpha x}}{x^{1+\gamma}}dx\int^{\nu t/2}_{0}e^{-u}u^{\gamma}du
=O\bigg(\frac{1}{\alpha^{2}\nu^{1+\gamma}t^{1+\gamma}}\bigg).$$

The same computation for the second integral of (7.7) yields
$$\frac{1}{\alpha\nu^{\gamma}t^{\gamma}}\int^{\infty}_{1}\frac{e^{-\alpha x}}{x^{\gamma}}dx\int^{\nu t/2}_{0}e^{-u}u^{\gamma-1}du
=O\bigg(\frac{1}{\alpha\nu^{\gamma}t^{\gamma}}\bigg).$$

Finally, we estimate the last term of (7.7) that is the worst. It yields a contribution $\frac{t}{\alpha}\sum^{\infty}_{l=1}e^{-\alpha l}\sum^{\infty}_{k=l}
\frac{e^{-\nu ltk}/(k+l)}{(k+l)^{2+\gamma}}.$ We compare this with the integral $\frac{t}{\alpha}\int^{\infty}_{1}e^{-\alpha x}\int^{\infty}_{x}
\frac{e^{-\nu ltx/(x+y)}}{(x+y)^{2+\gamma}}dydx,$
and the same change of variables as before equates this with
$$\frac{1}{\alpha\nu^{\gamma}t^{\gamma-1}}\int^{\infty}_{1}\frac{e^{-\alpha x}}{x^{\gamma}}dx\int^{\frac{\nu t}{2}}_{0}e^{-u}u^{\gamma-1}du
-\frac{1}{\alpha\nu^{1+\gamma}t^{\gamma}}\int^{\infty}_{1}\frac{e^{-\alpha x}}{x^{\gamma}}dx\int^{\frac{\nu t}{2}}_{0}e^{-u}u^{\gamma}du
=O\bigg(\frac{1}{\alpha\nu^{\gamma}t^{\gamma-1}}\bigg).$$
The proof of Proposition 7.2 follows by collecting all these bounds and keeping only the worst one.
To finish the growth control,  we have to prove the following result.
\begin{prop} With the same notation as in Proposition 7.2, for any $\gamma>1,$ we have
\begin{align}
\sup_{s\geq0}e^{\nu s}\int^{\infty}_{s}e^{-\nu t}K^{(\alpha),\gamma}(t,s)dt\leq C(\gamma)\bigg(\frac{1}{\alpha^{2}\nu}+
\frac{1}{\alpha\nu^{\gamma}}\bigg)
\end{align}
\end{prop}

$Proof.$ We first still reduce to $d=1,$ and split the integral as
$$e^{\nu s}\int^{\infty}_{s}e^{-\nu t}K^{(\alpha),\gamma}(t,s)dt
=e^{\nu s}\int^{\infty}_{2s}e^{-\nu t}K^{(\alpha),\gamma}(t,s)dt
+e^{\nu s}\int^{2s}_{s}e^{-\nu t}K^{(\alpha),\gamma}(t,s)dt=I_{1}+I_{2}.$$
For the first term $I_{1},$ we have $K^{(\alpha),\gamma}(t,s)\leq(1+s)\sum^{\infty}_{k=2}\sum_{l\in\mathbb{Z}_{\ast}}
e^{-\alpha|l|-\alpha|k-l|/2}\leq\frac{C(1+s)}{\alpha^{2}},$ and thus
$e^{\nu s}$
$\int^{\infty}_{s}e^{-\nu t}$
$K^{(\alpha),\gamma}(t,s)dt\leq\frac{C}{\nu\alpha^{2}}.$

We treat the second term $I_{2}$ as in the proof of Proposition 7.2:
$$e^{\nu s}\int^{\infty}_{s}e^{-\nu t}K^{(\alpha),\gamma}(t,s)dt\leq e^{\nu s}(1+s)\sum^{\infty}_{l=1}e^{-\alpha l}
\sum^{\infty}_{k=l}\int^{\frac{(k+l-1)s}{k-1}}_{\frac{(k+l+1)s}{k+1}}\frac{e^{-\alpha|k(t-s)-ls|}}{1+(k+l)^{\gamma}}e^{-\nu t}dt
\leq\frac{C}{\alpha\nu^{\gamma}},$$
where the last inequality is obtained by the same method in Proposition 7.2 with the change of variable $u=\frac{\nu x s}{y}.$
\begin{center}
\item\subsection{ Growth control}
%{\bf\large 1. \quad Introduction }
\end{center}
From now on, we will state the main result of this section that is the same with section 7.4 in [26,28], the detailed proof can be found in appendix ( also see [26,28] ).
We define $\|\Phi(t)\|_{\lambda}=\sum_{k\in\mathbb{Z}_{\ast}^{3}}|\Phi(k,t)|e^{2\pi\lambda|k|}.$

\begin{thm} Assume that $f^{0}(v),W=W(x)$ satisfy the conditions of Theorem 0.1, and the Stability condition holds.
Let $A\geq0,\mu\geq0$ and $\lambda\in(0,\lambda^{\ast}]$
with $0<\lambda^{\ast}<\lambda_{0}.$ Let $(\Phi(k,t))_{k\in\mathbb{Z}_{\ast}^{3},t\geq0}$ be a continuous functions of $t\geq0,$ valued in
$\mathbb{C}^{\mathbb{Z}^{3}_{\ast}},$ such that for all $t\geq0,$
\begin{align}\|\Phi(t)-\int^{t}_{0}e^{-(t-s)\nu}K^{0}(t-s)\Phi(s)ds\|_{\lambda t+\mu}\leq A+\int^{t}_{0}e^{-(t-s)\nu}(K_{0}(t,s)+K_{1}(t,s)+\frac{c_{0}}{(1+s)^{m}})
\|\Phi(s)\|_{\lambda s+\mu} ds,
\end{align}
where $c_{0}\geq0,m>1,$ and $K_{0}(t,s),K_{1}(t,s)$ are non-negative kernels. Let $\varphi(t)=\|\Phi(t)\|_{\lambda t+\mu}.$
Then we have the following:

(i) Assume that $\gamma>1$ and $K_{1}=cK^{(\alpha),\gamma}$ for some $c>0,\alpha\in(0,\bar{\alpha}(\gamma)),$ where $K^{(\alpha),\gamma},\bar{\alpha}(\gamma)$
are the same with that
 defined by Proposition 7.2. Then there are positive constants $C,\chi,$ depending only on $\gamma,\lambda^{\ast},\lambda_{0},\kappa,c_{0}, C_{W}$ and $ m,$ uniform
 as $\gamma\rightarrow1,$ such that if $\sup_{t\geq0}\int^{t}_{0}K_{0}(t,s)ds\leq\chi$ and $\sup_{t\geq0}$
 $\bigg(\int^{t}_{0}K_{0}(t,s)^{2}ds\bigg)^{\frac{1}{2}}+\sup_{t\geq0}\int^{\infty}_{t}K_{0}(t,s)dt\leq1,$ then for any $\nu\in(0,\alpha),$ for all $t\geq0,$
 \begin{align}\varphi(t)\leq CA\frac{1+c^{2}_{0}}{\sqrt{\nu}}e^{Cc_{0}}\bigg(1+\frac{c}{\alpha\nu}\bigg)e^{CT}e^{Cc(1+T^{2})}e^{\nu t}
 \end{align}
 where $T_{\nu}=C\max\bigg\{\bigg(\frac{c^{2}}{\alpha^{5}}\nu^{2+\gamma}\bigg)^{\frac{1}{\gamma-1}},
\bigg(\frac{c}{\alpha^{2}}\nu^{\frac{1}{2}+\gamma}\bigg)^{\frac{1}{\gamma-1}},\bigg(\frac{c^{2}_{0}}{\nu}\bigg)^{\frac{1}{2m-1}}\bigg\}.$

(ii) Assume that $K_{1}=\sum^{N}_{j=1}c_{j}K^{(\alpha_{j},1)}$ for some $\alpha_{j}\in(0,\bar{\alpha}(\gamma)),$ where $\bar{\alpha}(\gamma)$ also appears in proposition 7.2;
then there is a numeric constant $\Gamma>0$ such that whenever $1\geq\nu\geq\Gamma\sum^{N}_{j=1}\frac{c_{j}}{\alpha^{3}_{j}},$ with the same notation
 as in (I), for all $t\geq0,$  one has,

  \begin{align}\varphi(t)\leq CA\frac{1+c^{2}_{0}}{\sqrt{\nu}}e^{Cc_{0}}\bigg(1+\frac{c}{\alpha\nu}\bigg)e^{CT}e^{Cc(1+T^{2})}e^{\nu t}
 \end{align}
 where $c=\sum^{N}_{j=1}c_{j}$ and $T=\max\bigg\{\frac{1}{\nu^{2}}\sum^{N}_{j=1}\frac{c_{j}}{\alpha^{3}_{j}},
 \bigg(\frac{c^{2}_{0}}{\nu}\bigg)^{\frac{1}{2m-1}}\bigg\}.$

\end{thm}

\begin{cor} Assume that $f^{0}=f^{0}(v),$ under the assumptions of Theorem 0.1,
 we pick up $\nu_{n}$ such that $\lim_{n\rightarrow\infty}\nu_{n}=\nu;$
recalling that $\hat{\rho}(t,0)=0,$ and that our conditions imply an upper bound on $c_{n}$ and $c^{n}_{0},$ we have the uniform control,
$$\|\rho[h^{n+1}](t,\cdot)\|
_{\mathcal{F}^{\lambda'_{n}t+\mu'_{n}}}\leq\frac{C\delta^{2}_{n}(1+c^{n}_{0})^{2}}{\sqrt{\nu}(\lambda_{n}-\lambda'_{n})^{2}}
\bigg(1+\frac{1}{\alpha_{n}\nu_{n}^{\frac{3}{2}}}\bigg)e^{CT^{2}_{n}},$$
where $T_{n}=C\bigg(\frac{1}{\alpha^{5}_{n}\nu_{n}}\bigg)^{\frac{1}{\gamma-1}}.$
\end{cor}

$Proof.$ From Propositions 7.1-7.3, we know that
$$\int^{t}_{0}K^{n}_{0}(t,s)ds\leq C_{W}\sum^{n}_{i=1}\frac{\delta_{i}}{\pi(\lambda_{i}-\lambda'_{n})},\quad
\int^{\infty}_{s}K^{n}_{0}(t,s)ds\leq C_{W}\sum^{n}_{i=1}\frac{\delta_{i}}{\pi(\lambda_{i}-\lambda'_{n})},$$
$$\bigg(\int^{t}_{0}K^{n}_{0}(t,s)^{2}ds\bigg)^{\frac{1}{2}}\leq C_{W}\sum^{n}_{i=1}\frac{\delta_{i}}{\sqrt{2\pi(\lambda_{i}-\lambda'_{n})}}.$$
Here $\alpha_{n}=\pi\min\{(\mu_{n}-\mu'_{n}),(\lambda_{n}-\lambda'_{n})\},$ and assume $\alpha_{n}$ is smaller than $\bar{\alpha}(\gamma)$ in Theorem 7.4, and
that
$$\bigg(C^{4}_{\omega}\bigg(C'_{0}+\sum^{n}_{i=1}\delta_{i}+1\bigg)\bigg(\sum^{n}_{j=1}\frac{\delta_{j}}
{2\pi(\lambda_{j}-\lambda'_{j})^{6}}\bigg)\leq\frac{1}{8}\quad\quad (\mathbf{VII})$$
$$C_{W}\sum^{n}_{i=1}\frac{\delta_{i}}{\sqrt{2\pi(\lambda_{i}-\lambda'_{n})}}\leq\frac{1}{4},\quad
\sum^{n}_{i=1}\frac{\delta_{i}}{\pi(\lambda_{i}-\lambda'_{n})}\leq\max\{\chi,\frac{1}{8}\}.\quad (\mathbf{VIII})$$
Applying Theorem 7.4, we can deduce that for any $\nu_{n}\in(0,\alpha_{n})$ and $t\geq0,$
$$\|\rho[h^{n+1}](t,\cdot)\|
_{\mathcal{F}^{\lambda'_{n}t+\mu'_{n}}}$$
$$\leq\frac{C\delta^{2}_{n}(1+c^{n}_{0})^{2}}{\sqrt{\nu}(\lambda_{n}-\lambda'_{n})^{2}}
\bigg(1+\frac{1}{\alpha_{n}\nu_{n}^{\frac{3}{2}}}\bigg)e^{CT^{2}_{n}},$$
where $T_{n}=C\bigg(\frac{1}{\alpha^{5}_{n}\nu_{n}}\bigg)^{\frac{1}{\gamma-1}}.$
%\begin{align}
%&\|\rho[h^{n+1}](t,x)\|
%_{\mathcal{F}^{(\lambda'_{n}-B_{0})t+\mu'_{n}}}\leq \frac{C\delta^{2}_{n}(1+c^{n}_{0})^{2}}{\sqrt{\varepsilon}(\lambda_{n}-\lambda'_{n})^{2}}
%e^{Cc^{n}_{0}}\bigg(1+\frac{c_{n}}{\alpha_{n}\varepsilon}\bigg)e^{CT_{\varepsilon,n}}e^{Cc_{n}(1+T^{2}_{\varepsilon,n})}e^{\varepsilon t},\notag\\
%\end{align}
%where $c_{n}=2C\sum^{n}_{i=1}\delta_{i},$
%and $T_{\varepsilon,n}=C_{\gamma}\max\bigg\{\bigg(\frac{c^{2}_{n}}{\alpha^{5}_{n}}\varepsilon^{2+\gamma}\bigg)^{\frac{1}{\gamma-1}},
%\bigg(\frac{c_{n}}{\alpha^{2}_{n}}\varepsilon^{\frac{1}{2}+\gamma}\bigg)^{\frac{1}{\gamma-1}},\frac{(c^{n}_{0})^{\frac{2}{3}}}{\varepsilon^{\frac{1}{3}}}\bigg\}.$

\begin{center}
\item\subsection{ Estimates related to $h^{n+1}(t,X^{n}_{s,t}(x,v),V^{n}_{s,t}(x,v)))$}
%{\bf\large 1. \quad Introduction }
\end{center}

To finish Proposition 3.11 and Proposition 3.12,
 we shall again use the Vlasov equation. We rewrite it as
$$h^{n+1}(t,X^{n}_{s,t}(x,v),V^{n}_{s,t}(x,v))=\int^{t}_{0}e^{-(t-s)\nu}\Sigma^{n+1}(s,X^{n}_{s,\tau}(x,v),V^{n}_{s,\tau}(x,v))d\tau.$$
Then we get
$$\|h^{n+1}(t,X^{n}_{s,t}(x,v),V^{n}_{s,t}(x,v))\|_{\mathcal{Z}^{\lambda'_{n}(1+b),\mu'_{n};1}_{t-\frac{bt}{1+b}}}$$
$$\leq\int^{t}_{0}e^{-(t-s)\nu}\|\Sigma^{n+1}(\tau,X^{n}_{s,\tau}(x,v),V^{n}_{s,\tau}(x,v))\|_{\mathcal{Z}^{\lambda'_{n}(1+b),\mu'_{n};1}
_{t-\frac{bt}{1+b}}}d\tau$$
$$=\int^{t}_{0}e^{-(t-s)\nu}\|\Sigma^{n+1}(\tau,\Omega^{n}_{s,\tau}(x,v))\|_{\mathcal{Z}^{\lambda'_{n}(1+b),\mu'_{n};1}_{\tau-\frac{bt}{1+b}}}d\tau$$
$$\leq\int^{t}_{0}e^{-(t-s)\nu}\|\mathcal{E}^{n+1}_{s,t}\|_{\mathcal{Z}^{\lambda'_{n}(1+b),\mu'_{n};1}_{\tau-\frac{bt}{1+b}}}
\|G^{n}_{s,t}\|_{\mathcal{Z}^{\lambda'_{n}(1+b),\mu'_{n};1}_{\tau-\frac{bt}{1+b}}}ds+\int^{t}_{0}e^{-(t-s)\nu}
\bigg[\|H^{n}_{s,t}\|_{\mathcal{Z}^{\lambda'_{n}(1+b),\mu'_{n};1}_{\tau-\frac{bt}{1+b}}}
\|\mathcal{E}^{n}_{s,t}\|_{\mathcal{Z}^{\lambda'_{n}(1+b),\mu'_{n};1}_{\tau-\frac{bt}{1+b}}}
\bigg]ds,$$

therefore, from the induction assumptions, we obtain
 \begin{align}
\sup_{0\leq s\leq t}\|
h^{n+1}\circ\Omega^{n}_{t,\tau}\|_{\mathcal{Z}^{\lambda'_{n}(1+b),
\mu'_{n};1}_{s-\frac{bt}{1+b}}}\leq\frac{\delta^{2}_{n}}{(\lambda_{n}-\lambda'_{n}-\frac{\nu}{2\pi})^{2}}+\bigg(\sum^{n}_{i=1}\delta_{i}\bigg)
\sup_{0\leq s\leq t}\|\rho^{n+1}\|_{\mathcal{F}^{\lambda'_{n}s+\mu'_{n}}},
\end{align}
this is the conclusion of Proposition 3.11.
Next we show the control on $h^{i}.$
\begin{lem} For  any $ n \geq i\geq1,$
$$\| \nabla_{v}(h^{i}_{s}\circ\Omega^{i-1}_{t,s})-\langle\nabla_{v}(h^{i}_{s}\circ\Omega^{i-1}_{t,s})
\rangle)\|_{\mathcal{Z}^{\lambda_{i}(1+b),
\mu'_{i};1}_{s-\frac{bt}{1+b}}}\leq (1+s)\delta_{i}.$$
 \end{lem}
$Proof.$  First, we consider $i=1.$

 In fact,
$$\|\nabla_{v}h^{1}_{s}(t,x-v(t-s),v)-\langle\nabla_{v}h^{1}_{s}v
\rangle\|_{\mathcal{Z}^{\lambda'_{1}(1+b),
\mu'_{1};1}_{s-\frac{bt}{1+b}}}$$
$$\leq\|\nabla_{v}h^{1}
_{s}\|_{\mathcal{Z}^{\lambda'_{1}(1+b),
\mu'_{1};1}_{s-\frac{bt}{1+b}}}+\bigg\|\int_{\mathbb{T}^{3}}\nabla_{v}h^{1}
_{s}dx\bigg\|_{\mathcal{C}^{\lambda'_{1}(1+b);1}}$$
$$\leq\|\nabla_{v}h^{1}
_{s}\|_{\mathcal{Z}^{\lambda'_{1}(1+b),
\mu'_{1};1}_{s-\frac{bt}{1+b}}}\leq\|h^{1}
_{s}\|_{\mathcal{Z}^{\lambda_{1}(1+b),
\mu_{1};1}_{s-\frac{bt}{1+b}}}\leq (1+s)\delta_{1},$$
 where we use the property $(v)$ of Proposition 1.5.

$$\bigg\|\int_{\mathbb{T}^{3}}\nabla_{v}h^{1}
_{s}dx\bigg\|_{\mathcal{C}^{\lambda'_{1}(1+b);1}}=\|\langle\nabla_{v}h^{1}
_{s}\rangle\|_{\mathcal{C}^{\lambda'_{1}(1+b);1}}$$
$$=\|\langle(\nabla_{v}+s\nabla_{x})h^{1}
_{s}\rangle\|_{\mathcal{C}^{\lambda'_{1}(1+b);1}}
\leq\|(\nabla_{v}+s\nabla_{x})h^{1}
_{s}\|_{\mathcal{Z}^{\lambda'_{1}(1+b),\mu'_{1};1}_{s-\frac{bt}{1+b}}}\leq  \delta_{1},$$
where we use $(vi)$ of Proposition 1.5.

Suppose that $i=k,$ the conclusion holds, that is,
$$\bigg\|\nabla_{v}(h^{k}_{s}\circ\Omega^{k-1}_{t,s})-\langle\nabla_{v}(h^{k}_{s}\circ\Omega^{k-1}_{t,s})
\rangle\bigg\|_{\mathcal{Z}^{\lambda'_{k}(1+b),
\mu'_{k};1}_{s-\frac{bt}{1+b}}}$$
$$\leq\bigg\|h^{k}_{s}\circ\Omega^{k-1}_{t,s}-\langle h^{k}_{s}\circ\Omega^{k-1}_{t,s}
\rangle\bigg\|_{\mathcal{Z}^{\lambda_{k}(1+b),
\mu_{k};1}_{s-\frac{bt}{1+b}}}\leq (1+s)\delta_{k}.$$

$\mathbf{We}$ $\mathbf{need}$ $\mathbf{to}$ $\mathbf{show}$ the conclusion still holds for $i=k+1.$ We can get the estimate for
$h^{k+1}(t,X^{k}_{t,s}(x,v),V^{k}_{t,s}(x,v))$ from (3.11).
%For this, we have to reconsider $$h^{k+1}(t,X^{k}_{0,t}(x,v),V^{k}_{0,t}(x,v))=\int^{t}_{0}\Sigma^{k+1}(s,X^{k}_{0,s}(x,v),V^{k}_{0,s}(x,v))ds,$$
%then for $\tau\geq t,$ composing with $(X^{k}_{\tau,0}(x,v),V^{k}_{\tau,0}(x,v)),$ this gives
%$$h^{k+1}(t,X^{k}_{\tau,t}(x,v),V^{k}_{\tau,t}(x,v))=\int^{t}_{0}\Sigma^{k+1}(s,X^{k}_{\tau,s}(x,v),V^{k}_{\tau,s}(x,v))ds.$$

           Note that \begin{align}
           \left\{
           \begin{array}{l}
           (\nabla h)\circ\Omega=(\nabla\Omega)^{-1}\nabla(h\circ\Omega)\\
           (\nabla^{2} h)\circ\Omega=(\nabla\Omega)^{-2}\nabla^{2}(h\circ\Omega)-(\nabla\Omega)^{-1}\nabla^{2}\Omega(\nabla\Omega)^{-1}(\nabla h\circ\Omega).
           \end{array}
           \right.
           \end{align}
         Therefore, from (7.14),   we get
           \begin{align}
          &\|(\nabla h_{s}^{n+1})\circ\Omega_{t,s}^{n}\|_{\mathcal{Z}_{s-\frac{bt}{1+b}}
           ^{\lambda'^{\dagger}_{n}(1+b),\mu'^{\dagger}_{n};1}}
           \leq C\|\nabla (h_{s}^{n+1}\circ\Omega_{t,s}^{n})\|_{\mathcal{Z}_{s-\frac{bt}{1+b}}
           ^{\lambda'^{\dagger}_{n}(1+b),\mu'^{\dagger}_{n};1}}\notag\\
           &\leq\frac{C(1+s)}{\min\{\lambda'_{n}-\lambda'^{\dagger}_{n},\mu'_{n}-\mu'^{\dagger}_{n}\}}\|h_{s}^{n+1}\circ\Omega_{t,s}^{n}\|
           _{\mathcal{Z}_{s-\frac{bt}{1+b}}
           ^{\lambda'_{n}(1+b),\mu'_{n};1}},
           \end{align}
           and

           \begin{align}
          &\|(\nabla^{2} h_{s}^{n+1})\circ\Omega_{t,s}^{n}\|_{\mathcal{Z}_{s-\frac{bt}{1+b}}
           ^{\lambda'^{\dagger}_{n}(1+b),\mu'^{\dagger}_{n};1}}\notag\\
           &\leq C\bigg[\|\nabla^{2} (h_{s}^{n+1}\circ\Omega_{t,s}^{n})\|_{\mathcal{Z}_{s-\frac{bt}{1+b}}
           ^{\lambda'^{\dagger}_{n}(1+b),\mu'^{\dagger}_{n};1}}\notag\\
           & +\|\nabla^{2}\Omega^{n}_{t,s}\|_{\mathcal{Z}_{s-\frac{bt}{1+b}}
           ^{\lambda'^{\dagger}_{n}(1+b),\mu'^{\dagger}_{n}}}\|(\nabla h_{s}^{n+1})\circ\Omega_{t,s}^{n}\|_{\mathcal{Z}_{s-\frac{bt}{1+b}}
           ^{\lambda'^{\dagger}_{n}(1+b),\mu'^{\dagger}_{n};1}}\bigg]\notag\\
           &\leq\frac{C(1+s)^{2}}{\min\{\lambda'_{n}-\lambda'^{\dagger}_{n},\mu'_{n}-\mu'^{\dagger}_{n}\}}\|h_{s}^{n+1}\circ\Omega_{t,s}^{n}\|_{\mathcal{Z}_{s-\frac{bt}{1+b}}
           ^{\lambda'_{n}(1+b),\mu'_{n};1}}\notag\\
            &\leq\frac{C(1+s)^{2}}{\min\{\lambda'_{n}-\lambda'^{\dagger}_{n},\mu'_{n}-\mu'^{\dagger}_{n}\}}\|h_{s}^{n+1}\circ\Omega_{t,s}^{n}\|
            _{\mathcal{Z}_{s+\frac{bt}{1-b}}
           ^{\lambda'_{n}(1-b),\mu'_{n};1}},
           \end{align}

           We first write $$\nabla(h_{s}^{n+1}\circ\Omega_{t,s}^{n})-(\nabla h_{s}^{n+1})\circ\Omega_{t,s}^{n}=\nabla(\Omega_{t,s}^{n}-Id)
           \cdot[(\nabla h_{s}^{n+1})\circ\Omega_{t,s}^{n}],$$
           and we get $$\|\nabla(h_{s}^{n+1}\circ\Omega_{t,s}^{n})-(\nabla h_{s}^{n+1})\circ\Omega_{t,s}^{n}\|_{\mathcal{Z}^{\lambda_{n}'^{\dag}(1+b)
           ,\mu_{n}'^{\dag};1}_{s-\frac{bt}{1+b}}}$$
           $$\leq\|\nabla(\Omega_{t,s}^{n}-Id)\|_{\mathcal{Z}^{\lambda_{n}'^{\dag}(1+b)
           ,\mu_{n}'^{\dag}}_{s-\frac{bt}{1+b}}}\|(\nabla h_{s}^{n+1})\circ\Omega_{t,s}^{n}\|_{\mathcal{Z}^{\lambda_{n}'^{\dag}(1+b)
           ,\mu_{n}'^{\dag};1}_{s-\frac{bt}{1+b}}}$$
           $$\leq C\bigg(\frac{1+s}{\min\{\lambda_{n}'-\lambda_{n}'^{\dag},\mu_{n}'-\mu_{n}'^{\dag}\}}\bigg)^{2}
           \|\Omega_{t,s}^{n}-Id\|_{\mathcal{Z}^{\lambda_{n}'(1+b)
           ,\mu_{n}'}_{s-\frac{bt}{1+b}}}\| h_{s}^{n+1}\circ\Omega_{t,s}^{n}\|_{\mathcal{Z}^{\lambda_{n}'(1+b)
           ,\mu_{n}';1}_{s-\frac{bt}{1+b}}}$$
           $$\leq\frac{CC^{4}_{\omega}}{\min\{\lambda_{n}'-\lambda_{n}'^{\dag},\mu_{n}'-\mu_{n}'^{\dag}\}^{2}}\bigg(
           \sum^{n}_{k=1}\frac{\delta_{k}}{(2\pi(\lambda_{k}-\lambda'_{k}))^{6}}\bigg)(1+s)^{-2}\| h_{s}^{n+1}\circ\Omega_{t,s}^{n}\|_{\mathcal{Z}^{\lambda_{n}'(1+b)
           ,\mu_{n}';1}_{s-\frac{bt}{1+b}}},$$
           the above inequality implies $\nabla(h_{s}^{n+1}\circ\Omega_{t,s}^{n})\simeq(\nabla h_{s}^{n+1})\circ\Omega_{t,s}^{n}$ as $s\rightarrow\infty.$

           Since $$\|\nabla(h_{s}^{n+1}\circ\Omega_{t,s}^{n})\|_{\mathcal{Z}^{\lambda_{n}'^{\dag}(1+b)
           ,\mu_{n}'^{\dag};1}_{s-\frac{bt}{1+b}}}\leq C\bigg(\frac{1+s}{\min\{\lambda_{n}'-\lambda_{n}'^{\dag},\mu_{n}'-\mu_{n}'^{\dag}\}}\bigg)
           \| h_{s}^{n+1}\circ\Omega_{t,s}^{n}\|_{\mathcal{Z}^{\lambda_{n}'(1+b)
           ,\mu_{n}';1}_{s-\frac{bt}{1+b}}}$$
           and
 $$\|\nabla_{x}(h_{s}^{n+1}\circ\Omega_{t,s}^{n})\|_{\mathcal{Z}^{\lambda_{n}'^{\dag}(1+b)
           ,\mu_{n}'^{\dag};1}_{s-\frac{bt}{1+b}}}+\|(\nabla_{x}+s\nabla_{v})(h_{s}^{n+1}\circ\Omega_{t,s}^{n})\|_{\mathcal{Z}^{\lambda_{n}'^{\dag}(1+b)
           ,\mu_{n}'^{\dag};1}_{s-\frac{bt}{1+b}}}$$
          $$ \leq \frac{C}{\min\{\lambda_{n}'-\lambda_{n}'^{\dag},\mu_{n}'-\mu_{n}'^{\dag}\}}
           \| h_{s}^{n+1}\circ\Omega_{t,s}^{n}\|_{\mathcal{Z}^{\lambda_{n}'(1+b)
           ,\mu_{n}';1}_{s-\frac{bt}{1+b}}},$$
           we have $$\|(\nabla_{x}h_{s}^{n+1})\circ\Omega_{t,s}^{n}\|_{\mathcal{Z}^{\lambda_{n}'^{\dag}(1+b)
           ,\mu_{n}'^{\dag};1}_{s-\frac{bt}{1+b}}}+\|((\nabla_{x}+s\nabla_{v})(h_{s}^{n+1})\circ\Omega_{t,s}^{n}\|_{\mathcal{Z}^{\lambda_{n}'^{\dag}(1+b)
           ,\mu_{n}'^{\dag};1}_{s-\frac{bt}{1+b}}}$$
           $$\leq C\bigg(\frac{C^{4}_{\omega}}{\min\{\lambda_{n}'-\lambda_{n}'^{\dag},\mu_{n}'-\mu_{n}'^{\dag}\}^{2}}\bigg(
           \sum^{n}_{k=1}\frac{\delta_{k}}{(2\pi(\lambda_{k}-\lambda'_{k}))^{6}}\bigg)
           +\frac{1}{\min\{\lambda_{n}'-\lambda_{n}'^{\dag},\mu_{n}'-\mu_{n}'^{\dag}\}}\bigg) \| h_{s}^{n+1}\circ\Omega_{t,s}^{n}\|_{\mathcal{Z}^{\lambda_{n}'(1+b)
           ,\mu_{n}';1}_{s-\frac{bt}{1+b}}}.$$

\begin{center}
\item\subsection{Conclusion}
%{\bf\large 1. \quad Introduction }
\end{center}

If we define $$\lambda_{n+1}=\lambda'^{\dagger}_{n},\quad\quad \mu_{n+1}=\mu'^{\dagger}_{n},$$
then we see that the $n+1$th step of the inductive hypothesis have all been established with
\begin{align}
\delta_{n+1}=\frac{C_{F}(1+C_{F})(1+C^{4}_{\omega})e^{CT^{2}_{n}}}{\min\{\lambda'_{n}-\lambda_{n+1},\mu'_{n}-\mu_{n+1}\}^{9}}
\max\bigg\{1,\sum^{n}_{i=1}\delta_{k}\bigg\}\bigg(1+\sum^{n}_{i=1}\frac{\delta_{i}}{\nu^{6}}\bigg)\delta^{2}_{n}.
\end{align}

For any $n\geq1,$ we set
$\lambda_{n}-\lambda'_{n}=\lambda'_{n}-\lambda_{n+1}=\mu_{n}-\mu'_{n}=\mu'_{n}-\mu_{n+1}=\frac{\Lambda}{n^{2}}$
for some $\Lambda>0.$ By choosing $\Lambda$ small enough, we can make sure that the conditions
$2\pi(\lambda_{k}-\lambda'_{k})<1\quad \textmd{and}\quad 2\pi(\mu_{k}-\mu'_{k})<1$
are satisfied for all $k,$ as well as the other smallness assumptions made throughout this section.  We also have
$\lambda_{k}-\lambda'_{k}\geq\frac{\Lambda}{k^{2}}.$  $(\mathbf{I})-(\mathbf{VIII})$  will be satisfied if  we  choose constants $\Lambda,\omega>0$
such that $\sum^{\infty}_{i=1}i^{12}\delta_{i}\leq\Lambda^{6}\omega.$

Then we have that $T_{n}\leq C_{\gamma}(n^{2}/\Lambda)^{\frac{7+\gamma}{\gamma-1}},$ so  the induction relation on $\delta_{n}$ gives
$\delta_{1}\leq C\delta \quad \textmd{\textmd{and}} \quad\delta_{n+1}=C(\frac{n^{2}}{\Lambda})^{9}e^{C(n^{2}/\Lambda)^{(14+2\gamma)/(\gamma-1)}}\delta^{2}_{n}.$

To make this relation hold, we also assumed that $\delta_{n}$ is bounded below by $C_{F}\zeta_{n},$ the error coming from the short-time iteration; but this
follows easily by construction, since the constraints imposed on $\delta_{n}$ are much worse than those on $\zeta_{n}.$

\begin{center}
%\item\section{Linear}
{\bf\large 1. \quad Appendix }
\end{center}
$Proof$ $of$  $Theorem$  $7.4.$  Here we only prove (i), the  proof of (ii) is similar. We decompose the proof into three step.

$Step$ 1. Crude pointwise bounds. From (7.9), we have
$$\varphi(t)=\sum_{k\in\mathbb{Z}^{3}_{\ast}}|\Phi(k,t|e^{2\pi(\lambda t+\mu)|k|}
\leq A+\sum_{k\in\mathbb{Z}^{3}_{\ast}}\int^{t}_{0}|K^{0}(k,t-s)|e^{2\pi(\lambda t+\mu)|k|}|\Phi(t,s)|d\tau$$
$$+\int^{t}_{0}(K_{0}(t,s)+K_{1}(t,s)+\frac{c_{0}}{(1+s)^{m}})\varphi(s)ds$$
$$\leq A+\int^{t}_{0}(K_{0}(t,s)+K_{1}(t,s)+\frac{c_{0}}{(1+s)^{m}}+\sup_{k\in\mathbb{Z}^{3}_{\ast}}|K^{0}(k,t-s)|e^{2\pi\lambda(t-s)|k|})
\varphi(s)d\tau.$$
We note that for any $k\in\mathbb{Z}^{3}_{\ast}$ and $t\geq0,$
$$|K^{0}(k,t-s)|e^{2\pi\lambda|k|(t-s)}\leq 4\pi^{2}|\widehat{W}(k)|C_{0}e^{-2\pi(\lambda_{0}-\lambda)|k|t}|k|^{2}t
\leq\frac{CC_{0}C_{W}}{\lambda_{0}-\lambda},$$
where (here and below) $C$ stands for a numeric constant which may change from line to line. Assuming that $\int^{t}_{0}K_{0}(t,\tau)d\tau\leq\frac{1}{2},$
we deduce that
$$\varphi(t)\leq A+\frac{1}{2}\sup_{0\leq s\leq t}\varphi(s)+C\int^{t}_{0}\bigg(\frac{C_{0}C_{W}}{\lambda_{0}-\lambda}+c(1+s)
+\frac{c_{0}}{(1+s)^{m}}\bigg)\varphi(s)ds,$$
and, by Gr$\ddot{\textmd{o}}$nwall's lemma,
\begin{align}
\varphi(t)\leq 2Ae^{C(C_{0}C_{W}t/(\lambda_{0}-\lambda)+c(t+t^{2})+c_{0}C_{m})},
\end{align}
where $C_{m}=\int^{\infty}_{0}(1+s)^{-m}d\tau.$

$Step$ 2. $L^{2}$ bound. For all $k\in\mathbb{Z}^{3}_{\ast}$ and $t\geq0,$ we define
$\Psi_{k}(t)=e^{-\varepsilon t}\Phi(k,t)e^{2\pi(\lambda t+\mu)|k|}, \mathcal{K}^{0}_{k}(t)=e^{-\varepsilon t}K^{0}(k,t)$
$e^{2\pi(\lambda t+\mu)|k|},$
$R_{k}(t)$
$=e^{-\varepsilon t}\bigg(\Phi(k,t)-\int^{t}_{0}K^{0}(k,t-s)\Phi(k,s)ds\bigg)e^{2\pi(\lambda t+\mu)|k|}=(\Psi_{k}-\Psi_{k}\ast
\mathcal{K}^{0}_{k})(t),$ and we extend all these functions by $0$ for negative values of $t.$ Taking Fourier transform in the time-variable yields
$\hat{R}_{k}=(1-\widehat{\mathcal{K}}^{0}_{k})\widehat{\Psi_{k}}.$ Since the Stability condition implies that $|1-\widehat{\mathcal{K}}^{0}_{k}|\geq\kappa,$
 we can deduce that $\|\hat{\Psi}_{k}\|_{L^{2}}\leq\kappa^{-1}\|\hat{R}_{k}\|_{L^{2}},$ i.e., $\|\Psi_{k}\|_{L^{2}}\leq\kappa^{-1}\|R_{k}\|_{L^{2}}.$ So we have
 \begin{align}\|\Psi_{k}-R_{k}\|_{L^{2}(dt)}\leq\kappa^{-1}\|\mathcal{K}^{0}_{k}\|_{L^{1}(dt)}\|R_{k}\|_{L^{2}(dt)} \quad \textmd{for}\quad \textmd{all}
 \quad k\in\mathbb{Z}_{\ast}^{3}.
 \end{align}

 Then \begin{align}
& \|\varphi(t)e^{-\varepsilon t}\|_{L^{2}(dt)}=\|\sum_{k\in\mathbb{Z}^{3}}|\Psi_{k}|\|_{L^{2}(dt)}\leq\|\sum_{k\in\mathbb{Z}^{3}}|R_{k}|\|_{L^{2}(dt)}
 +\sum_{k\in\mathbb{Z}^{3}}\|R_{k}-\Psi_{k}\|_{L^{2}(dt)}\notag\\
 &\leq\|\sum_{k\in\mathbb{Z}^{3}}|R_{k}|\|_{L^{2}(dt)}(1+\frac{1}{\kappa})
 \end{align}

 Next, we note that
 $$\|\mathcal{K}^{0}_{k}\|_{L^{1}(dt)}\leq 4\pi^{2}|\widehat{W}(k)|\int^{\infty}_{0}C_{0}e^{-2\pi(\lambda_{0}-\lambda)|k|t}|k|^{2}t dt
 \leq 4\pi|\widehat{W}(k)|\frac{C_{0}}{(\lambda_{0}-\lambda)^{2}},$$
 so $\sum_{ k\in\mathbb{Z}_{\ast}^{3}}\|\mathcal{K}^{0}_{k}\|_{L^{1}(dt)}\leq 4\pi(\sum_{ k\in\mathbb{Z}_{\ast}^{3}}|\widehat{W}(k)|)\frac{C_{0}}
 {(\lambda_{0}-\lambda)^{2}}.$
 Furthermore, we get
 \begin{align}
 &\|\varphi(t)e^{-\varepsilon t}\|_{L^{2}(dt)}\leq\bigg(1+\frac{CC_{0}C_{W}}{\kappa(\lambda_{0}-\lambda)^{2}}\bigg)
 \|\sum_{k\in\mathbb{Z}_{\ast}^{3}}\|_{L^{2}(dt)}\notag\\
 &\leq \bigg(1+\frac{CC_{0}C_{W}}{\kappa(\lambda_{0}-\lambda)^{2}}\bigg)\bigg(\int^{\infty}_{0}e^{-2\varepsilon t}
 \bigg(A+\int^{t}_{0}\bigg(K_{1}+K_{0}+\frac{c_{0}}{(1+s)^{m}}\bigg)\varphi(s)ds\bigg)^{2}\bigg)^{\frac{1}{2}}.
 \end{align}

 By Minkowski's inequality, we separate (7.15) into various contributions which we estimate separately. First,
 $\bigg(\int^{\infty}_{0}e^{-2\varepsilon t}A^{2}dt\bigg)^{\frac{1}{2}}=\frac{A}{\sqrt{2\varepsilon}}.$ Next, for any $T\geq1,$ by Step 1 and
 $\int^{t}_{0}K_{1}(t,\tau)ds\leq\frac{Cc(1+t)}{\alpha},$ we have
 \begin{align}
 &\bigg(\int^{T}_{0}e^{-2\varepsilon t}\bigg(\int^{t}_{0}K_{1}(t,s)\varphi(s)\bigg)^{2}\bigg)^{\frac{1}{2}}
 \leq (\sup_{0\leq t\leq T}\varphi(t)) \bigg(\int^{T}_{0}e^{-2\varepsilon t}\bigg(\int^{t}_{0}K_{1}(t,s)\bigg)^{2}\bigg)^{\frac{1}{2}}\notag\\
 &\leq CAe^{C(C_{0}C_{W}T/(\lambda_{0}-\lambda)+c(T+T^{2}))}\frac{c}{\alpha}\bigg(\int^{\infty}_{0}e^{-2\varepsilon t}(1+t)^{2}dt\bigg)^{\frac{1}{2}}\notag\\
 &\leq CA\frac{c}{a\varepsilon^{\frac{3}{2}}}e^{C(C_{0}C_{W}T/(\lambda_{0}-\lambda)+c(T+T^{2}))}.
 \end{align}

 Invoking Jensen's inequality and Fubini's theorem, we also have
 \begin{align}
 &\int^{\infty}_{T}e^{-2\varepsilon t}\bigg(\int^{t}_{0}K_{1}(t,s)\varphi(s)ds\bigg)^{2}dt\bigg)^{\frac{1}{2}}
 =\int^{\infty}_{T}\bigg(\int^{t}_{0}K_{1}(t,s)e^{-2\varepsilon (t-s)}e^{-2\varepsilon s}\varphi(s)ds\bigg)^{2}dt\bigg)^{\frac{1}{2}}\notag\\
 &\leq\int^{\infty}_{T}\bigg(\int^{t}_{0}K_{1}(t,s)e^{-\varepsilon (t-s)}ds\bigg)\bigg(
\int^{\infty}_{T}\bigg(\int^{t}_{0}K_{1}(t,s)e^{-\varepsilon (t-s)}e^{-2\varepsilon s}\varphi(s)^{2}ds\bigg)dt\bigg)^{\frac{1}{2}}\notag\\\notag\\
&\leq\bigg(\sup_{t\geq T}\int^{t}_{0}K_{1}(t,s)e^{-\varepsilon (t-s)}ds\bigg)^{\frac{1}{2}}\bigg(
\int^{\infty}_{T}\bigg(\int^{t}_{0}K_{1}(t,s)e^{-\varepsilon (t-s)}e^{-2\varepsilon s}\varphi(s)^{2}ds\bigg)dt\bigg)^{\frac{1}{2}}\notag\\\notag\\
&=\bigg(\sup_{t\geq T}\int^{t}_{0}K_{1}(t,s)e^{-\varepsilon (t-s)}ds\bigg)^{\frac{1}{2}}\bigg(
\int^{\infty}_{0}\int^{\infty}_{\max\{s,T\}}K_{1}(t,\tau)e^{-\varepsilon (t-\tau)}e^{-2\varepsilon s}\varphi(s)^{2}dtds \bigg)^{\frac{1}{2}}\notag\\
&\leq\bigg(\sup_{t\geq T}\int^{t}_{0}K_{1}(t,s)e^{-\varepsilon (t-s)}d\tau\bigg)^{\frac{1}{2}}\bigg(
\sup_{\tau\geq0}\int^{\infty}_{s}K_{1}(t,s)e^{-\varepsilon (t-s)}e^{-2\varepsilon s}dt\bigg)^{\frac{1}{2}}
\bigg(\int^{\infty}_{0}\varphi(s)^{2}e^{-2\varepsilon s}ds \bigg)^{\frac{1}{2}}.\notag\\
\end{align}
Similarly,
\begin{align}
 &\int^{\infty}_{T}e^{-2\varepsilon t}\bigg(\int^{t}_{0}K_{0}(t,s)\varphi(s)ds\bigg)^{2}dt\bigg)^{\frac{1}{2}}
 \leq\bigg(\sup_{t\geq T}\int^{t}_{0}K_{0}(t,s)ds\bigg)^{\frac{1}{2}}\bigg(
\sup_{s\geq0}\int^{\infty}_{\tau}K_{0}(t,s)dt\bigg)^{\frac{1}{2}}
\bigg(\int^{\infty}_{0}\varphi(s)^{2}ds \bigg)^{\frac{1}{2}}
\end{align}

The last term is also split, this time according to $\tau\leq T$ or $\tau> T:$
\begin{align}&\bigg(\int^{\infty}_{0}e^{-2\varepsilon t}\bigg(\int^{T}_{0}\frac{c_{0}\varphi(s)}{(1+s)^{m}}ds\bigg)^{2}dt\bigg)^{\frac{1}{2}}
\leq c_{0}(\sup_{0\leq s\leq T}\varphi(s))
\bigg(\int^{\infty}_{0}e^{-2\varepsilon t}\bigg(\int^{T}_{0}\frac{ds}{(1+s)^{m}}\bigg)^{2}dt\bigg)^{\frac{1}{2}}\notag\\
 &\leq c_{0}\frac{CA}{\sqrt{\varepsilon}}e^{C(C_{0}C_{W}T/(\lambda_{0}-\lambda)+c(T+T^{2}))}C_{m},
\end{align}
and \begin{align}
&\bigg(\int^{\infty}_{0}e^{-2\varepsilon t}\bigg(\int^{t}_{T}\frac{c_{0}\varphi(s)}{(1+s)^{m}}ds\bigg)^{2}dt\bigg)^{\frac{1}{2}}
\leq c_{0}\bigg(\int^{\infty}_{0}e^{-2\varepsilon t}\varphi(t)^{2}\bigg)^{\frac{1}{2}}
\bigg(\int^{\infty}_{0}\int^{t}_{T}\frac{e^{-2\varepsilon (t-s)}}{(1+s)^{2m}}ds dt\bigg)^{\frac{1}{2}}\notag\\
&=c_{0}\bigg(\int^{\infty}_{0}e^{-2\varepsilon t}\varphi(t)^{2}\bigg)^{\frac{1}{2}}
\bigg(\int^{\infty}_{T}\frac{ds}{(1+s)^{2m}}\bigg)^{\frac{1}{2}}
\bigg(\int^{\infty}_{0}e^{-2\varepsilon s} ds\bigg)^{\frac{1}{2}}
=\frac{C_{2m}^{\frac{1}{2}}c_{0}}{T^{m-\frac{1}{2}}\sqrt{\varepsilon}}\bigg(\int^{\infty}_{0}e^{-2\varepsilon t}\varphi(t)^{2}\bigg)^{\frac{1}{2}}.
\end{align}
Gathering estimates (7.16)-(7.20), we deduce from (7.15) that
\begin{align}&\|\varphi(t)e^{-\varepsilon t}\|_{L^{2}(dt)}\leq\bigg(1+\frac{CC_{0}C_{W}}{\kappa(\lambda_{0}-\lambda)^{2}}\bigg)\frac{CA}{\sqrt{\varepsilon}}
\bigg(1+\frac{c}{a\varepsilon}+c_{0}C_{m}\bigg)e^{C(C_{0}C_{W}T/(\lambda_{0}-\lambda)+c(T+T^{2}))}\notag\\
&\quad \quad \quad \quad\quad\quad\quad+a\|\varphi(t)e^{-\varepsilon t}\|_{L^{2}(dt)},
\end{align}
where $$a=\bigg(1+\frac{CC_{0}C_{W}}{\kappa(\lambda_{0}-\lambda)^{2}}\bigg)\bigg[
\bigg(\sup_{t\geq T}\int^{t}_{0}e^{-\varepsilon t}K_{1}(t,s)e^{\varepsilon s}ds\bigg)^{\frac{1}{2}}
\bigg(\sup_{s\geq 0}\int^{\infty}_{s}e^{\varepsilon s}K_{1}(t,s)e^{-\varepsilon t}dt\bigg)^{\frac{1}{2}}$$
$$+\bigg(\sup_{t\geq T}\int^{t}_{0}K_{0}(t,\tau)d\tau\bigg)^{\frac{1}{2}}\bigg(
\sup_{s\geq0}\int^{\infty}_{\tau}K_{0}(t,s)dt\bigg)^{\frac{1}{2}}+\frac{C_{2m}^{\frac{1}{2}}c_{0}}{T^{m-\frac{1}{2}}\sqrt{\varepsilon}}\bigg].$$

Using Proposition 7.2 and 7.3, together with the assumptions of Theorem 7.4, we see that $a\leq\frac{1}{2}$ for $\chi$ sufficiently small. Then we have
$$\|\varphi(t)e^{-\varepsilon t}\|_{L^{2}(dt)}\leq\bigg(1+\frac{CC_{0}C_{W}}{\kappa(\lambda_{0}-\lambda)^{2}}\bigg)\frac{CA}{\sqrt{\varepsilon}}
\bigg(1+\frac{c}{a\varepsilon}+c_{0}C_{m}\bigg)e^{C(C_{0}C_{W}T/(\lambda_{0}-\lambda)+c(T+T^{2}))}.$$

$Step$ 3. For $t\geq T,$ using (7.9) we get
\begin{align}& e^{-\varepsilon t}\varphi(t)\leq Ae^{-\varepsilon t}+\bigg[\bigg(\int^{t}_{0}\bigg(\sup_{k\in\mathbb{Z}^{3}_{\ast}}|K^{0}(k,t-s)|
e^{2\pi\lambda(t-s)|k|}\bigg)^{2}ds\bigg)^{\frac{1}{2}}\notag\\
&+\bigg(\int^{t}_{0}K_{0}(t,\tau)^{2}d\tau\bigg)^{\frac{1}{2}}+\bigg(\int^{\infty}_{0}\frac{c^{2}_{0}}{(1+s)^{2m}}ds\bigg)^{\frac{1}{2}}
+\bigg(\int^{t}_{0}e^{-2\varepsilon t}K_{1}(t,\tau)^{2}e^{2\varepsilon \tau}d\tau\bigg)^{\frac{1}{2}}\bigg]
\bigg(\int^{\infty}_{0}\varphi(s)e^{-\varepsilon s}ds\bigg)^{\frac{1}{2}}.
\end{align}

We note that, for any $k\in\mathbb{Z}^{3}_{\ast},$ $(|K^{0}(k,t)|e^{2\pi\lambda|k|t})^{2}\leq C\pi^{4}|\widehat{W}(k)|^{2}|\hat{f}^{0}(kt)|^{2}|k|^{4}t^{2}
\leq \frac{CC_{0}}{(\lambda_{0}-\lambda)^{2}}C^{2}_{W}e^{-2\pi(\lambda_{0}-\lambda)t},$ so we get $\int^{t}_{0}\bigg(\sup_{k\in\mathbb{Z}_{\ast}^{3}}
|K^{0}(k,t-\tau)|e^{2\pi\lambda(t-\tau)|k|}\bigg)^{2}d\tau\leq\frac{CC^{2}_{0}C^{2}_{W}}{(\lambda_{0}-\lambda)^{3}}.$

From Proposition 7.2,(7.22), the conditions of Theorem 7.4 and Step 2, the conclusion is finished.

Having fixed $\Lambda,$ we will check that for $\delta$ small enough, the above relation hold and the fast convergence of $\{\delta_{i}\}^{\infty}_{i=1}.$
 The details are similar to that of the local-time case,and it can be also found in [21], here we omit it.
 $$$$

 Acknowledgements: The author is grateful for the comfortable and superior academic environment of Yau Mathematical Science Center, Tsinghua University.  I also thank professor Yifei Wu in Tianjin University for providing the financial support and hospitality.

\end{document}